\documentclass{amsart}
\usepackage{amssymb,graphicx}

\usepackage{color}

\newcommand{\R}{{\mathbb R}}

\theoremstyle{remark}

\newtheorem{remark}{Remark}

\theoremstyle{definition}

\theoremstyle{theorem}

\newtheorem{theorem}{Theorem}

\newtheorem*{merci}{Acknowledgements}
\newtheorem{conj}{Conjecture}

\begin{document}
\title[Dispersive perturbations of Burgers'  equation]{A numerical 
approach to Blow-up issues for dispersive perturbations of Burgers' equation}
\author{Christian Klein}
\address{Institut de Math\'ematiques de Bourgogne\\
                Universit\'e de Bourgogne, 9 avenue Alain Savary, 21078 Dijon
                Cedex, France\\
		Tel.: +33380395858, Fax.: +33380395869\\
    E-mail Christian.Klein@u-bourgogne.fr}

\author{Jean-Claude Saut}
\address{Laboratoire de Math\' ematiques, UMR 8628\\
Universit\' e Paris-Sud et CNRS\\ 91405 Orsay, France\\ E-mail jean-claude.saut@math.u-psud.fr}
\date{January 15th  2014}
\maketitle
\large

\begin{abstract}
We provide a detailed numerical study of various issues pertaining to the 
dynamics of the Burgers equation perturbed by a weak dispersive term: 
blow-up in finite time versus global existence, nature of the 
blow-up, existence for ``long" times, and the decomposition of the 
initial data into solitary waves plus radiation. We 
numerically construct solitons for fractionary Korteweg-de Vries 
equations.
\end{abstract}

\section{Introduction}
This paper is concerned with the numerical study of the dynamics and 
blow-up issues for ``weak" dispersive perturbations of Burgers' (inviscid) equation.
As in \cite{LPS2} the motivation is to study the influence of dispersion on  the dynamics of solutions to the Cauchy problem for \lq\lq weak\rq\rq \ dispersive perturbations of hyperbolic quasilinear equations or systems, as for instance the  Boussinesq  systems for surface water waves.

We thus want  to investigate the competition  between nonlinearity 
and dispersion. Usually this problem is attacked by fixing the 
dispersion, {\it eg} that of the Korteweg- de Vries (KdV) equation 
and varying the nonlinearity $u^p u_x$ in the  generalized KdV (gKdV)
equation

\begin{equation}\label{gKdv}
u_t+u^pu_x+u_{xxx}=0,
\end{equation}

or the generalized Benjamin-Ono
equation (gBO)
\begin{equation}\label{gBO}
u_t+u^pu_x-Hu_{xx}=0,
\end{equation}
where $p\geq 2$ and $H$ is the Hilbert transform. The case $p=2$ is the $L^2$ critical case for the gBO equation.
Numerical simulations in \cite{BK} are performed for $p=3$ suggesting   finite time blow-up in this case. No rigorous proof of blowup is known for \eqref
{gBO} when $p\geq 2$ though, contrary to the $L^2$ critical gKdV equation (that is \eqref {gKdv} with $p=4)$ for which a proof of existence of blowing-up solutions has been given  in \cite{MM} and revisited in \cite{MMR}.

However it is  probably more physically relevant  to fix the 
quadratic nonlinearity ({\it eg} $uu_x$) and to vary (lower) the 
dispersion. In fact in many problems arising from Physics or 
Continuum Mechanics, the nonlinearity is quadratic, with terms like 
$(u\cdot \nabla) u$ and the dispersion is in some sense weak. On the 
other hand, dispersive terms occuring in the Schr\"{o}dinger or in the KdV 
equations are obtained by a Taylor expansion of the true dispersion at a 
given wave number $k_0$ or  in the long wave limit. This allows 
strong dispersive effects (and nice mathematical properties!), but 
restricts the range of frequencies for which the model is relevant. 
For instance the dispersion relation of the KdV equation is a bad 
approximation of the one of the water waves for not too small frequencies. We refer to    \cite{La} for examples of water wave models with ``full dispersion", that is models which share the same dispersion as that of the water waves system. In particular the original dispersion is not strong enough for yielding the dispersive estimates that allows to solve the Cauchy problem in relatively large functional classes (like the KdV or Benjamin-Ono equation in particular), down to the energy level for instance.\footnote{And thus obtaining {\it global well-posedness} from the conservation laws.}

In fact many nonlinear dispersive systems have the following structure

\begin{equation}\label{basic}
\partial_t U+\mathcal B U+\epsilon\; \mathcal A(U, \nabla U)+\epsilon \mathcal LU=0,
\end{equation}

where the order $0$ part $\partial_t U+\mathcal B U$ is linear hyperbolic, $ \mathcal L$ being  a linear (not necessarily skew-adjoint) dispersive operator and $\epsilon>0$ is a small parameter which measures the (comparable) nonlinear and dispersive effects. Both the linear part and the dispersive part may involve nonlocal terms (see {\it eg} \cite{X}, \cite{SX2}).

Boussinesq systems for surface water waves  (see \cite{BCL, BCS1}) 
are  important examples of somewhat similar  systems. Note however 
that the Boussinesq systems \eqref{Bsq} below cannot be reduced exactly to 
the form \eqref {basic} except when $b=c=0.$ Otherwise the presence 
of a ``BBM like" term induces a smoothing effect on one or both 
nonlinear terms. They read

\begin{equation} \label{Bsq}
\left\{ \begin{array}{l}   \partial_t\eta+\text{div}\,
\textbf{v}+\epsilon\;\text{div} \, (\eta\textbf{v})+\epsilon( a\;\text{div}\Delta \textbf{v}-b\Delta \eta_t)=0 \\
\partial_t \textbf{v} +\nabla \eta +\epsilon\frac12 \nabla(|\textbf{v}|^2)+\epsilon(c\nabla \Delta
\eta-d\Delta {\bf v}_t)=0 \end{array} \right., \quad (x_1,x_2) \in \mathbb R^2, \ t \in
\mathbb R.
\end{equation}
where $a,b,c,d$ are modelling constants satisfying the constraint 
$a+b+c+d=\frac{1}{3}$ and ad hoc conditions implying the 
well-posedness of the linearized system at the trivial solution $(0,{\bf 0}).$ 

When  $b>0,d>0$, the dispersion in \eqref{Bsq} is ``weak" (the 
corresponding linear operator is of order  $-1,0$ or $1$ (see 
\cite{BCS1}) contrary to the case $b=d=0,$ $a<0, c<0$ when it is of order $3$ as in  the KdV equation).

We will not address the blow-up issues for \eqref{Bsq}, which are 
completely open, but instead focus on  one dimensional scalar equations.

More precisely, the major question addressed here will be whether or 
not smooth solutions of the Cauchy problem for weak dispersive 
perturbations of Burgers' equation  develop singularities in finite 
time and if blow-up occurs, what is its nature: a shock like in 
solutions to Burgers' equation, a  blow-up similar to the one of the 
$L^2$ critical or supercritical generalized KdV equations, or an ``energy critical or super critical" blow-up, a case that never occurs for the generalized KdV equations (see below for a more detailed discussion).

We will also consider the {\it long time existence} issues, that is 
what is the qualitative behavior of the solution when it is global, 
or before the blow-up time if not; in particular does it show a 
decomposition into solitary waves plus radiation, {a typical behavior of solutions to integrable equations such as the KdV equation? Also when global 
existence is not assured, how does the presence of a ``weak" 
dispersive term affect the life span of the solution to the 
underlying Burgers equation when a small parameter $\epsilon$ appears in front of the quadratic term  and possibly in front of the dispersive terms?

 The last question is particularly important for the complete rigorous  
 justification of water waves models (see \cite{La}). It is worth 
 noticing that except when the solution is global (which can be proven essentially only for scalar one-way models) the only long time existence results for water waves models such as the Boussinesq systems \eqref{Bsq} are established on the ``hyperbolic" time scale $1/\epsilon$ (see \cite{MSZ, X, SX, SX2}). We will not address this issue in the context of relevant water waves models but again for the fractionary KdV or BBM (fKdV or fBBM) equations which will serve as toy models.
 
The paper is organized as follows: in section 2 we collect some 
analytically known facts about the equations to be studied numerically. 
The above questions will  then be addressed for these equations 
numerically in section 3. We add some concluding remarks in section 4. 

\subsection*{Notations}

The following notations will be used throughout this article. The Fourier transform of a function $f$  will be denoted $\hat{f}$ or $\mathcal {F}f.$ For any $s\in \R,$ we define $D^s f$ by $\widehat{D^sf}(\xi)=|\xi|^s\hat f (\xi).$

For $1 \le p \le \infty$, $L^p(\mathbb R)$ is the usual Lebesgue space with the norm $\|\cdot \|_{L^p}$, and for $s \in \mathbb R$, the Sobolev space $H^s(\mathbb R)$ is defined via its usual norm $\|f \|_{H^s}= (\int_\R(f^2+|D^s f|^2)dx)^{1/2}$.

\section{Theoretical preliminaries}
In this section we gather known facts about various dispersive 
regularizations of Burgers' equation. In particular we consider 
fractionary KdV and BBM equations. We formulate analytic questions
which will then be addressed numerically in the following section.

\subsection{Fractionary KdV equations}

 We will thus focus as a paradigm on   one -dimensional  model equations. 
 A first one  (introduced by Whitham \cite{W}) is of KdV type with a weak dispersion.
\begin{equation}\label{Whit}
u_t+uu_x+\int_{-\infty}^{\infty}k(x-y)u_x(y,t)dy=0.
\end{equation}
This equation can also be written in the form
\begin{equation}\label{Whibis}
u_t+uu_x-Lu_x=0,
\end{equation}
where the Fourier multiplier operator $L$ is defined by 
$$\widehat{Lf}(\xi)=p(\xi)\hat{f}(\xi),$$
with $p=\hat{k}.$

 In the original Whitham equation, the kernel $k$ was given by 
\begin{equation}\label{tanh}
k(x)=\frac{1}{2\pi}\int_\R \left( \frac{\tanh \xi}{\xi} \right)^{1/2} e^{ix\xi} d\xi,
\end{equation}
that is $p(\xi)=\left( \frac{\tanh \xi}{\xi} \right)^{1/2},$ which 
corresponds to the phase velocity of purely gravitational waves. 
The Whitham equation is also the one-dimensional version of the Full dispersion Kadomtsev-Petviashvili equation (FDKP) studied in \cite{LaSa}.
When surface tension is added, the above $p$ has to be changed to 
$p_S(\xi)=(1+\beta|\xi|^2)^{1/2}\left( \frac{\tanh \xi}{\xi} \right)^{1/2},$ where $\beta\geq 0$ measures the surface tension effects. This leads to the extended Whitham equation.


We will consider here a family of equations where 
$p(\xi)=|\xi|^\alpha$, that we call ``fractionary KdV equations". The 
case $\alpha =2$ corresponds to the usual KdV equation, $\alpha =1$ to the Benjamin-Ono equation. It is well known (see \cite{HIKK}) that for $\alpha\geq 1$ the solutions of the Cauchy problem (in appropriate functional spaces) are global and therefore no finite time blow-up occurs. In the KdV case, it is known that the solution decomposes into solitons traveling to the right and radiation going to the left. 

The case $\alpha <1$ is more delicate, and we will focus on it. 

When $-1<\alpha <0,$ a blow-up occurs, that is there is a finite time 
blow-up  of the solution of the Cauchy problem corresponding to 
suitable  smooth initial data (see \cite{CCG} and\cite{ CE, VMH, NS} 
for related equations). The proof in \cite{CCG} extends easily to the 
Whitham equation (see \cite{LaSa}).  By  contradiction one proves that the $C^{1+\alpha}$ norm of a solution blows up in finite time, but the proof does not indicate whether the gradient only or both the solution and its gradient blow up. Our numerical simulations (see Section 3.5) suggest that only the gradient blows up.

The occurrence of a possible blow-up when $0<\alpha <1$ is an open problem. It is claimed in \cite{KZ} without proof that then a shock formation is not possible.

The situation appears anyhow very different from that of the {\it fractal Burgers equation} that is the Burgers equation perturbed by a fractionary dissipation term $D^\alpha u.$ There (see \cite {ADV, KNS}), a possible shock formation persists when $0<\alpha <1/4$.}

We briefly recall here the known results for the associated Cauchy problem, that is

\begin{equation}\label{Cauchy}
u_t+uu_x-D^\alpha u_x=0,\quad u(.,0)=u_0,
\end{equation}

where $\widehat{D^\alpha f}(\xi)=|\xi|^\alpha \hat f(\xi).$

The following quantities are formally conserved by the flow associated to \eqref{Cauchy},
\begin{equation} \label{M}
M(u)=\int_{\mathbb R}u^2(x,t)dx,
\end{equation}
and the Hamiltonian
\begin{equation} \label{H}
H(u)=\int_{\mathbb R}\big( \frac{1}{2} |D^{\frac{\alpha}2}u(x,t)|^2-\frac{1}{6}u^3(x,t)\big) dx.
\end{equation}

Note that by the Sobolev embedding $H^{\frac{1}{6}}(\R)\hookrightarrow L^3(\R)$,  $H(u)$ is well-defined when $\alpha \geq \frac{1}{3} $ but it does not make sense  when $u\in H^{\alpha /2}(\R),\alpha <\frac{1}{3},$ the {\it energy super critical case}.

Moreover, equation \eqref{Cauchy} is invariant under the scaling transformation 
\begin{equation} \label{scaling}
u_{\lambda} (x,t)=\lambda^{\alpha}u(\lambda x,\lambda^{\alpha+1}t),
\end{equation}
for any positive number $\lambda$. A straightforward computation shows that $\|u_{\lambda}\|_{\dot{H}^s}=\lambda^{s+\alpha-\frac{1}{2}}\|u_{\lambda}\|_{\dot{H}^s}$, and thus the critical index corresponding to \eqref{Cauchy} is $s_{\alpha}=\frac{1}{2}-\alpha$. In particular, equation \eqref{Cauchy} is $L^2$-critical for $\alpha=\frac{1}{2}$. This corresponds to $p=4$ in the generalized KdV equation \eqref{gKdv}.


Compactness arguments  (see \cite{S}) prove that \eqref {Cauchy} 
admit global weak solutions (without uniqueness) in $L^{\infty}(\R; H^{\alpha /2}(\R))$ when $\alpha >\frac{1}{2}$ for initial data in $H^{\alpha /2}(\R)$(with a smallness condition when $ \alpha =\frac{1}{2})$ and, thanks to a Kato type smoothing effect global weak solutions in $L^{\infty}(\R; L^2(\R))\cap L^2_{\text{loc}}(\R;H_{\text{loc}}^{\alpha/2}(\R)),$ see \cite{GV, GV2}.

By using standard energy methods and the fact that the dispersive 
term is skew-adjoint (thus without using explicitely the dispersion), 
one can on the other hand prove that the Cauchy problem associated to \eqref{Cauchy} is locally well-posed in $H^s(\mathbb R)$ for $s>\frac{3}{2},$ which is the same ``hyperbolic type " result as for the Burgers equation.

Taking into account the dispersion, it has been established in \cite{LPS2} that the Cauchy problem for \eqref {Cauchy} is locally well posed for initial data  in $H^{s}(\R),$ for $s>s_\alpha=\frac{3}{2}-\frac{3\alpha}{2}>\frac{\alpha}{2},$ which does not allow to globalize the solution using the conservation laws.

It is well known that the solitary waves play a significant role in 
the dynamics of the KdV equation and one can ask whether or not this is the case for the dispersive Burgers equation.

A  (localized) solitary wave solution of \eqref{Cauchy} of the form
$u(x,t)=Q_c(x-ct)$  must satisfy the equation
\begin{equation} \label{solitarywave}
D^{\alpha}Q_c + cQ_c -\frac12Q_c^2=0,
\end{equation}
where $c>0$.

One does not expect solitary waves to exist when
$\alpha< \frac13$ since then the Hamiltonian
does not make sense (see a formal argument in \cite{KZ} and a 
rigorous proof in \cite{LPS2} which also proves nonexistence of solitary waves when $\alpha <0$). 

 On the other hand solitary waves 
exist when $\alpha >\frac{1}{3}$ (see  \cite{EGW, FQT, F, FL}) and 
they have a slow decay for $|x|\to\infty$ as 
$\frac{1}{x^{1+\alpha}}.$ They are expected to be orbitally stable when $\alpha >\frac{1}{2},$ that is in the $L^2$ subcritical case (this would correspond to $p<4$ for  the generalized KdV equations \eqref {gKdv}). It is also worth noticing that solitary waves exist for the original Whitham equation (\cite{EGW}). The proof uses in a crucial way that the dispersion relation of the Whitham equation approaches that of the KdV equation for small frequencies. Existence of {\it periodic} traveling waves to the Whitham equation has been proven in \cite{EK} and their stability properties studied in \cite{HJ} and numerically in \cite{CKKD}.

Note that the exponent $\alpha=\frac{1}{3}$ corresponds to the 
so-called ``energy critical case" that never occurs for the generalized KdV equations where the energy (Hamiltonian) always makes sense in the energy space $H^1(\R).$

At this stage, one could make the following conjectures for \eqref {Cauchy} when $0<\alpha <1$:

\vspace{0.5cm}

{\bf 1.} No hyperbolic blow-up exists (no blow-up of the spatial gradient with bounded sup-norm).

{\bf 2.} When $\frac{1}{2}<\alpha <1,$ the solution is global (no blow-up).

{\bf 3.} When $\frac{1}{3}<\alpha \leq \frac{1}{2},$ which is the 
$L^2$ supercritical case, one has a ``nonlinear dispersive blow-up", 
with a kind of self-similar structure. This would correspond to 
$p\geq 4$ for the generalized KdV equation. Recall that this type of 
blow-up for \eqref {gKdv} is supported by numerical simulations in 
the super critical case $p>4$ (\cite {BDKM, KP2013}, and in the 
critical case  $p=4$ \cite{KP2013} and rigorously proven in the 
critical case $p=4$ (\cite {MM, MMR}). The dynamics of the blow-up is 
different in the $L_{2}$-critical and in $L_{2}$-supercritical cases.

{\bf 4.} When $0<\alpha<\frac{1}{3},$  the energy critical 
case  which has no counterpart for the generalized KdV equations, one expects a blow-up, but with a different structure 
of that of the previous case.

\subsection{Fractionary BBM equations}

We now turn to the  BBM version of the dispersive Burgers equation, namely

\begin{equation}\label{fracBBM}
\partial_tu+\partial_xu+u\partial_xu+D^{\alpha}\partial_tu=0,
\end{equation}
where the operator $D^{\alpha}$ is defined as previously.

The case $\alpha=2$ corresponds to the classical BBM equation, $\alpha=1$ to the BBM version of the Benjamin-Ono equation. A class of equations containing \eqref{fracBBM} has been introduced in the Appendix 1 of  \cite{BBM}.

For any $\alpha$ the energy
\begin{equation}
    E(t)=\int_{\R}(u^2+|D^{\frac{\alpha}{2}}u|^2) dx
    \label{fBBMenergy}
\end{equation}
is formally conserved. By a standard compactness method, this implies that the Cauchy problem for \eqref{fracBBM} admits a global weak solution in $L^{\infty}(\R;H^{\frac{\alpha}{2}}(\R))$ for any initial data $u_0=u(\cdot,0)$ in $H^{\frac{\alpha}{2}}(\R).$

One can also use the equivalent form

\begin{equation}\label{fracBBM2}
\partial_t u+\partial_x(I+D^{\alpha})^{-1}\left(u+\frac{u^2}{2}\right)=0,
\end{equation}

which gives the Hamiltonian formulation

$$u_t+J_{\alpha}\nabla_u H(u)=0$$

where the skew-adjoint operator $J_{\alpha}$ is given by $J_{\alpha}=\partial_x(I+D^{\alpha})^{-1}$ and $H(u)=\frac{1}{2}\int_\R (u^2+\frac{1}{3}u^3).$ Note that the Hamiltonian makes  sense and is formally conserved for $u\in H^{\frac{\alpha}{3}}(\R)$ if and only if $\alpha\geq \frac{1}{3}.$

 
\vspace{0.3cm}
We will again focus on the case  $0<\alpha <1.$ Actually when $\alpha 
\geq1,$ \eqref{fracBBM2} is an ODE in the Sobolev space $H^s(\R)$, 
$s>\frac{1}{2},$ and one obtains by standard arguments (see \cite {Ma, ABS2}) the local well-posedness of the Cauchy problem in $H^s(\R)$, $s>\frac{1}{2}.$ When $\alpha =1$ (the Benjamin-Ono BBM equation), the conservation of energy and an ODE argument as in \cite{S}
 or the Br\' ezis-Gallou\"{e}t inequality implies that this local solution is in fact global.
 
The following local well-posedness result is established in \cite{LPS2}.

  \begin{theorem}
 Let $0<\alpha<1.$
 Then the Cauchy problem for \eqref {fracBBM} or \eqref {fracBBM2} is locally well-posed for initial data in $H^r(\R),\; r>r_{\alpha}= \frac{3}{2}-\alpha .$
 \end{theorem}


It is clear from the formulation \eqref{fracBBM2} that the fractionary BBM equation is for $0<\alpha <1$ a different dispersive perturbation of the Burgers equation of that given by the fractionary KdV equation: roughly speaking, one replaces the $\partial_ x$ derivative of $u+\frac{u^2}{2}$ by a derivative of order $1-\alpha.$


The question is now:
Is there a blow-up when $0<\alpha< 1$ and if yes, what is its nature?

\subsection{Large time existence issues}

An important issue, when the solution is global, is 
to describe its qualitative behavior. In the range of fKdV equations, the only fully understood case is $\alpha=2,$ the KdV equation where Inverse Scattering techniques yield a full description of the solution. For the other values of $\alpha,$ such a description is not known (the only other integrable equation in the family corresponds to $\alpha =1,$ the Benjamin-Ono equation, but then Inverse Scattering techniques only work for small initial data). One aim of our numerical simulations is to investigate what happens when $\frac{1}{2}<\alpha<1.$ They actually suggest that a kind of {\it decomposition into solitary waves plus radiation} occurs,
despite the fact that none of those equations are integrable. A similar behavior seems  to hold for the fBBM equations, at least when $\alpha >\frac{1}{3}.$

Another interesting and widely open question is to investigate the influence of a weak dispersion on the existence time of solutions to hyperbolic equations or systems (see \cite{SX, SX2, X} for some partial answers in the case of various water wave systems, in particular the Boussinesq systems).

Again a simple relevant example is the dispersive Burgers equation

\begin{equation}\label{Cauchybis}
u_t+\epsilon uu_x-\epsilon D^\alpha u_x=0,\quad u(.,0)=u_0,
\end{equation}

  where $\epsilon$ is a small positive parameter (in the water waves context it could model the comparable effects of dispersion and nonlinearity).
  
  Since the dispersive term in \eqref{Cauchybis} is skew-adjoint, one immediately deduces that for initial data $u_0$ of order $1$ in say $H^s(\R),s>\frac{3}{2},$  the solution exists on time scales of order $\frac{1}{\epsilon}.$

By changing the time variable as $\tau=\epsilon t$ one eliminates the $\epsilon's$ from \eqref{Cauchybis}, that is one obtains

\begin{equation}\label{tau}
u_\tau+uu_x- D^\alpha u_x=0,\quad u(.,0)=u_0,
\end{equation}

and we have the following dichotomy:

$\bullet$ Either the solution to \eqref{tau} is global and so is that of \eqref{Cauchybis},

$\bullet$ or the solution to \eqref{tau} has a lifespan of order 
$O(1)$ and that of \eqref{Cauchybis} has a lifespan of order $O(1/\epsilon),$ that is the same as for the Burgers equation.

Our simulations below suggest that when $\alpha >1/2$ the 
first situation occurs while when $\alpha\leq 1/2$ the second one holds.


 Things are more delicate when the small parameter appears only in front of the nonlinear term  as shows the  striking example of the Burgers-Hilbert equation
  
  \begin{equation}\label{BH}
  u_t+\epsilon uu_x+\mathcal H u=0,\quad u(\cdot, t)=u_0
  \end{equation}
  
  where $\mathcal H$ is the Hilbert transform. It is established in \cite{HI, HITW} that for initial data of order $O(1)$ the solution of \eqref{BH} exists on time scales of order $\frac{1}{\epsilon^2}$ while the corresponding Burgers solution exists on time scales of order $\frac{1}{\epsilon}.$

One should thus consider the Cauchy problem

\begin{equation}\label{ter}
u_t+\epsilon uu_x- D^\alpha u_x=0,\quad u(.,0)=u_0,
\end{equation}

  An equivalent formulation (by setting $v=\epsilon u$) is to consider the Cauchy problem
  
  \begin{equation}\label{terbis}
  v_t+vv_x- D^\alpha v_x=0,\quad v(\cdot,0)=\epsilon u_0,
  \end{equation}
    
    and the question is to see how the existence time 
    $O(\frac{1}{\epsilon})$ is enlarged by the dispersive term $ 
    D^\alpha v_x, $ in particular, is it possible to prove the existence of global small solutions to \eqref{terbis}? This would give an example of initial data leading to a shock for the Burgers equation and to a global solution for fKdV.
    
    Of course similar questions can be addressed for the BBM version, that is
    
    \begin{equation}\label{BBMbis}
u_t+u_x+\epsilon uu_x+\epsilon D^\alpha u_t=0,\quad u(.,0)=u_0,
\end{equation}

or 

 \begin{equation}\label{BBMter}
u_t+u_x+\epsilon uu_x+ D^\alpha u_t=0,\quad u(.,0)=u_0,
\end{equation}

Note that in \eqref {BBMbis} one cannot eliminate anymore the $\epsilon's$ by a change of time scale because of the transport term.

  Similar issues for  the Boussinesq systems \eqref{Bsq} are widely 
  open. Here (and for other relevant water waves models)  no 
  existence results beyond the ``hyperbolic" time 
  $\frac{1}{\epsilon}$ seem to be known (see \cite{MSZ,X,SX, SX2}). Recall that for the Boussinesq systems, the error estimates with the full water waves system being $O(\epsilon ^2 t),$ (see \cite {BCL}), and the solutions are expected to exist at least up to time scales of order $\frac{1}{\epsilon ^2}.$ 
There are here difficulties that are not present in the toy model fKdV or fBBM. First the change of time scale
$\tau =\epsilon t$ does not eliminate anymore the $\epsilon's$ in the system because of the order zero hyperbolic part. Moreover, as previously noticed, when $b$ and/or $d$ are strictly positive, the nonlinear term in the corresponding equation is smoothed by the operator $(I-b\Delta)^{-1},$ (resp. $(I-d\Delta)^{-1}).$

  The case studied in \cite{A, Sc} (that is $a=b=c=0, d=\frac{1}{3}$) is specially interesting since it was considered in those references as a perturbation of the (hyperbolic) Saint-Venant system.

\section{Numerical study of weakly dispersive regularizations of 
Burgers' equation}
In this section we will study numerically solutions to fKdV, fBBM  
and Whitham
equations. The focus will be on smooth, localized initial data for 
which we study the decomposition into solitons and a possible blow-up 
in finite time. We also address the long time behavior of small 
initial data. First we will present the used numerical tools and 
discuss how accuracy is ensured. Then we will study concrete examples.

\subsection{Numerical Methods for the time evolution}

Since the equations to be studied in this paper contain fractionary 
derivatives being defined as Fourier multipliers, it is convenient to 
use Fourier spectral methods in the numerical treatment. We 
concentrate  on rapidly decreasing, smooth initial data which can 
be analytically continued within the finite numerical precision as 
periodic functions. In other words, we use large enough 
computational domains that the Fourier coefficients of the 
periodically continued initial data decrease to machine precision 
($10^{-16}$ in our case). The fractional derivatives for Fourier 
series are defined in analogy to the definition for Fourier 
transforms. The discrete Fourier transform is computed via a 
\emph{fast Fourier transform} (fft).

The fKdV equation for the Fourier transform $\hat{u}$ of $u$ has the 
form
\begin{equation}
    \hat{u}_{t}=\mathcal{F}(u)+\mathcal{L}\hat{u},
    \label{gKdVfourier}
\end{equation}
where $\mathcal{L}=i\xi|\xi|^{\alpha}$ is a linear 
operator\footnote{For the Whitham equation (\ref{Whit}), the symbol 
$|\xi|^{\alpha}$ is simply replaced with $\sqrt{\tanh(\xi)/\xi}$.}, and where 
$\mathcal{F}(u)=-i\xi \widehat{u^{2}}/2$ 
denotes the nonlinear terms. It is an advantage of Fourier methods 
that the $x$-derivatives and thus the operator $\mathcal{L}$ are 
diagonal.  For equations of 
the form (\ref{gKdVfourier}) with diagonal $\mathcal{L}$, there are many efficient high-order time 
integrators, see 
e.g.~\cite{KassamTrefethen2005,HO,Klein2008,KleinRoidot2011} and 
references therein. This allows for an efficient integration in time.

In \cite{KP2013} it was shown that an implicit Runge-Kutta method of 
fourth order (IRK4), a two-stage Gauss scheme, is both very efficient 
and accurate up to blow-up for generalized KdV equations. 
For the initial value problem
$y'=f(y,t)$, $y(t_0)=y_0$ and constant time steps $t_{m}$, $m=0,1,\ldots$ with $t_{m+1}-t_{m}=h$ 
and $y(t_{m}) =y_{m}$, this scheme has the form
\begin{eqnarray}
 y_{m+1} = y_{m} + \frac{h}{2}(K_{1}+K_{2}), \\
 K_{i} = f\left(t_{m}+c_ {i}h,\,y_{m}+h  
 \underset{j=1}{\overset{s}{\sum}} \, a_{ij}K_{j}\right), \quad i=1,2,
\end{eqnarray}
where 
$c_{1}=\frac{1}{2}-\frac{\sqrt{3}}{6}$, 
$c_{2}=\frac{1}{2}+\frac{\sqrt{3}}{6}$ and $a_{11}=a_{22}=1/4$,
$a_{12}=\frac{1}{4}-\frac{\sqrt{3}}{6}$, 
$a_{21}=\frac{1}{4}+\frac{\sqrt{3}}{6}$. 
The implicit equations for $K_{1}$ and $K_{2}$ are solved iteratively with a 
simplified Newton scheme, 
\begin{equation}
    \begin{pmatrix}
        K_{1} \\
        K_{2}
    \end{pmatrix}=
    \begin{pmatrix}
        1-ha_{11}\mathcal{L} &-h a_{12}\mathcal{L} \\
        -ha_{21}\mathcal{L} & 1-ha_{22}\mathcal{L}
    \end{pmatrix}^{-1}
    \begin{pmatrix}
        \mathcal{F}(\hat{u}_{m}+h(a_{11}K_{1}+a_{12}K_{2}))+\mathcal{L}\hat{u}_{m} \\
        \mathcal{F}(\hat{u}_{m}+h(a_{21}K_{1}+a_{22}K_{2}))+\mathcal{L}\hat{u}_{m}
    \end{pmatrix}
    \label{Ki}.
\end{equation}
Since  the operator $\mathcal{L}$ is 
diagonal, the inverse matrix on the right-hand side of (\ref{Ki}) 
can be given explicitly. In this form the iteration converges rapidly 
(at early times in 3 to 4 iterations). 

The fBBM equation 
(\ref{fracBBM}) reads in Fourier space
$$\hat{u}_{t}=-\frac{i\xi}{1+|\xi|^{\alpha}}(\hat{u}+\widehat{u^{2}}/2),$$
which is again integrated with the IRK4 method.

Accuracy of the numerical solution is controlled as discussed in 
\cite{Klein2008,KleinRoidot2011,KP2013} via the numerically computed 
energies (\ref{H}) and (\ref{fBBMenergy}) 
which will depend on time due to unavoidable numerical errors. We use 
the quantity
\begin{equation}
    \Delta=|E(t)/E(0)-1|
    \label{Delta}
\end{equation}
as an indicator of the numerical 
accuracy. It was shown in \cite{Klein2008,KleinRoidot2011} that the numerical accuracy 
of this quantity overestimates the $L_{\infty}$ norm of the 
difference between numerical and exact solution by two to three 
orders of magnitude. A precondition for the usability of this 
quantity is a sufficient resolution in Fourier space.

\paragraph{\textbf{Dynamic rescaling}}
To study blow-up in the solutions to fKdV, the scale invariance 
(\ref{scaling}) can be used as for gKdV and NLS equations in the form 
of a dynamic rescaling, see for instance \cite{SS99,KP2013} and 
references therein. To this end the constant $\lambda$ in 
(\ref{scaling}) is replaced with a $t$-dependent function $L(t)$
\begin{equation}
    y = \frac{x-x_{m}}{L},\quad 
    \frac{d\tau}{dt}=\frac{1}{L^{1+\alpha}},\quad U = L^{\alpha}u
    \label{resc}.
\end{equation}
As in \cite{KP2013} for gKdV, we take into account the fact that the 
peak at $x=x_{m}(t)$ developing eventually into a blow-up is travelling with 
increasing speed. To address this, we choose a commoving frame. 
Equation (\ref{resc}) implies for (\ref{Cauchy})
\begin{equation}
    U_{\tau}-(\ln L)_{\tau}\left(\alpha U+y U_{y}\right)
    -\frac{x_{m,\tau}}{L}U_{y}
    +UU_{y}-D^{\alpha}_{y}U_{y}=0
    \label{sautalphar}.
\end{equation}
The scaling function $L$ can be chosen to keep certain norms constant, for 
instance the $L_{\infty}$-norm. It appears that the choice to keep 
the $L_{2}$ norm of $U_{y}$ constant is numerically preferred for 
stability reasons. This choice leads to 
\begin{equation}
    a:=(\ln 
    L)_{\tau}=\frac{1}{(2\alpha+1)||U_{y}||_{2}^{2}}\int_{\mathbb{R}}^{}U^{2} U_{yyy}dy
    \label{L2x}.
\end{equation}
Since $y=0$ is supposed to be an extremum of the solution, we 
have $U_{y}(0,\tau)=0$ and $U_{y}(0,\tau)=0$.
Putting $U(0,0,\tau)=U_{0}=const$, we get from (\ref{sautalphar}) for the speed 
$v=x_{m,\tau}/L$
\begin{equation}
    v = U_{0}+\frac{1}{U_{yy}(0,\tau)}(D^{\alpha}U_y)(0,\tau).
        \label{vx}
\end{equation}
Thus all quantities in (\ref{sautalphar}) can be expressed in terms of $U$ 
alone. 
In the case of an $L_{\infty}$ blow-up for $t=t^{*}$, the scaling function $L$ is 
expected to vanish for $\tau\to\infty$. Thus for $t\to t^{*}$, both 
$x$ and $t$ are dynamically rescaled. Note that for solutions of
this equation the energy
\begin{equation}
    E[U]=  
    \frac{1}{L^{1-3\alpha}}\int_{\R}(U^2+|D^{\frac{\alpha}{2}}U|^2) dy
    \label{Eresc}
\end{equation}
is conserved. It can be seen that the case $\alpha=1/3$ is energy 
critical in the sense that it is invariant under the rescaling 
(\ref{sautalphar}). 

For $\tau\to\infty$, it is expected that the solution $U$ to the 
rescaled equation (\ref{sautalphar}) as well as $a\to a_{\infty}$ and 
$v\to v_{\infty}$ become 
$\tau$ independent (this behavior was proven in \cite{MMR} for the 
$L_{2}$ critical gKdV equation for initial data in the vicinity of a 
soliton). In this case equation (\ref{sautalphar}) reduces to an 
ordinary fractionary differential equation 
\begin{equation}
    -a_{\infty}\left(\alpha U^{\infty}+y U^{\infty}_{y}\right)
    -v_{\infty}U^{\infty}_{y}
    +U^{\infty}U^{\infty}_{y}-D^{\alpha}_{y}U^{\infty}_{y}=0
    \label{sautalpharinfty}.
\end{equation}
If $a_{\infty}$ vanishes, this is exactly the equation 
(\ref{solitarywave}) for the soliton. A blow-up of this type would be 
asymptotically given by a rescaled solution of equation 
(\ref{sautalpharinfty}).

In \cite{KP2013} it was shown numerically that generic rapidly decreasing hump-like 
initial data for gKdV lead to a 
tail of dispersive oscillations towards infinity with slowly 
decreasing amplitude. Due to the imposed periodicity, these 
oscillations reappear after some time
on the opposing side of the computational domain 
and lead to numerical instabilities in the dynamically rescaled 
equation. The source of these problems is the term $yU_{y}$ in 
(\ref{sautalphar}) since $y$ is large at the boundaries of the 
computational domain. This implies that this term is very sensitive 
to numerical errors.  For gKdV this could be addressed 
by using high resolution in time and large computational domains. It 
turns out that for fKdV, the dispersive oscillations have an 
amplitude that decreases even more slowly towards infinity. Thus we 
could not compute long enough to get conclusive results. Instead we 
integrate fKdV as described above directly, and then we use some 
postprocessing to characterize the type of blow-up via the above 
rescaling. We read off the 
time evolution of the quantity $L$ as defined via the $L_{2}$ norm of 
the gradient, 
 \begin{equation}
    \frac{||U_{y}||_{2}}{||u_{x}||_{2}}=L^{\alpha+1/2}
    \label{scal},
\end{equation}
where the constant $||U_{y}||_{2}$ is chosen to be 
$||u_{x}(x,0)||_{2}$.  This allows to study the type of the 
blow-up for fKdV in a similar way as for gKdV in \cite{KP2013} 
without actually solving equation (\ref{sautalphar}). We ensure that 
the numerically computed energy of the solution
is conserved to a 
certain precision ($\Delta <10^{-3}$ in (\ref{Delta})), that there is
sufficient resolution in Fourier space 
and in time for all computed times close to the blow-up time $t^{*}$. 
We generally choose the time step in blow-up scenarios such that the 
accuracy is limited by the resolution in Fourier space, i.e., that a 
further reduction of the time step for a given number of Fourier 
modes does not change the final result within numerical accuracy. 

For gKdV, essentially two cases were relevant for the 
$\tau$-dependence of the scaling factor $L$, and we want to test 
whether the same is true for fKdV. In the $L_{2}$-critical 
case, one has $L\propto \tau^{-\gamma}$ with $\gamma=1$. An algebraic decrease 
of $L$ close to blow-up implies with (\ref{resc}) for fKdV
\begin{equation}
    ||u_{x}||_{2}^{2}\propto (t^{*}-t)^{-\frac{2\alpha+1}{1+\alpha-1/\gamma}},\quad
    ||u||_{\infty}\propto 
    (t^{*}-t)^{-\frac{\alpha}{1+\alpha-1/\gamma}}.
    \label{L2scal}
\end{equation}
In the supercritical case for gKdV ($n>4$), there does not exist a proof yet, 
but it is expected that $L$ vanishes exponentially with $\tau$ for 
$\tau\to\infty$, 
$L\propto \exp(-\kappa \tau)$ with $\kappa$ a positive constant. This 
implies with (\ref{resc})
\begin{equation}
    ||u_{x}||_{2}^{2}\propto (t^{*}-t)^{-\frac{2\alpha+1}{1+\alpha}},\quad
    ||u||_{\infty}\propto 
    (t^{*}-t)^{-\frac{\alpha}{1+\alpha}}.    
    \label{genscal}
\end{equation}
It will be tested whether such scalings can be observed in the 
numerical experiments for fKdV.

\paragraph{\textbf{Singularity tracing in the complex plane}}
An important question in the context of the numerical solution of 
nonlinear dispersive PDEs is the identification and the 
characterization of a blow-up of the solution. The task is to 
identify the appearence of a singularity of the solution. Since we 
use a Fourier approach 
here, a possible tool in this context is 
a method from asymptotic 
Fourier analysis first applied to 
numerically identify singularities in solutions to PDEs in 
\cite{SSF}. The idea is to use the fact that a singularity of a real 
function in the complex plane of the form $U\sim (z-z_{0})^{\mu}$ 
($\mu$ not an integer)
leads for $|\xi|$ large to a Fourier transform of the form (if this 
is the only singularity of this type in the complex plane)
\begin{equation}
    |\hat{U}|\sim 
    \frac{1}{\xi^{\mu+1}} e^{-\xi\delta},\quad \xi\ll 1,
    \label{fourasymp}
\end{equation}
where $\delta=\Im z_{0}$. In \cite{KR2013} it was  discussed how this 
approach can be used to quantitatively identify the time where the 
singularity hits the real axis, i.e., where the real solution becomes 
singular. In addition this method gives the quantity $\mu$ which 
characterizes the type of the singularity. It was shown in 
\cite{KR2013} that the quantity $\delta$ can be identified reliably 
from a fitting of the Fourier coefficients, whereas there is a larger 
inaccuracy in the quantity $\mu$. Thus the 
numerically determined values for $\mu$ have to be taken with a grain 
of salt.

\paragraph{\textbf{Tests}}
Contrary to the case of the KdV equation, it is not known if for the BO and the BBM equation the solutions of the initial value problems for these 
equations decompose for large times into solitons and 
radiation. The reason is that the Inverse Scattering Transform 
technique for the BO equation works only for small initial data and 
that the BBM equation is not completely integrable. Nevertheless a 
similar behavior can be expected for fKdV and fBBM solution,  and it 
is tempting to examine this issue. Solitons are very important in this context, and 
an accurate reproduction of these solutions for the cases where they 
are explicitly known is crucial. Equation 
(\ref{solitarywave}) has for $\alpha=1$ the solution
\begin{equation}
    Q_{c}=\frac{4c}{1+(cx)^{2}}
    \label{Qc}.
\end{equation}
This solution decreases as $x^{-2}$ for $|x|\to\infty$ which is 
numerically problematic for the Fourier methods we use here. The 
algebraic decrease implies that the solution cannot be analytically 
continued within the given finite precision even on large 
domains. Thus a Gibbs phenomenon will appear which has the 
consequence that the Fourier coefficients do not decrease to machine 
precision. In other words the numerical precision will be limited 
because of this algebraic decrease. This is in contrast to gKdV equations, for 
which the solitons are explicitly known, and where they decrease 
exponentially for $|x|\to\infty$. We recall that for $\alpha<1$ the soliton 
solutions of (\ref{solitarywave}) decrease as $|x|^{-(1+\alpha)}$ for 
$|x|\to\infty$, see \cite{FL}. 

To numerically test the propagation of the soliton, we use $Q_{c}$ in 
(\ref{Qc}) as initial data for fKdV with $\alpha=1$ and $c=2$. The 
test is performed with $N=2^{14}$ Fourier modes for 
$x\in100[-\pi,\pi]$. The modulus of the Fourier coefficients 
decreases in this case to $10^{-8}$. We propagate the initial data 
with $N_{t}=10^{4}$ time steps to time $t=1$. The numerically 
computed energy is conserved to the order of $10^{-14}$, i.e., 
machine precision. But since 
the Fourier coefficients just decrease to the order of $10^{-8}$, the 
accuracy of the solution cannot be higher than this value. In fact we find that the 
$L_{\infty}$ norm of the difference between exact and numerical 
solution is of the order of $10^{-7}$ as indicated by the resolution 
in Fourier space. This shows that the code reproduces the exact solution 
with the precision available in Fourier space, and that it is 
crucial to provide enough resolution there to distinguish the formation of 
solitons from the appearance of a blow-up. 

Travelling wave solutions of the fBBM equation satisfy again equation 
(\ref{solitarywave}) after a change $c\to \tilde{c}=1-1/c$ and $u\to 
u/c$. Thus for $\alpha=1$, the fBBM equations have soliton solutions 
of the form
\begin{equation}
    u(x,t)=\frac{4(c-1)}{1+\tilde{c}^{2}(x-ct)^{2}},\quad 
    \tilde{c}=1-\frac{1}{c}.
    \label{QfBBM}
\end{equation}
If we choose the same $c$ and the same parameters for the numerical 
solution as above, the modulus of the Fourier coefficients decreases 
only to $10^{-6}$ in this case since the initial data decrease even 
more slowly for $|x|\to\infty$ because of $\tilde{c}=1/2$. The 
numerically computed energy is again conserved to the order of 
machine precision. But since the resolution in Fourier space is only 
of the order of $10^{-6}$, the difference between the exact and the 
numerical solution for $t=1$ is of the order of $10^{-5}$ as 
expected.  

\subsection{Numerical construction of solitons}
In this subsection we numerically solve equation (\ref{solitarywave}) to 
obtain solitons of the studied equations. Note that we can 
concentrate on the case $c=1$ since the solution for general $c$ 
follows from the former via the simple rescaling 
\begin{equation}
    Q_{c}(z)=cQ_{1}(zc^{1/\alpha})
    \label{solresc}.
\end{equation}
To simplify the notation we suppress the index of $Q$ for $c=1$.

We solve equation (\ref{solitarywave}) for $c=1$  as in the 
previous subsection in Fourier space, where it has the form
\begin{equation}
    (|\xi|^{\alpha}+1)\hat{Q}-\frac{1}{2}\widehat{Q^{2}}=0
    \label{solfourier}.
\end{equation}
In the numerical solution the Fourier transform is approximated as 
before via a discrete Fourier series computed via an fft. The task is 
thus to find the (nontrivial) zero of the function $F= 
(|\xi|^{\alpha}+1)\hat{Q}-\frac{1}{2}\widehat{Q^{2}}$ depending on 
the vector $\hat{Q}$ of the discrete Fourier transform of $Q$, which 
can in general only be done iteratively. A 
potential problem in this context is that equation (\ref{solfourier}) 
has the trivial solution $\hat{Q}=0$, and a straight-forward fixed 
point iteration typically converges to this solution, even if one 
starts with an exact solution. Thus we use a Newton method in the 
standard form 
$\hat{Q}_{n+1}=\hat{Q}_{n}-\mbox{Jac}^{-1}F(\hat{Q}_{n})$, where 
$\mbox{Jac}$ is the Jacobian of (\ref{solfourier}). 

A technical problem in this context is the slow decrease of the 
soliton solution for $|x|\to\infty$ which is known to be of the order 
$O(|x|^{-(1+\alpha)})$, see \cite{FL}. Since the function $Q$ is treated as 
essentially periodic, it will not be smooth at the boundaries of the 
computational domain. This loss of smoothness implies an algebraic 
decrease of the modulus of the Fourier coefficients and thus a slower 
convergence rate of the numerical scheme. To address this we choose a 
large domain, $x\in 100[-\pi,\pi]$ and high resolution in Fourier 
space, $N>2^{14}$. Since the Jacobian entering the Newton iteration 
is a $N\times N$ matrix, its inverse cannot be computed with 
conventional methods on a standard computer. Therefore we use 
Krylov-subspace techniques, here GMRES \cite{gmres}. The idea is to 
compute the inverse of a matrix iteratively. The advantage of GMRES 
is that the to be inverted matrix does not have to be known 
explicitly, just its action on a vector. Thus to compute 
$\mbox{Jac}^{-1}F(\hat{Q}_{n})$, we need only 
$$\mbox{Jac} \hat{X}=(|\xi|^{\alpha}+1)\hat{X}  - \widehat{QX},$$
where $X$ is some vector.  The Newton iteration is stopped when the 
$L_{\infty}$ norm of $F$ is smaller than $10^{-8}$. Generally we 
reach machine precision after 4-8 iterations if the initial iterate 
is chosen as explained below. 

\begin{figure}[htb!]
  \includegraphics[width=0.49\textwidth]{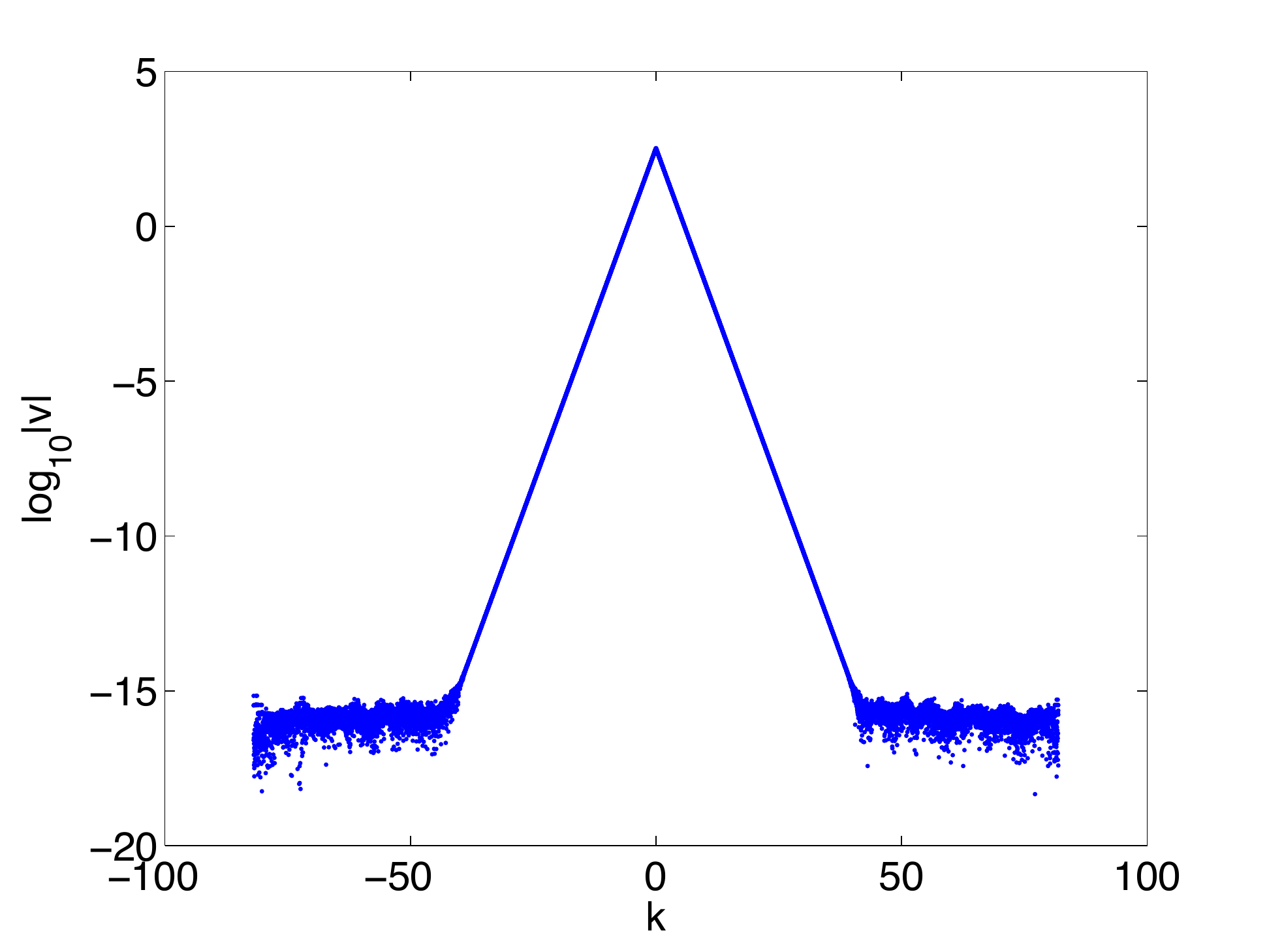}
  \includegraphics[width=0.49\textwidth]{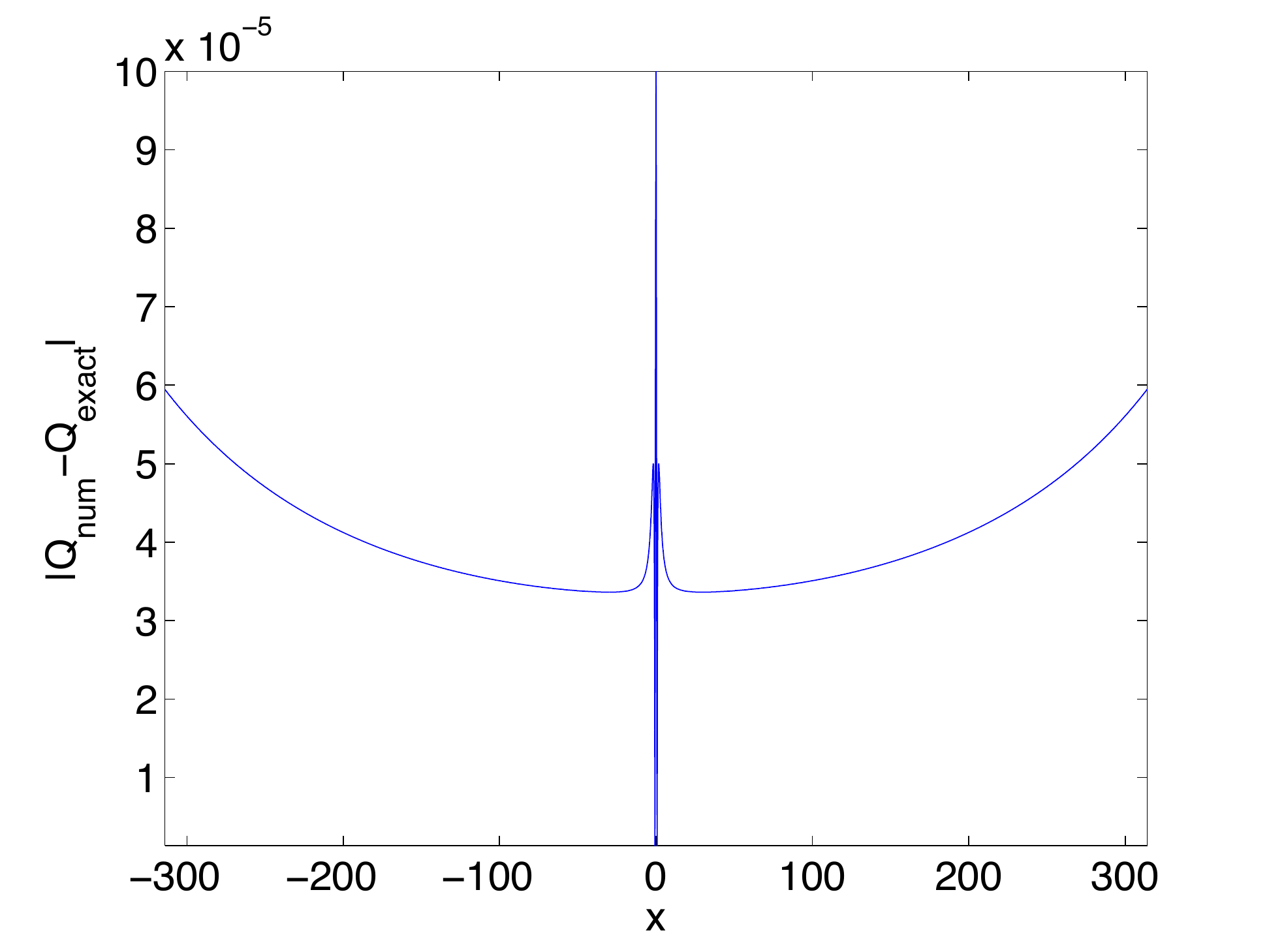}
 \caption{Modulus of the Fourier coefficients for the soliton 
 (\ref{solitarywave}) of the BO equation ($\alpha=1$) with $c=1$ on the left, and 
 the difference of the numerical and the exact solution on the right.}
 \label{BOsol}
\end{figure}

To test this approach we consider the BO soliton (\ref{Qc}), i.e., 
the case $\alpha=1$. We choose $N=2^{14}$ and $1.1Q$ as the initial 
iterate, i.e., the exact solution multiplied with a factor $1.1$. After 3 
iterations we reach a residual $||F||_{\infty}$ of the order of 
$10^{-12}$. It can be seen in Fig.~\ref{BOsol} that the Fourier 
coefficients of the solution decrease to machine precision. In the 
same figure we also show the difference between numerical and exact 
solution which is of the order of $10^{-4}$, the largest difference 
being observed at the maximum of the solution. The reason for this 
difference despite resolution in Fourier space and the solution of 
equation (\ref{solfourier}) to machine precision is the periodic 
setting we are using here. This also implies the increase of the 
difference towards the boundary of the domain. But the example indicates that we should 
be able to obtain the solitons of (\ref{solitarywave}) to the order 
of plotting accuracy also for $\alpha<1$.

\begin{figure}[htb!]
  \includegraphics[width=\textwidth]{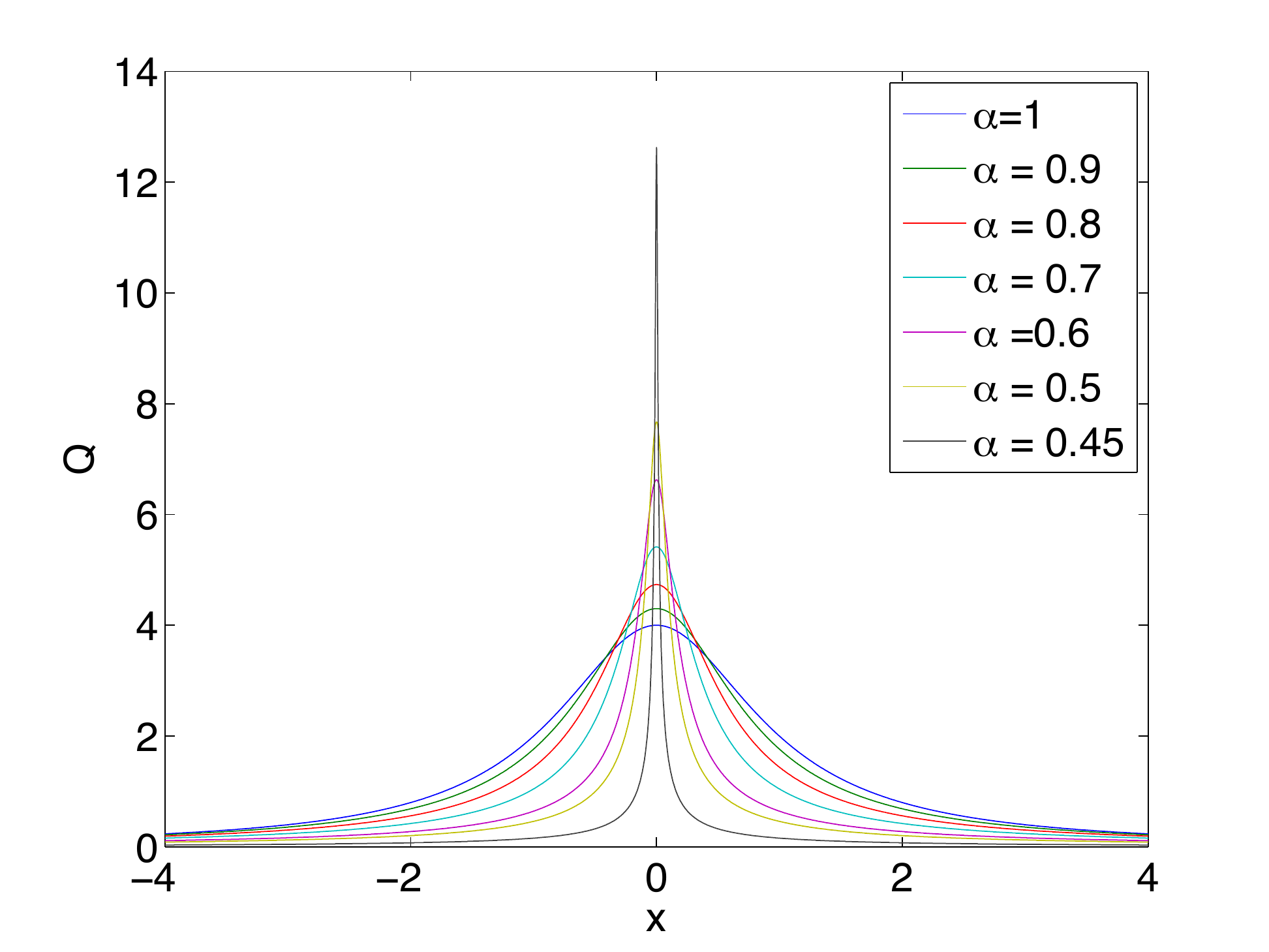}
 \caption{Solitary waves (\ref{solitarywave}) for $c=1$ and different 
 values of $\alpha$.}
 \label{solQleg}
\end{figure}

To construct solitons for $\alpha<1$, we use the solution for a 
given value of $\alpha$ as the initial iterate for equation 
(\ref{solfourier}) for a 
smaller value $\alpha$: to construct the solution for $\alpha=0.9$ 
for instance, we start with the BO soliton. By lowering $\alpha$ by 
$0.1$ in this way, we can reach values of $0.6$ without problems. For 
even smaller values of $\alpha$, we have to take smaller steps in 
$\alpha$, of the order of $0.01$. The resulting solutions are shown 
in Fig.~\ref{solQleg} for values of $x$ close to the maximum. 

\begin{figure}[htb!]
  \includegraphics[width=0.49\textwidth]{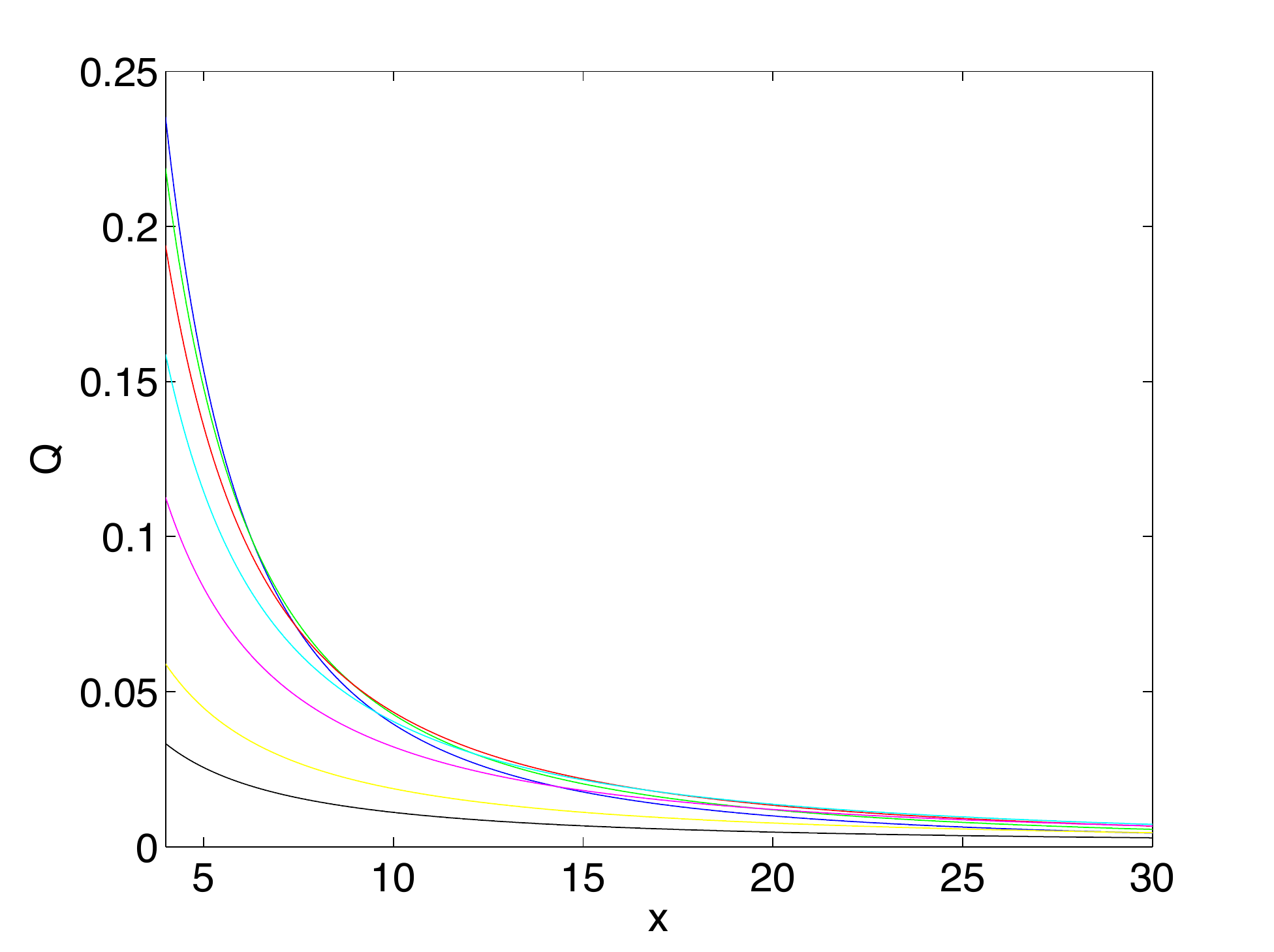}
  \includegraphics[width=0.49\textwidth]{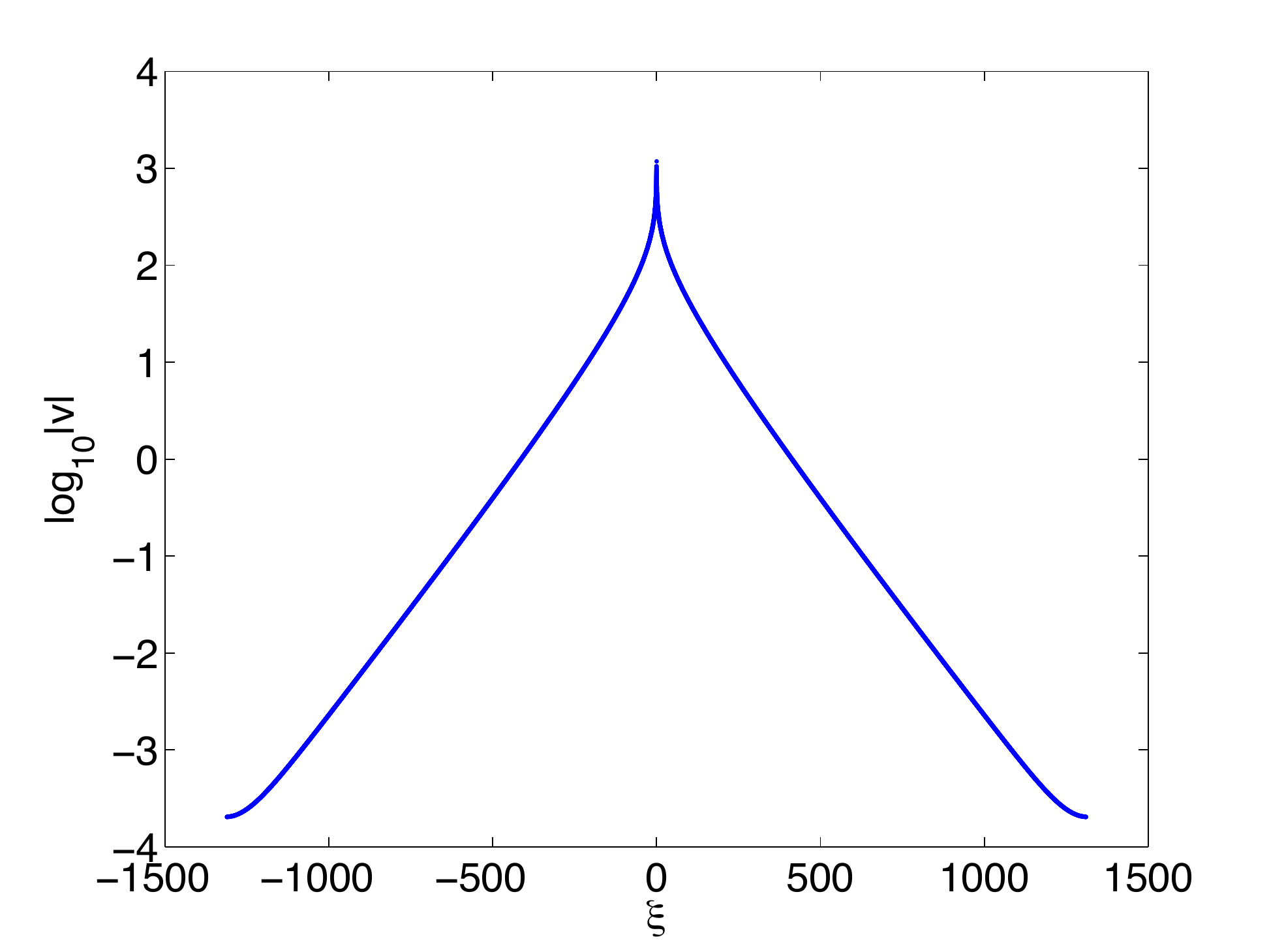}
 \caption{Solitons of Fig.~\ref{solQleg} for larger values of $x$ on 
 the left, 
 and the modulus of the Fourier coefficients for the soliton 
 (\ref{solitarywave}) with $c=1$ and $\alpha=0.45$ on the right.}
 \label{solmax}
\end{figure}

It can 
be seen that the solutions become more and more peaked the smaller 
$\alpha$ is, and the decrease towards infinity becomes slower; the 
latter can be seen more clearly in Fig.~\ref{solmax}, where the 
BO soliton intersects for larger $x$ the solitons for smaller $\alpha$. This 
is challenging for the numerical method, and we use $N=2^{18}$ 
Fourier modes to provide the necessary resolution. In 
Fig.~\ref{solmax}, it can be also seen that the modulus of the Fourier 
coefficients still goes down to $10^{-4}$ in this case. To be able to 
increase the resolution further, it would be necessary at one point 
to use higher precision, i.e., quadruple instead of the used double 
precision. But given the strong increase of the maximum of the 
solution for decreasing 
$\alpha$, it seems doubtful that one can reach much smaller values of 
$\alpha$ with this approach. The limit $\alpha=1/3$ appears to be 
inaccessible in this way. 

\begin{figure}[htb!]
  \includegraphics[width=0.49\textwidth]{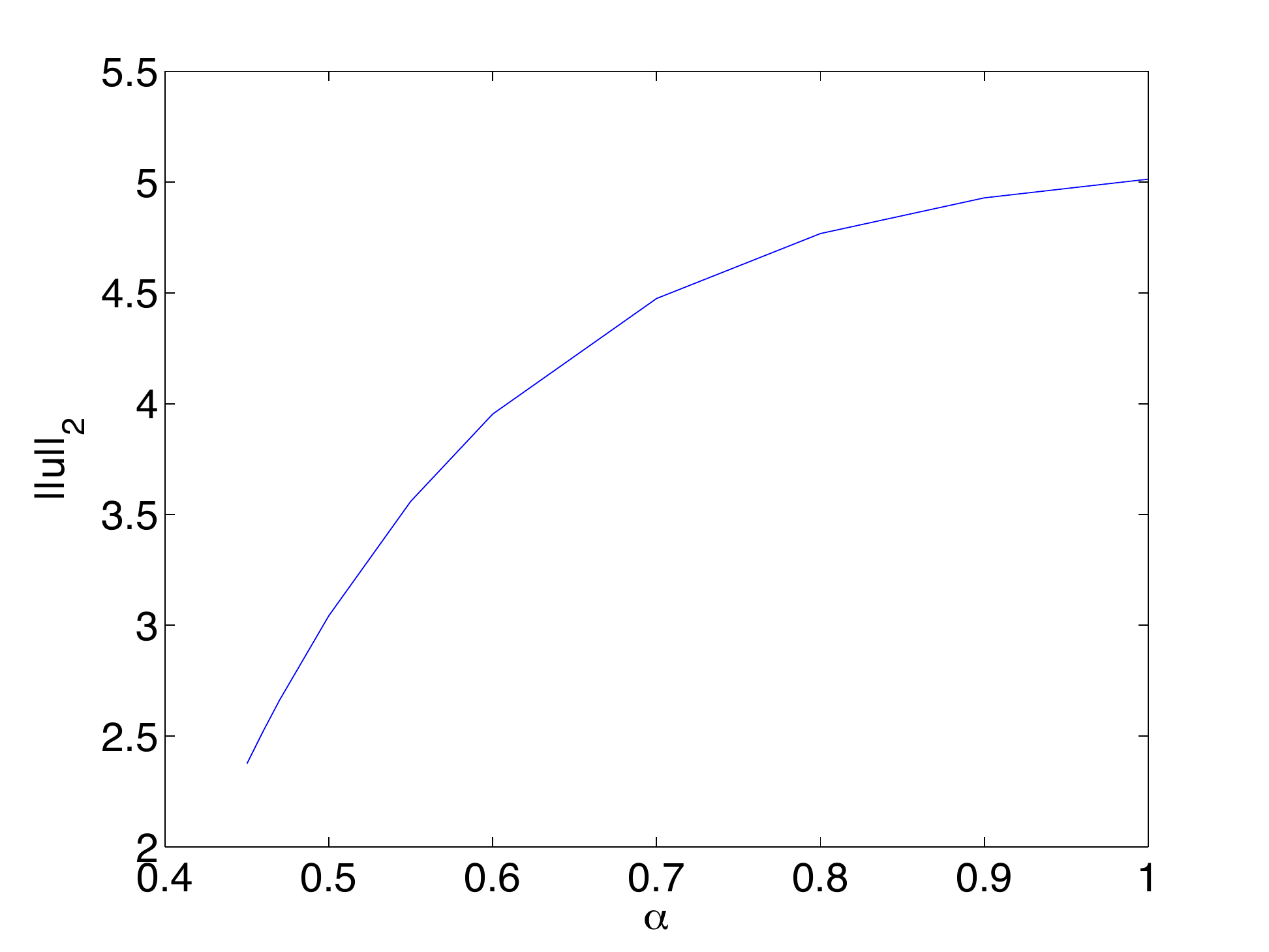}
  \includegraphics[width=0.49\textwidth]{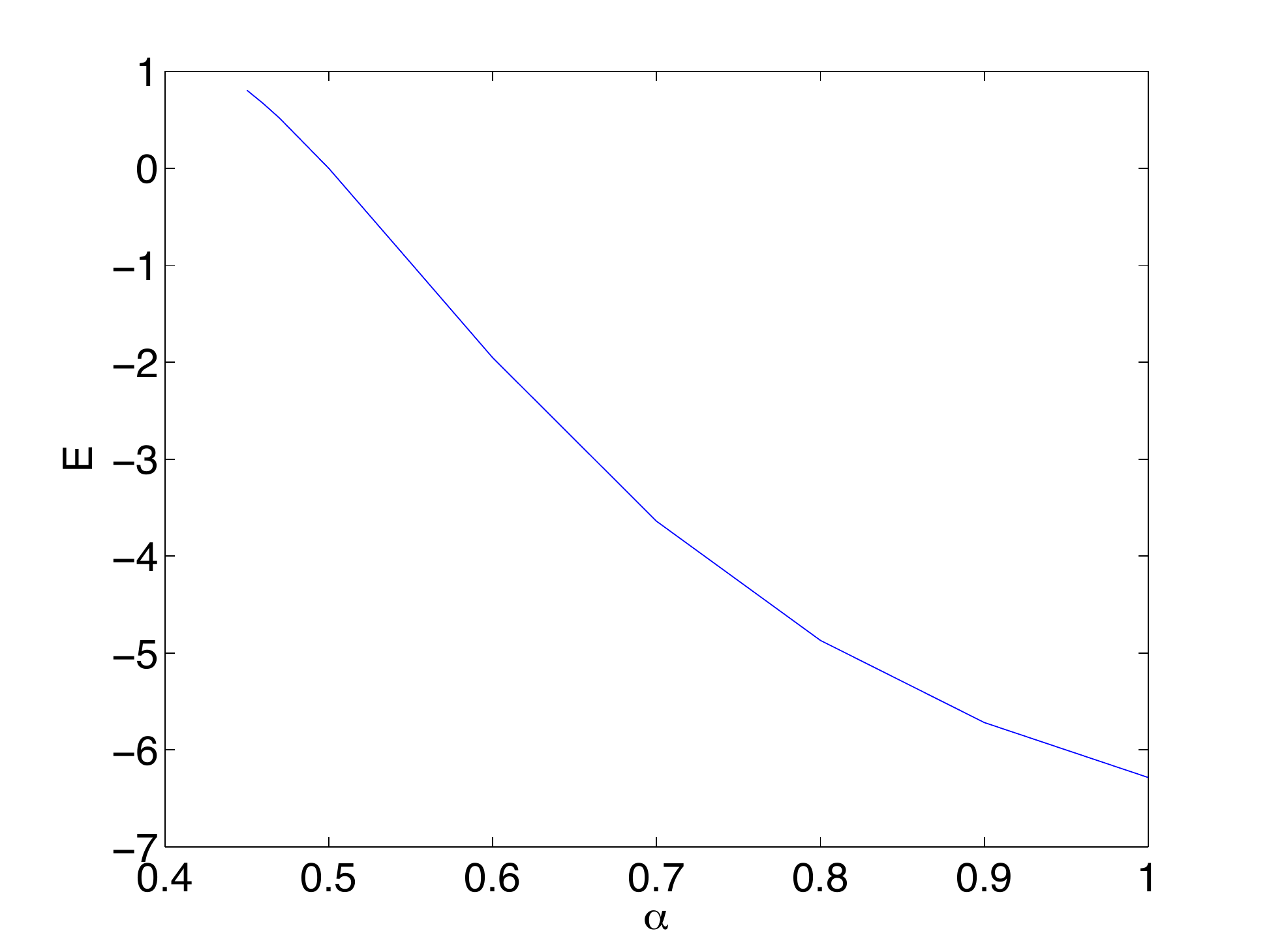}
 \caption{Mass (left) and energy (right) of the solitons in 
 Fig.~\ref{solQleg} in dependence of $\alpha$.}
 \label{solmass}
\end{figure}

In Fig.~\ref{solmass} we trace the mass and the energy of the 
solitons for $c=1$ in dependence of $\alpha$. It can be seen that the 
mass is monotonically increasing with $\alpha$, whereas the energy 
decreases. For $\alpha=0.5$ we get $M\sim3.043$  and $E=-0.002$. The 
latter value is compatible with zero for the precision with which we 
determine the soliton. 

\begin{figure}[htb!]
  \includegraphics[width=0.6\textwidth]{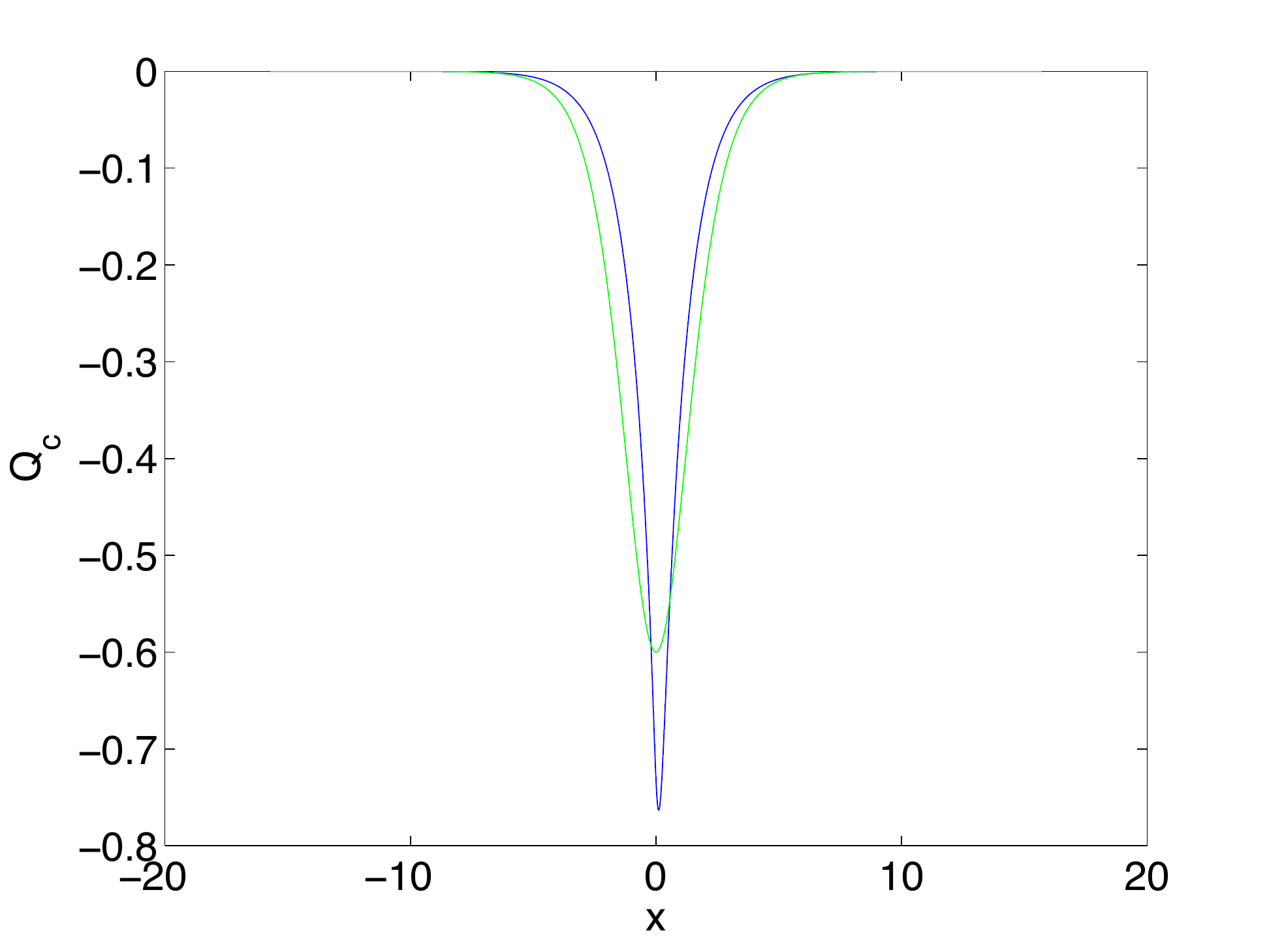}
 \caption{Soliton of the Whitham equation for $c=-1.2$ in blue and 
 the corresponding KdV soliton in green.}
 \label{solwhitham}
\end{figure}

It is possible to construct in a similar way solitons to the Whitham 
equation (\ref{Whit}). These solitons satisfy equation
(\ref{solitarywave}) if one replaces $D^{\alpha}$ with the operator 
$L$ in (\ref{Whibis}) having the Fourier symbol 
$\left(\frac{\tanh(\xi)}{\xi}\right)^{1/2}$. For small $\xi$, the 
Whitham equation reduces to the KdV equation in the form 
$u_{t}+uu_{x}-\frac{1}{6}u_{xxx}=0$ which has the solitonic solution 
\begin{equation}
    Q_{c}=3(c+1)\mbox{sech}^{2}(\sqrt{-1.5(c+1)}x).
    \label{kdvsol}
\end{equation}
Thus for real 
solutions one needs $c+1<0$ which implies that the soliton is 
negative and travelling to the left in contast to the fKdV solitons 
considered above. For $c\to-1$ the amplitude of the soliton tends to 
zero. If we use the KdV soliton (\ref{kdvsol}) for negative $c+1$ close 
to zero, the iteration for the Whitham soliton converges quickly. In 
Fig.~\ref{solwhitham} we show the Whitham and the KdV soliton for 
$c+1=-1.2$. For smaller $c$ the iteration stops converging. 
A numerical treatment of the solitary waves of the 
Whitham equation along similar lines has been already presented in 
\cite{EK}. There periodic solutions were treated which tend to 
solitons in the limit of large periods. In our case, the period is 
chosen large enough that the found solution (which is as the KdV 
soliton rapidly decreasing in contrast to the fKdV solitons) should 
be equal to the soliton within numerical accuracy. In \cite{EK} it 
was also shown that there is a maximal $|c|$ of the soliton due to 
the fact that the dispersion of the Whitham equation decreases with 
$|c|$. Note that the Whitham equation in \cite{EK} is not written in 
dimensionless coordinates as here, and that the 
dispersion has a different sign. It is clear that the maximal $|c|$ 
we observe for the iteration to converge is related to the 
nonexistence of solitons of the Whitham equation for large speeds. 
But we cannot say how close we get to the limiting velocity since 
failure of an iteration to converge can also be related to an 
inadequate initial iterate. 

\subsection{Numerical study of blow-up in fractionary KdV equations}
In this subsection we will study numerically the possibility of blow-up 
in fKdV equations. We concentrate  on the case 
$0<\alpha<1$ and the initial data $u_{0}=\beta \mbox{sech}^{2}x$ 
which are motivated by the KdV soliton. Since these initial data are 
rapidly decreasing, they can be treated as essentially periodic 
as discussed above. We find that solutions for sufficiently large 
$\beta$ decompose for $\alpha>1/2$ asymptotically into solitons and 
radiation. For $\alpha=1/2$, we find blow-up for initial data with 
negative energy. The blow-up profile is given by a dynamically 
rescaled soliton as for the $L_{2}$ critical case for gKdV. In the 
supercritical case with $1/3<\alpha<1/2$, we find blow-up for initial data with 
sufficiently large mass similar to the supercritical blow-up in 
gKdV.  In the energy supercritical case $\alpha\leq 1/3$, there will 
be again blow-up, but since there are no solitons in this case, there 
will be no formation of a soliton which finally blows up, but the 
initial hump steepens as known from the Burgers' equation and finally 
blows up. 

\begin{figure}[htb!]
  \includegraphics[width=\textwidth]{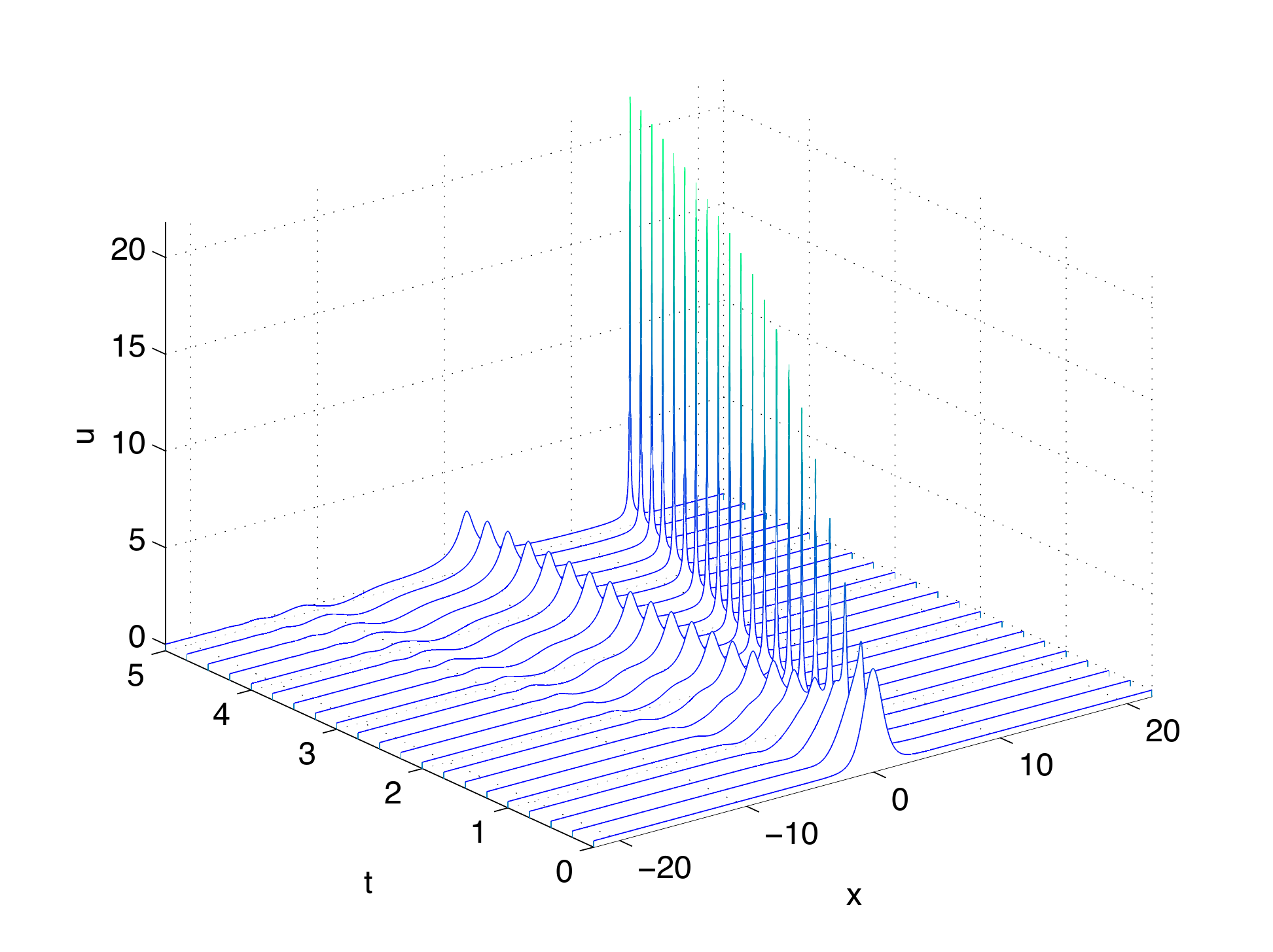}
 \caption{Solution to the fKdV equation (\ref{Cauchy}) for 
 $\alpha=0.6$ and the initial data $u_{0}=5\mbox{sech}^{2}x$.}
 \label{gBO3sechalpha064t}
\end{figure}

We first study an example for $\alpha=0.6$, i.e., $1/2<\alpha<1$. No 
blow-up is expected in this case. For the initial data 
$u_{0}=5\mbox{sech}^{2}x$, we get the solution shown in 
Fig.~\ref{gBO3sechalpha064t}. The initial hump gets laterally 
compressed and increases in height. At a given point it splits into 
two humps which both continue to grow in height. The humps seem to 
approach solitons since their speed becomes almost constant after 
some time and their shape stops changing. Thus it appears that an 
initial pulse of sufficient size decomposes into solitons. The 
remaining energy is radiated away in the form of 
oscillations propagating to the left whilst the humps travel to the 
right. The computation is carried out with $N=2^{14}$ Fourier modes 
for $x\in7[-\pi,\pi]$ and $N_{t}=10^{4}$ time steps. The relative 
energy is conserved to the order of $10^{-12}$. 

\begin{figure}[htb!]
  \includegraphics[width=0.49\textwidth]{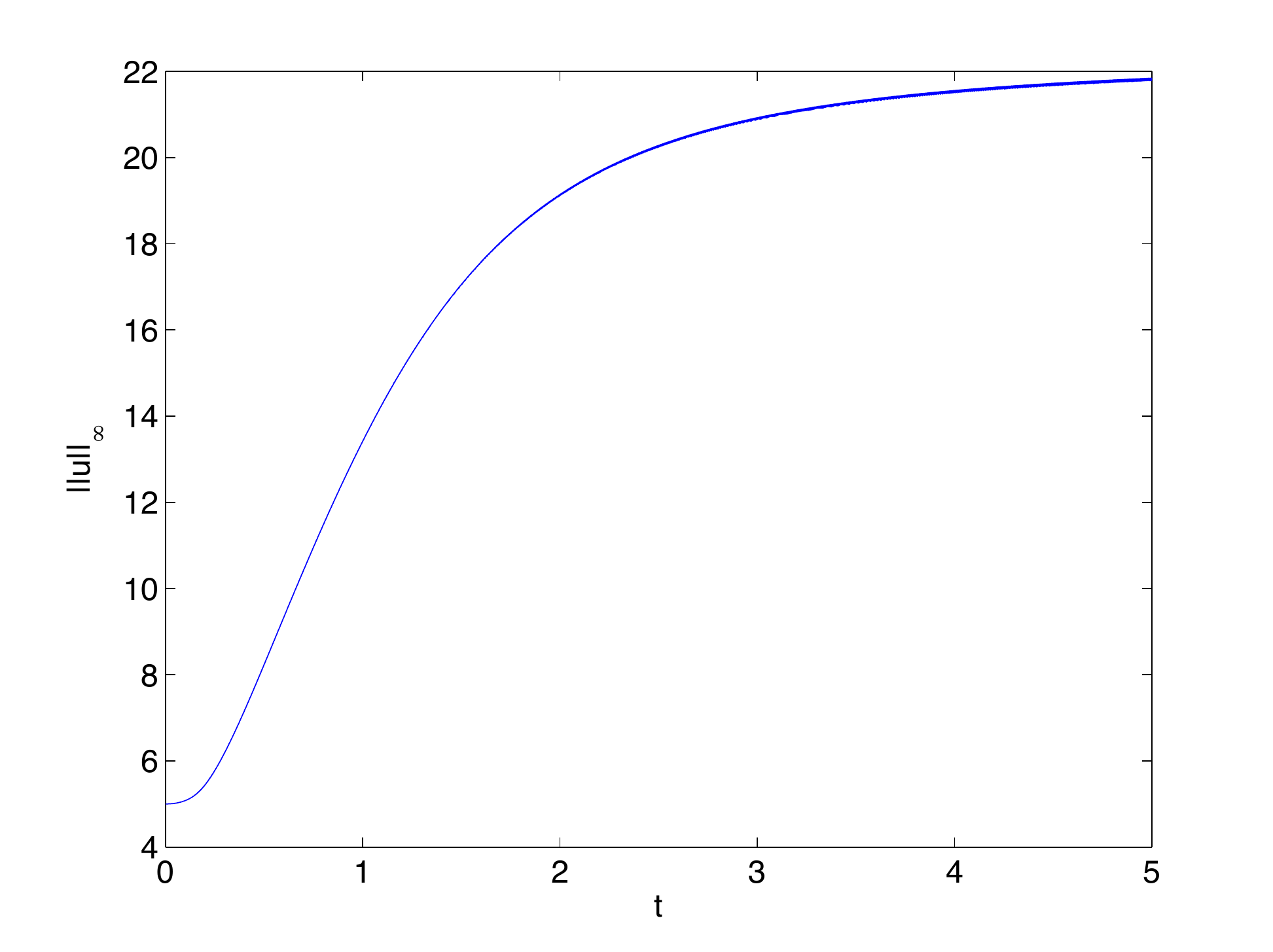}
  \includegraphics[width=0.49\textwidth]{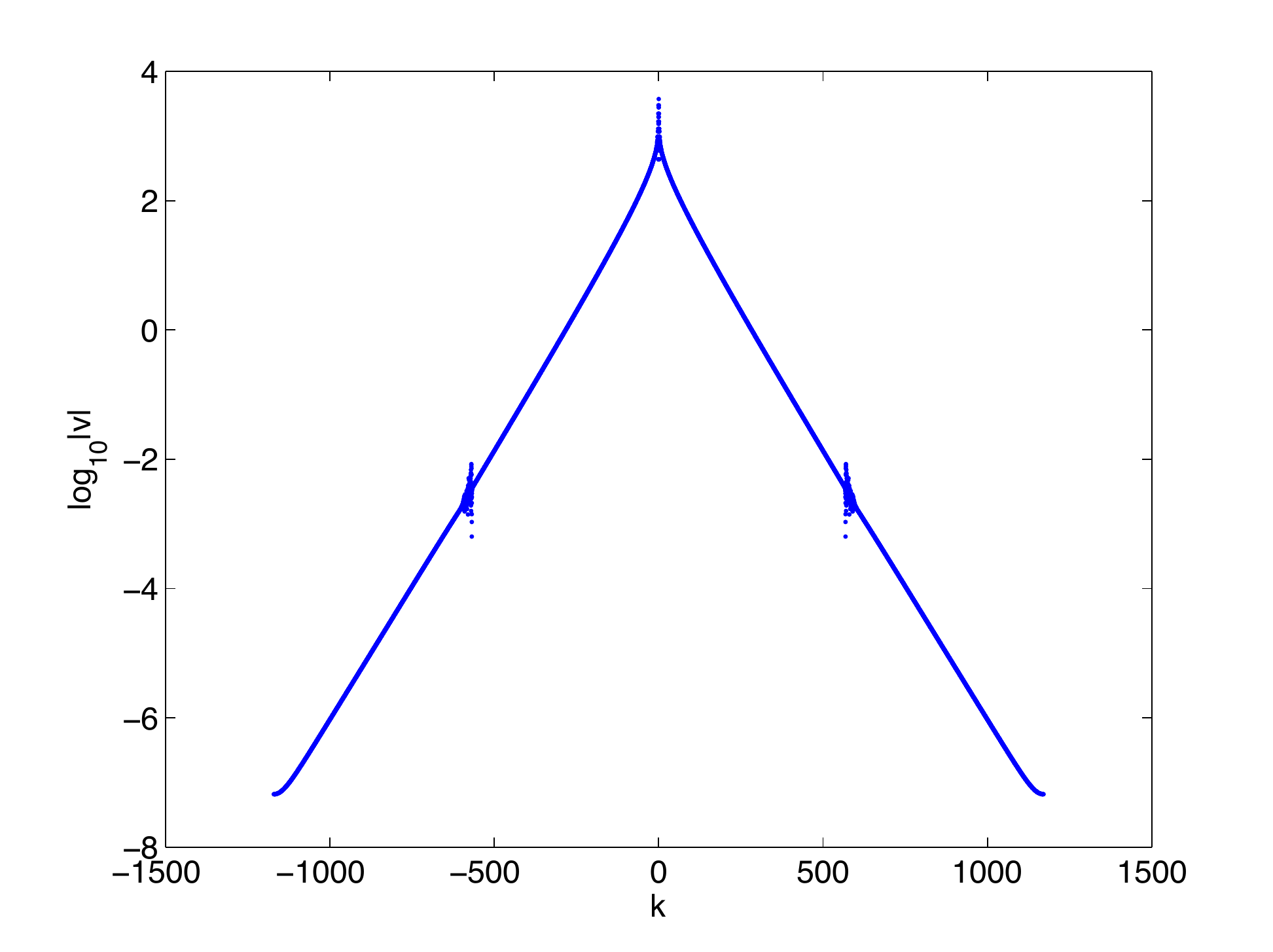}
 \caption{$L_{\infty}$ norm of the solution  
to the fKdV equation (\ref{Cauchy}) for 
 $\alpha=0.6$ and the initial data $u_{0}=5\mbox{sech}^{2}x$
 in dependence of time on the left, and the modulus of the Fourier 
 coefficients of the solution for $t=5$ on the right.}
 \label{gBO3sechalpha06}
\end{figure}

Note that the  dispersive oscillations in 
Fig.~\ref{gBO3sechalpha064t} decay only very slowly in amplitude. 
Due to the imposed periodic boundary conditions, 
this leads to oscillations also to the right of the initial hump 
where there would be none on an infinite domain. A consequence of 
these oscillations are small oscillations in the $L_{\infty}$ norm of the 
solution in Fig.~\ref{gBO3sechalpha064t} which is shown in 
Fig.~\ref{gBO3sechalpha06}. Nonetheless it is clear that the 
$L_{\infty}$ norm of the solution after some strong initial increase 
appears to reach a plateau, probably corresponding to some 
asymptotically solitonic solution. Though the energy of the initial data is 
negative, there is no reason to assume that there 
will be blow-up in this case. Since there is 
sufficient resolution in Fourier space (also shown in 
Fig.~\ref{gBO3sechalpha06}) and since the computed energy is conserved 
to the order of $10^{-12}$, the numerical evidence for this is quite 
strong. 

\begin{figure}[htb!]
  \includegraphics[width=0.7\textwidth]{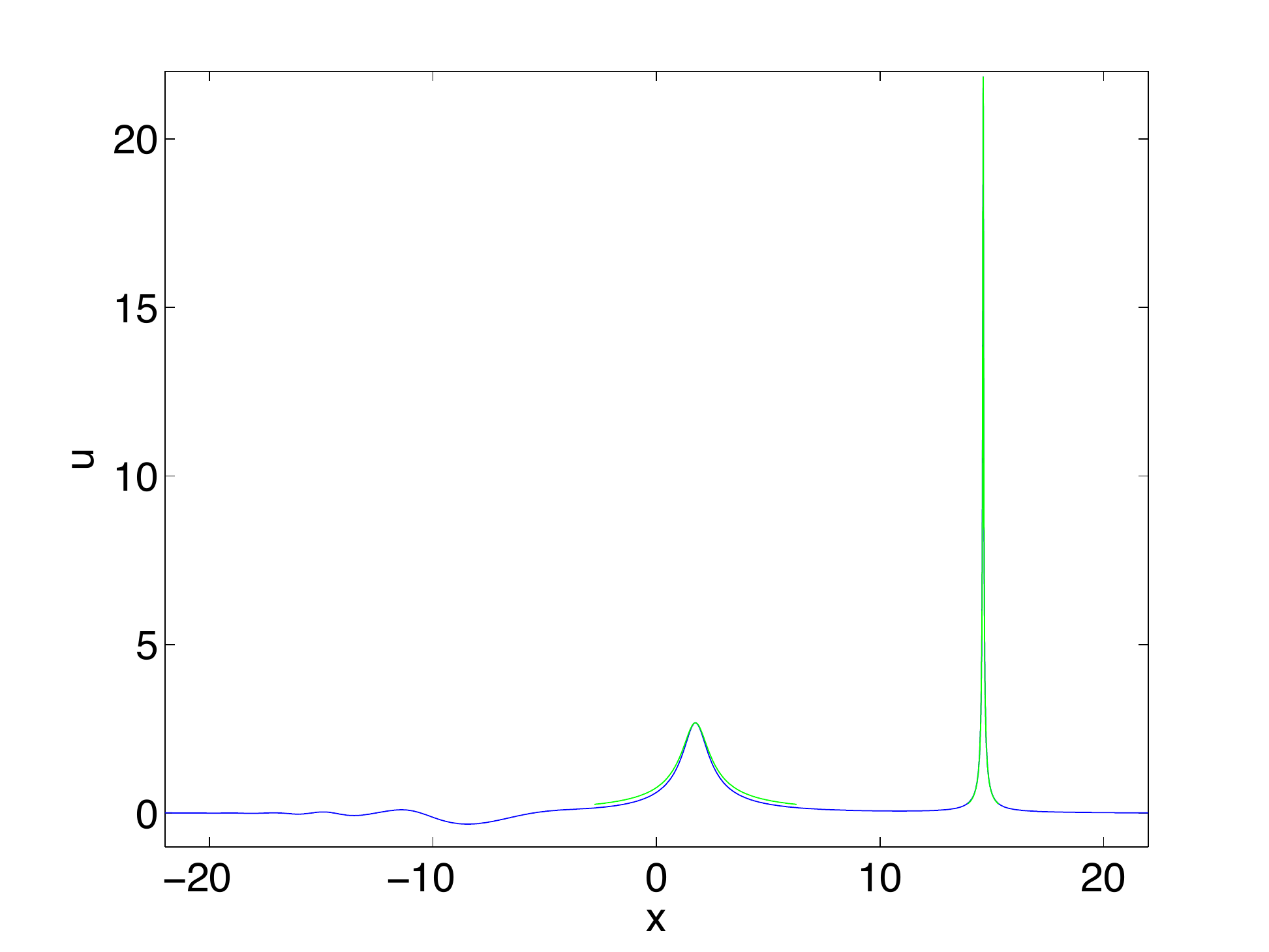}
 \caption{Solution to the fKdV equation (\ref{Cauchy}) for 
 $\alpha=0.6$ and the initial data $u_{0}=5\mbox{sech}^{2}x$ for 
 $t=5$ in blue, fitted solitons at the humps in green.}
 \label{gBO5sechalpha06solfit}
\end{figure}

To test further whether the humps in Fig.~\ref{gBO3sechalpha064t} 
are related to solitons, we fit the solitonic solutions constructed 
numerically in the previous subsection to the humps at the final 
computed time. To this end we determine the locations and the value of 
the maxima, shift the numerically constructed soliton solution to 
these values of $x$ and rescale via (\ref{solresc}). This leads to 
Fig.~\ref{gBO5sechalpha06solfit}. Obviously the final time is not yet 
fully in the asymptotic regime for the smaller soliton on the left, 
but the agreement is already very good there. For the larger soliton, 
the solutions are hardly distinguishable. Thus it appears that 
solutions to the fKdV equation decompose for large times into 
solitons and radiation.

\begin{figure}[htb!]
  \includegraphics[width=\textwidth]{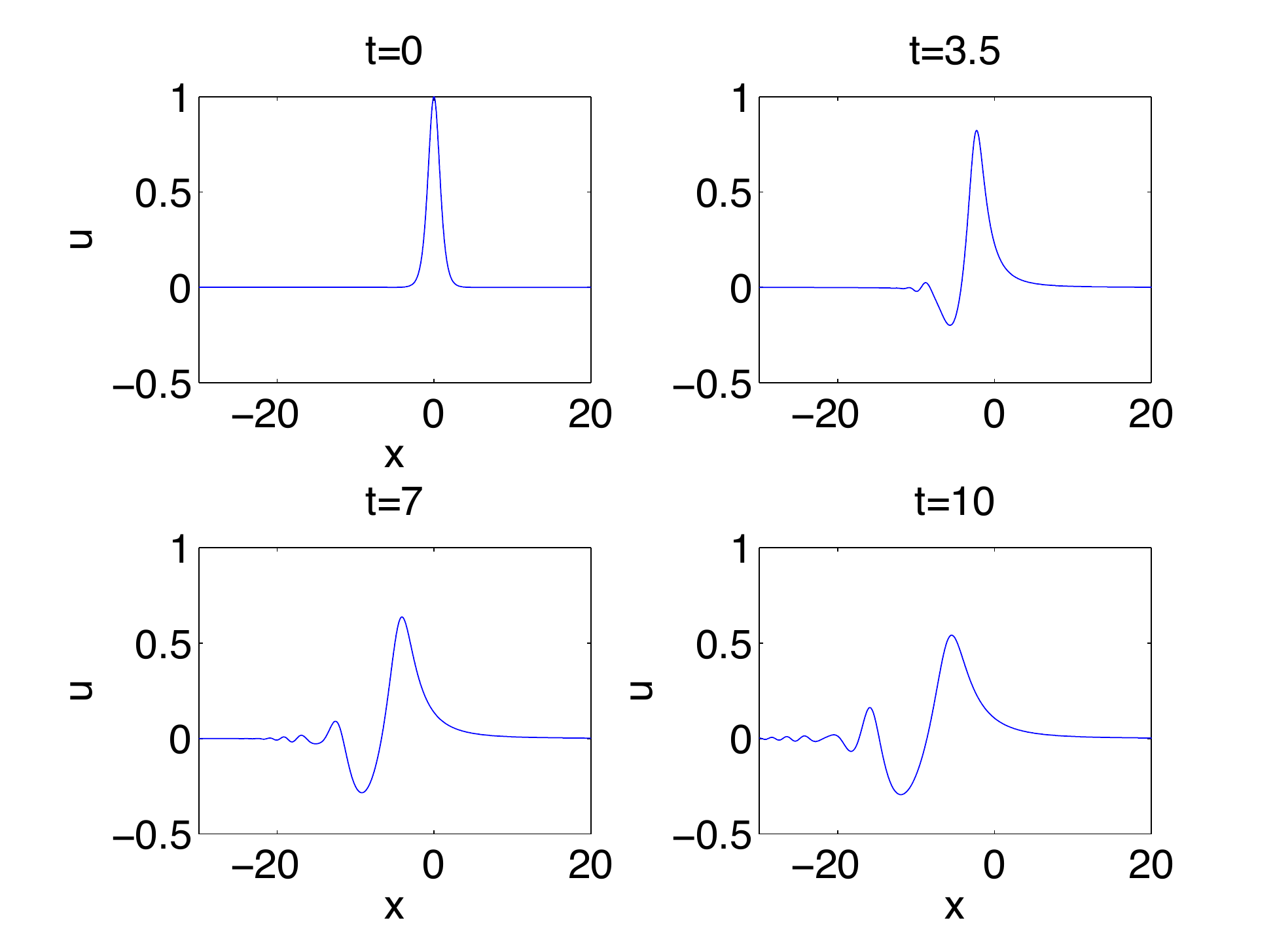}
 \caption{Solution to the fKdV equation (\ref{Cauchy}) for 
 $\alpha=0.5$ and the initial data $u_{0}=\mbox{sech}^{2}x$ 
 for several values of $t$.}
 \label{gBOsechalpha054t}
\end{figure}

There are also no indications for blow-up for the critical case $\alpha=0.5$ and the initial data 
$u_{0}=\mbox{sech}^{2} x$ which can be seen in Fig.~\ref{gBOsechalpha054t} 
and for which the 
energy is positive (and thus larger than the soliton energy which is 
equal to zero) and for which the mass ($M[u_{0}]=4/3$) is smaller  
than the soliton mass $M\sim3.043$. We use $N=2^{14}$ Fourier modes 
for $x\in10[-\pi,\pi]$ and $N_{t}=10^{3}$ time steps for the computation. 
The solution is well 
resolved in Fourier space and the computed energy is conserved to the 
order of $10^{-12}$. Dispersive oscillations form 
immediately, and after a short increase the $L_{\infty}$ norm of the 
solution appears to decrease monotonically. The initial hump seems to 
be radiated away, it is below the threshold to form a 
soliton. 

\begin{figure}[htb!]
  \includegraphics[width=\textwidth]{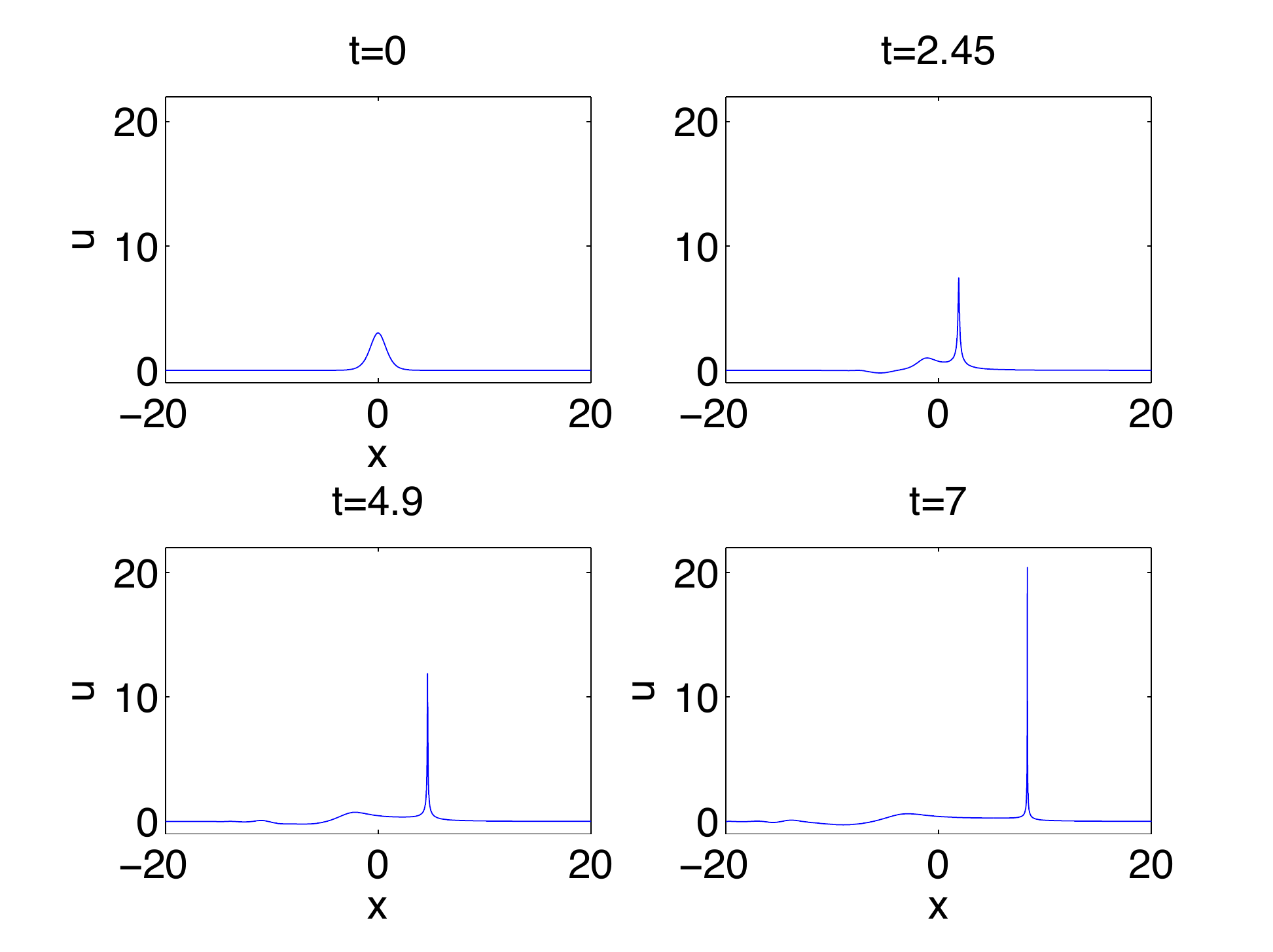}
 \caption{Solution to the fKdV equation (\ref{Cauchy}) for 
 $\alpha=0.5$ and the initial data $u_{0}=3\mbox{sech}^{2}x$ 
 for several values of $t$.}
 \label{gBO3sechalpha054t}
\end{figure}

The picture changes completely for the same $\alpha$ and the initial 
data $u_{0}=3\mbox{sech}^{2}x$ for which the energy is negative and 
for which the mass $M[u_{0}]=12$ is larger than the soliton mass. The 
solution for this case can be seen in Fig.~\ref{gBO3sechalpha054t}. 
The initial hump gets again laterally compressed and increases in 
height, a soliton seems to emerge. 
But this time it appears that the oscillations forming at the 
same time and propagating to the left cannot compensate the increase 
in height, the `soliton' seems to be unstable and eventually blows up. 
We compute with  $N=2^{16}$ Fourier modes 
for $x\in20[-\pi,\pi]$ and $N_{t}=10^{4}$ time steps.

\begin{figure}[htb!]
   \includegraphics[width=0.49\textwidth]{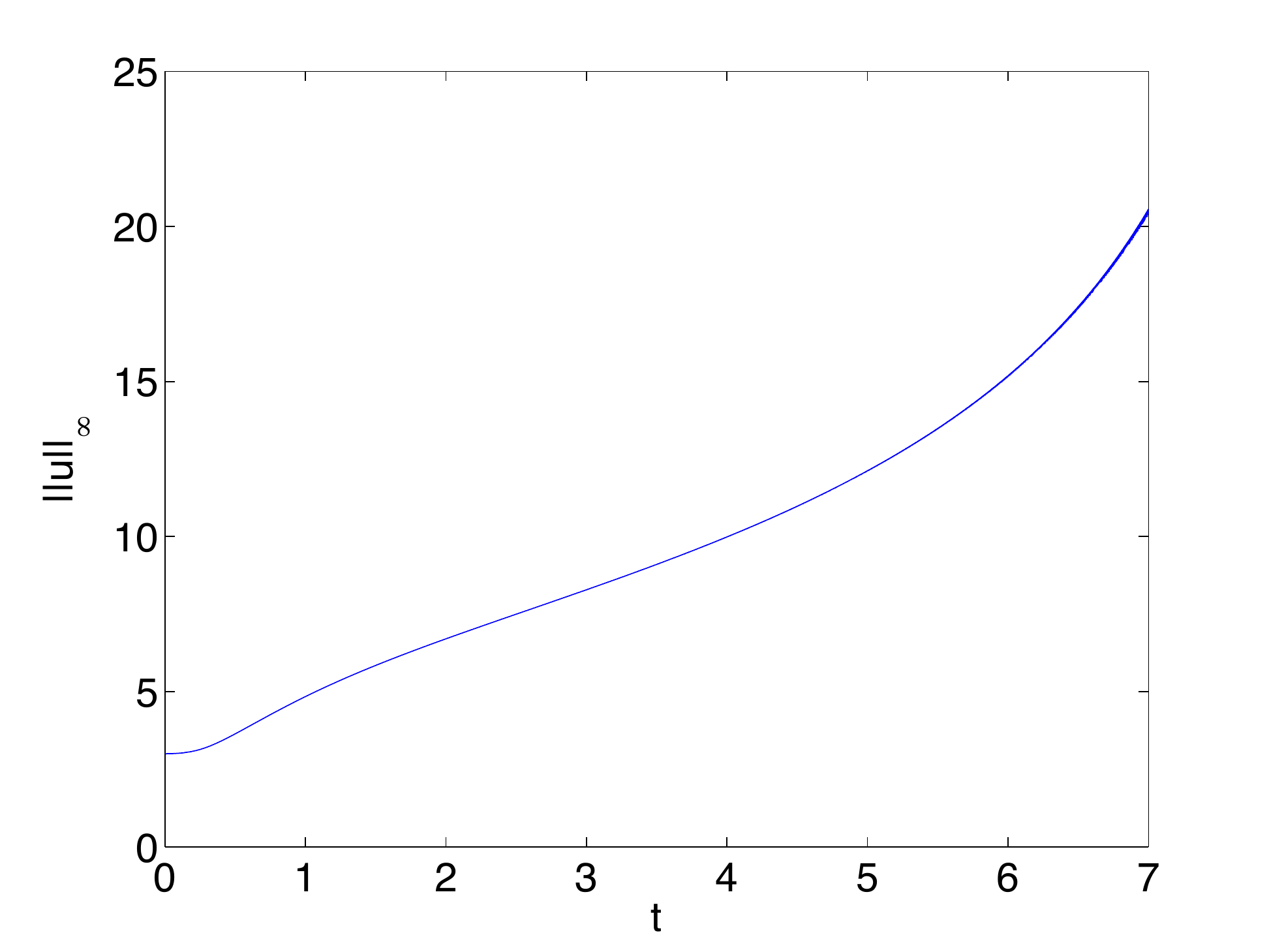}
  \includegraphics[width=0.49\textwidth]{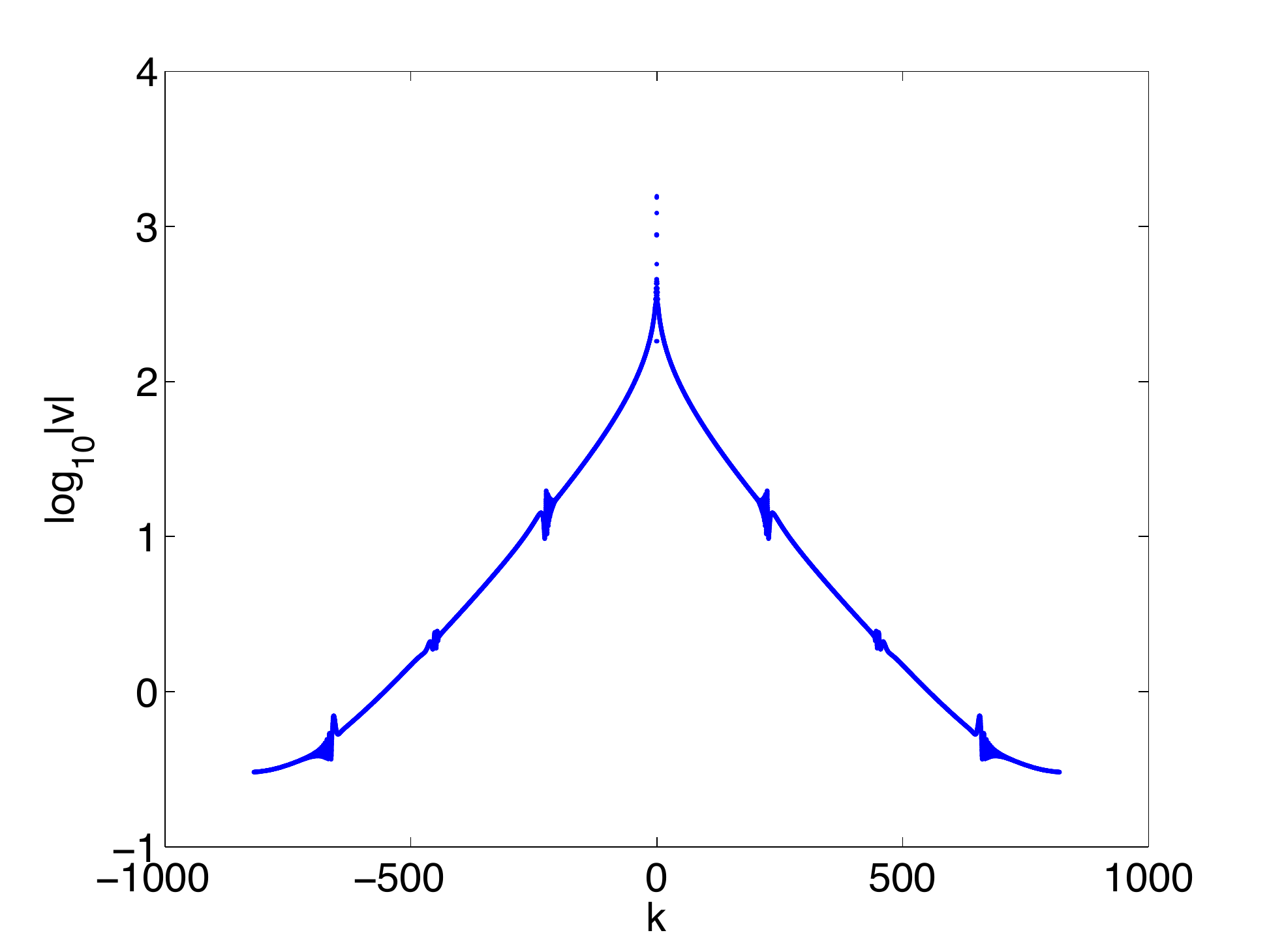}
 \caption{$L_{\infty}$ norm of the solution  
to the fKdV equation (\ref{Cauchy}) for 
 $\alpha=0.5$ and the initial data $u_{0}=3\mbox{sech}^{2}x$
 in dependence of time on the left, and the modulus of the Fourier 
 coefficients of the solution for $t=7$ on the right.}
 \label{gBO3sechalpha05}
\end{figure}

We stop the 
code at $t=7$ where the computed relative energy drops below $10^{-3}$ 
which implies that the solution is no longer reliable for larger 
times. Note that the time where the code is stopped for this reason 
is not taken as the blow-up time. The latter is determined below via 
a fitting to a theoretical 
blow-up scenario. As can be seen in Fig.~\ref{gBO3sechalpha05}, the 
loss of accuracy due to a 
lack of resolution in Fourier space due to the strongly peaked hump, 
not to a lack in resolution in time. It is also clear from 
Fig.~\ref{gBO3sechalpha05} that the $L_{\infty}$ norm of the solution 
increases monotonically. As in the $L_{2}$ critical case for gKdV in 
\cite{KP2013}, this increase appears to be algebraic in time. The peak 
which appears to blow up eventually gets more and more compressed 
laterally, grows in height and propagates faster. Thus it 
seems to be similar to the blow-up in the $L_{2}$ critical case for 
gKdV, also in the sense that the solution could blow up  at $x^{*}$ 
with $x^{*}\to\infty$. 

\begin{figure}[htb!]
   \includegraphics[width=0.49\textwidth]{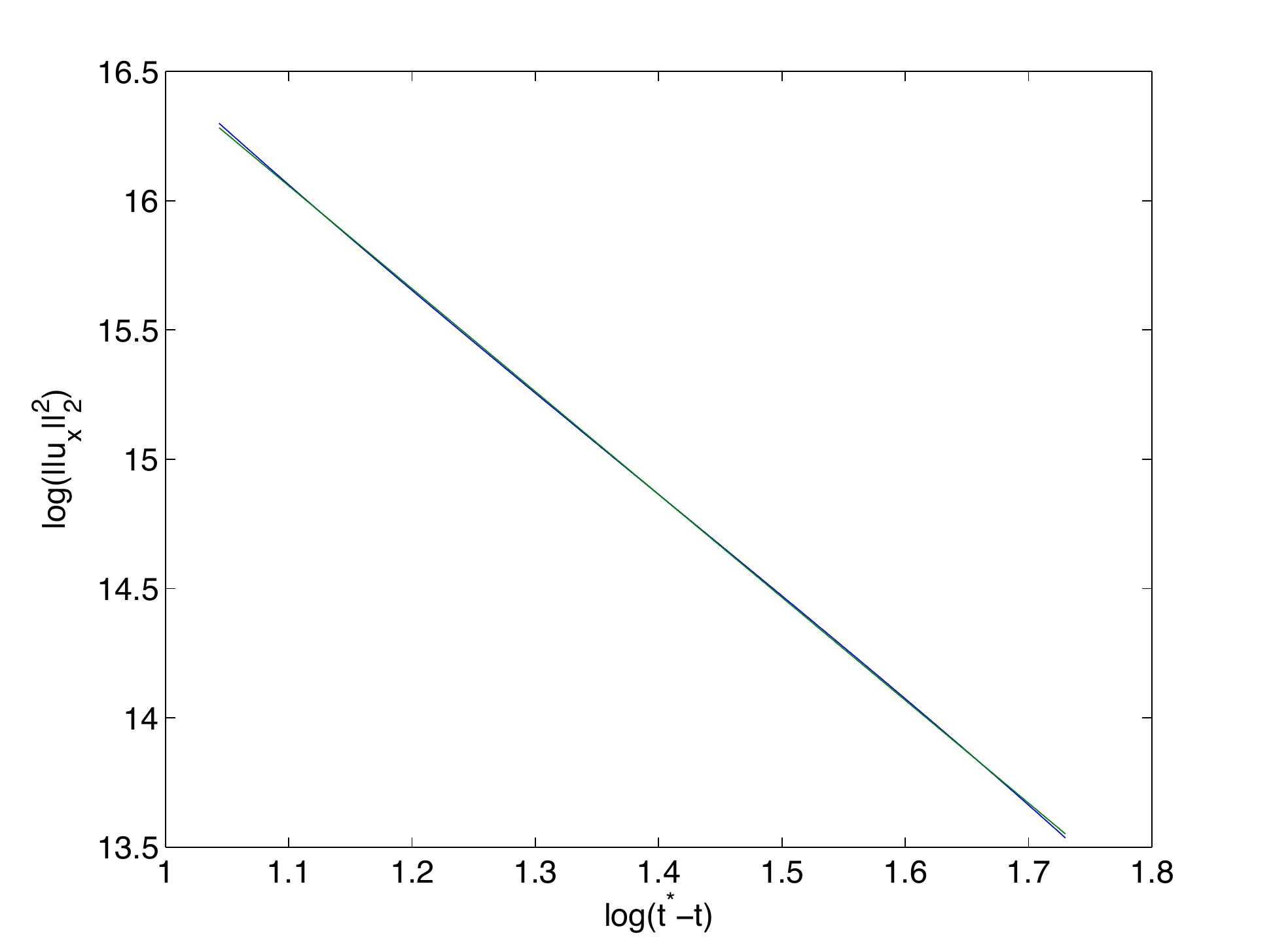}
  \includegraphics[width=0.49\textwidth]{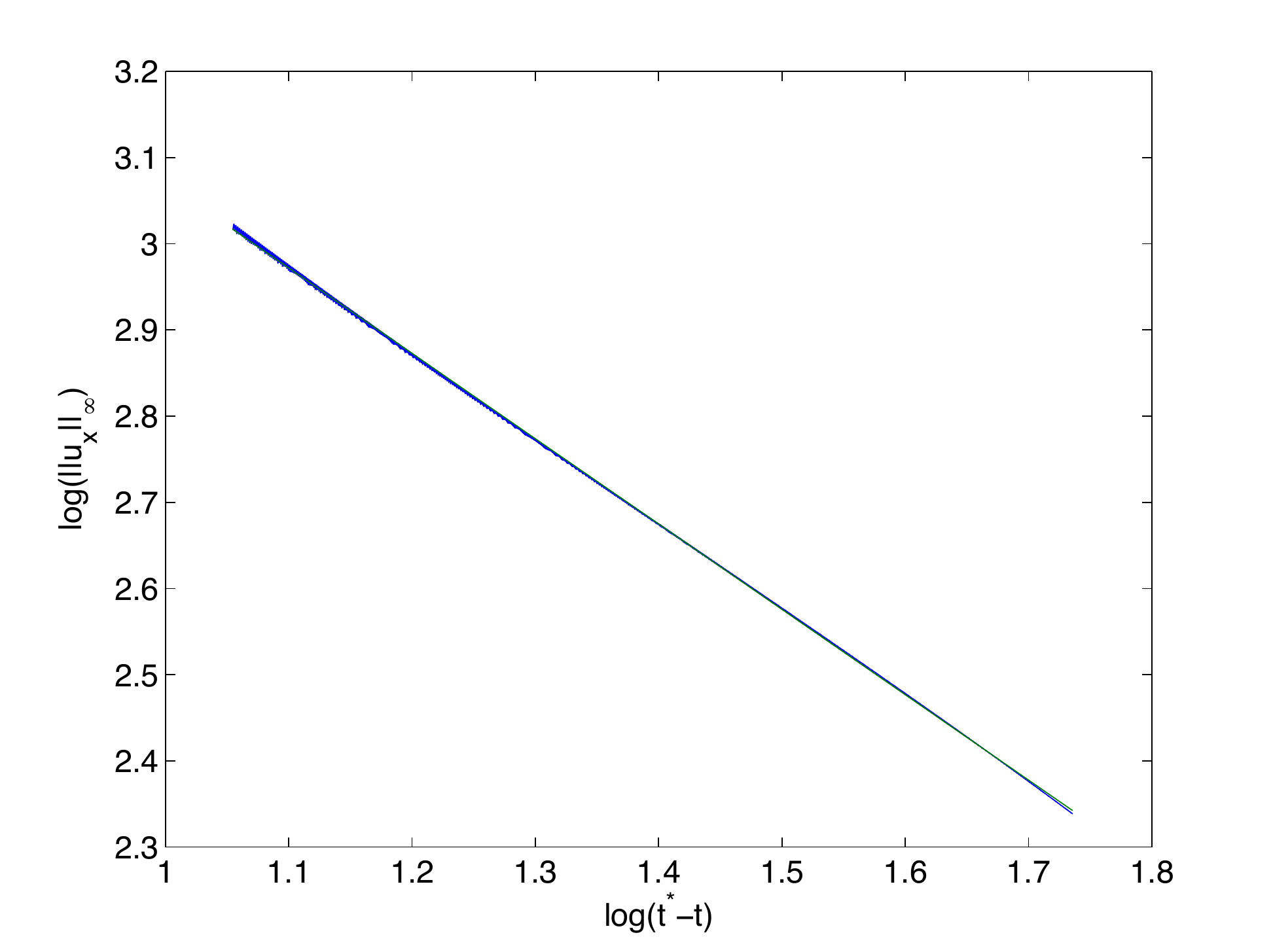}
 \caption{$L_{2}$ norm of the gradient (left) and 
 $L_{\infty}$ norm of the solution  
to the fKdV equation (\ref{Cauchy}) for 
 $\alpha=0.5$ and the initial data $u_{0}=3\mbox{sech}^{2}x$ (right) 
 for $t>4.1993$ in blue and the fitted lines $\kappa_{1}\ln 
 (t^{*}-t)+\kappa_{2}$ in green.}
 \label{gBO3sechalpha05fit}
\end{figure}

For $t$ close enough to the blow-up time $t^{*}$, it is possible to 
fit the data to the expected scalings (\ref{L2scal}) and 
(\ref{genscal}) to characterize the type of blow-up 
further. The task is to find the range of the values of $t$ where the 
asymptotic behavior can be already observed whilst still allowing for 
a large enough interval to have sufficient computed values for the 
fitting. The 
problem is that the blow-up time $t^{*}$ is not known and has to be determined in 
the process. In \cite{KP2013} it was shown that this can be done with the 
optimization algorithm
\cite{fminsearch} distributed with Matlab as \emph{fminsearch}. In 
Fig.~\ref{gBO3sechalpha05fit} we show the fits for the $L_{2}$ norm 
of $u_{x}$ and the $L_{\infty}$ norm of $u$ for $4.1993<t<7$. The 
results of the fitting do not change much if the lower bound is 
slightly changed. For both quantities the logarithms are fitted to 
$\kappa_{1}\ln (t^{*}-t)+\kappa_{2}$. For $||u||_{\infty}$ we 
find $t^{*}=9.8712$, $\kappa_{1}=-0.9901$ and    $\kappa_{2}= 
4.0609    $, for $||u_{x}||_{2}^{2}$ we obtain  $t^{*}=9.8393$, 
$\kappa_{1}=-3.9785$ and $\kappa_{2}=   20.4335$. Thus both fittings 
are consistent within numerical accuracy, and are compatible with 
(\ref{L2scal}) for $\gamma=1$ which is exactly the value for the 
$L_{2}$ critical gKdV.     

\begin{figure}[htb!]
   \includegraphics[width=0.7\textwidth]{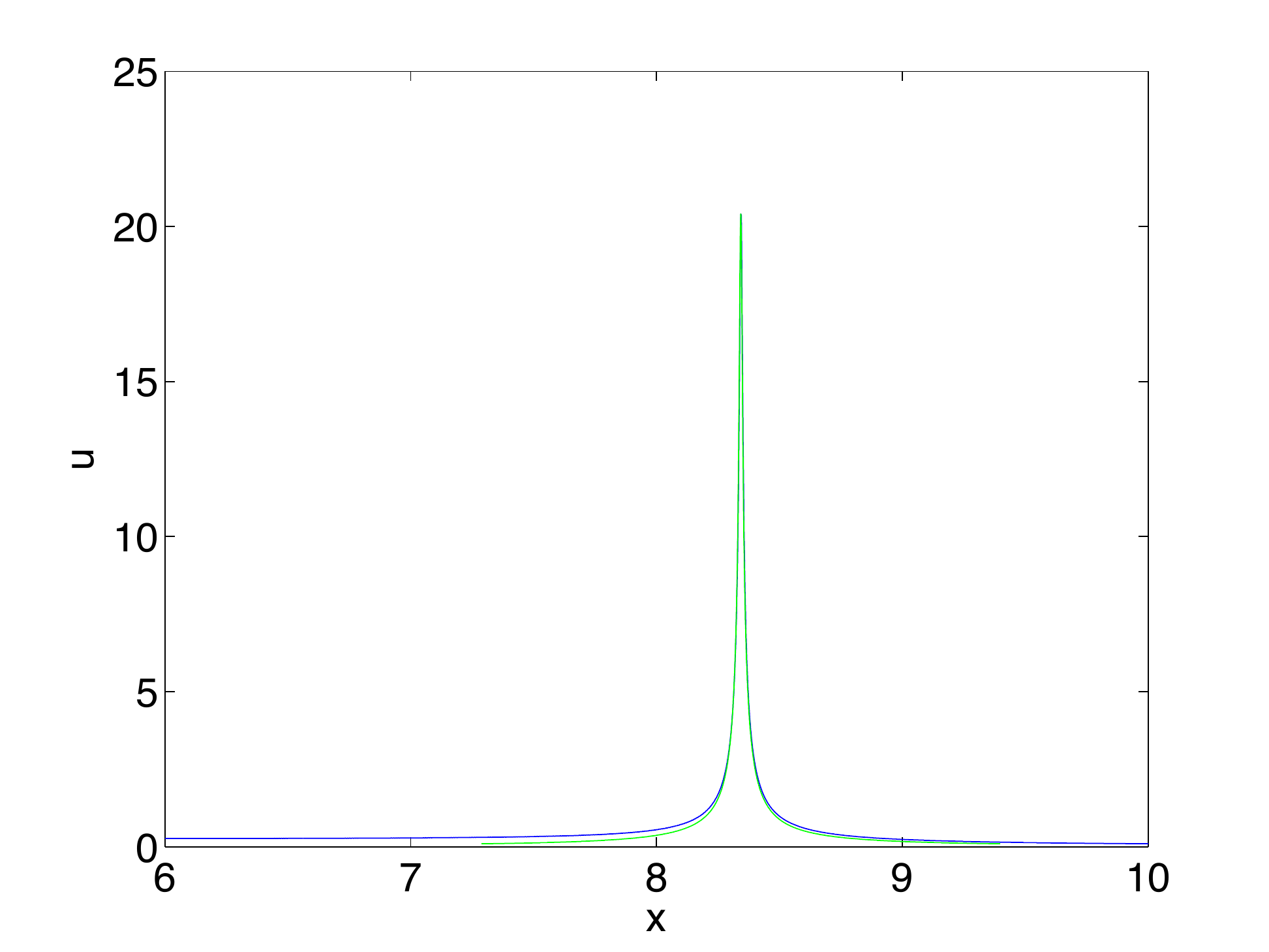}
 \caption{Solution to the fKdV equation (\ref{Cauchy}) for 
 $\alpha=0.5$ and the initial data $u_{0}=3\mbox{sech}^{2}x$ 
 for  $t=7$ in blue and the fitted soliton (\ref{solitarywave}) 
 rescaled according to (\ref{resc}) in green.}
 \label{gBO3sechalpha05solfit}
\end{figure}

There are, however, more implications from the dynamic rescaling 
(\ref{resc}). Since the dependence of the norms in 
Fig.~\ref{gBO3sechalpha05fit} appears to be algebraic in the rescaled 
time $\tau$, the blow-up profile should be asymptotically given by a 
rescaled soliton (\ref{solitarywave}). To test this we fit the hump 
in the last frame of Fig.~\ref{gBO3sechalpha054t} according to 
(\ref{resc}) by reading the value of $L$ off from the maximum of the 
solution. The 
result of this fitting to the numerically obtained soliton can be 
seen in Fig.~\ref{gBO3sechalpha05solfit}. This fitting is analogous 
to the one for the case $\alpha=0.6$ in 
Fig.~\ref{gBO5sechalpha06solfit}, the difference being that for 
$\alpha>0.5$, humps of sufficient size asymptotically approach an 
exact fKdV soliton, whereas for $\alpha=0.5$ a hump of sufficient size will simply 
blow up with a profile of a dynamically rescaled soliton. The good 
agreement of the fitting indicates that blow-up in the solution to the $L_{2}$ 
critical fKdV equations is indeed given by a dynamically rescaled 
soliton. This implies that most of the $L_{2}$ norm of the initial data will 
be asymptotically concentrated in the rescaled soliton. 

Note that
it cannot be decided numerically whether energy zero marks a dividing 
line between radiation and blow-up. We can only state that the 
solution behaves decisively differently for clearly positive and 
negative energies. For the $L_{2}$ critical case $\alpha=0.5$, it 
appears that negative energy (the soliton has vanishing energy) or a 
mass greater than the soliton mass are as for gKdV with $n=4$ the 
criterion for blow-up.

\begin{figure}[htb!]
   \includegraphics[width=0.5\textwidth]{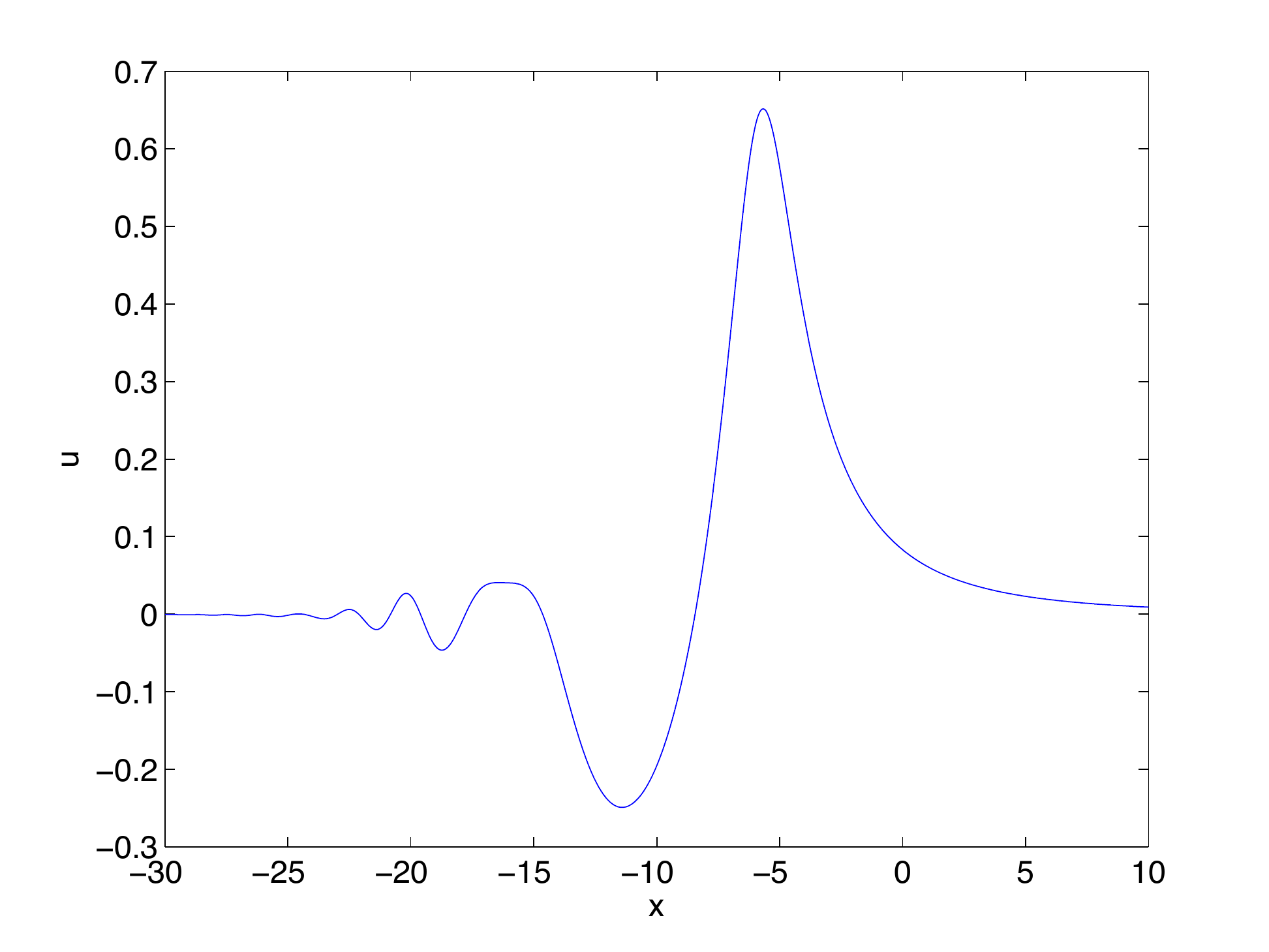}
 \caption{Solution to the fKdV equation (\ref{Cauchy}) for 
 $\alpha=0.45$ for the initial data $u_{0}=\mbox{sech}^{2}x$  
 at $t=10$.}
 \label{gBOsechalpha04}
\end{figure}

For $1/3<\alpha<0.5$, i.e., the $L_{2}$ supercritical case, there are 
certain similarities to the $L_{2}$ critical case. For initial data 
with a mass smaller than the soliton mass, the initial hump just 
appears to be radiated away as can be seen for 
$u_{0}=\mbox{sech}^{2}x$ and $\alpha=0.45$ in 
Fig.~\ref{gBOsechalpha04}. The mass of these data 
$M[u_{0}]=4/3$ is clearly smaller than the soliton mass which is 
roughly  $2.375$.

\begin{figure}[htb!]
   \includegraphics[width=\textwidth]{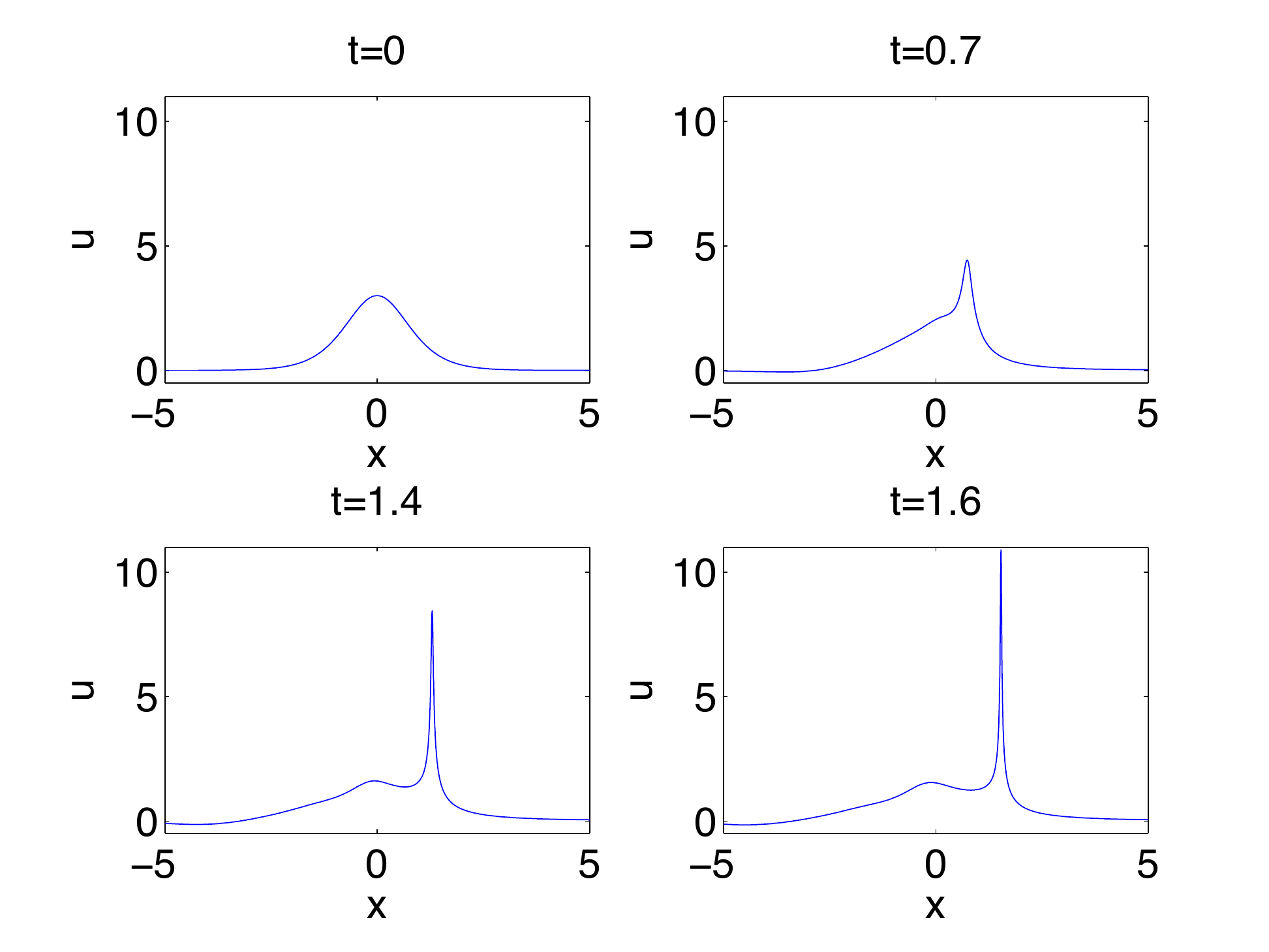}
 \caption{Solution to the fKdV equation (\ref{Cauchy}) for 
 $\alpha=0.45$ for the initial data $u_{0}=3\mbox{sech}^{2}x$  
 for several values of $t$.}
 \label{gBO3sechalpha045}
\end{figure}

In Fig.~\ref{gBO3sechalpha045} we show the solution to the fKdV 
equation for the initial data $u_{0}=3\mbox{sech}^{2}x$ for which the 
mass is $M[u_{0}]=12$ and thus much larger than the soliton mass 
$2.375$. It can be seen that a soliton forms in the evolution of the 
initial hump which separates from the remainder of the initial data 
to finally  blow up. Visibly the type of blow-up is different from 
the mass critical case in 
Fig.~\ref{gBO3sechalpha05solfit} since the $L_{2}$ norm of the 
initial data is no longer concentrated in the blow-up profile. This 
is the same type of behavior known from supercritical blow-up in gKdV 
equations, see for instance \cite{KP2013}. 

\begin{figure}[htb!]
   \includegraphics[width=0.49\textwidth]{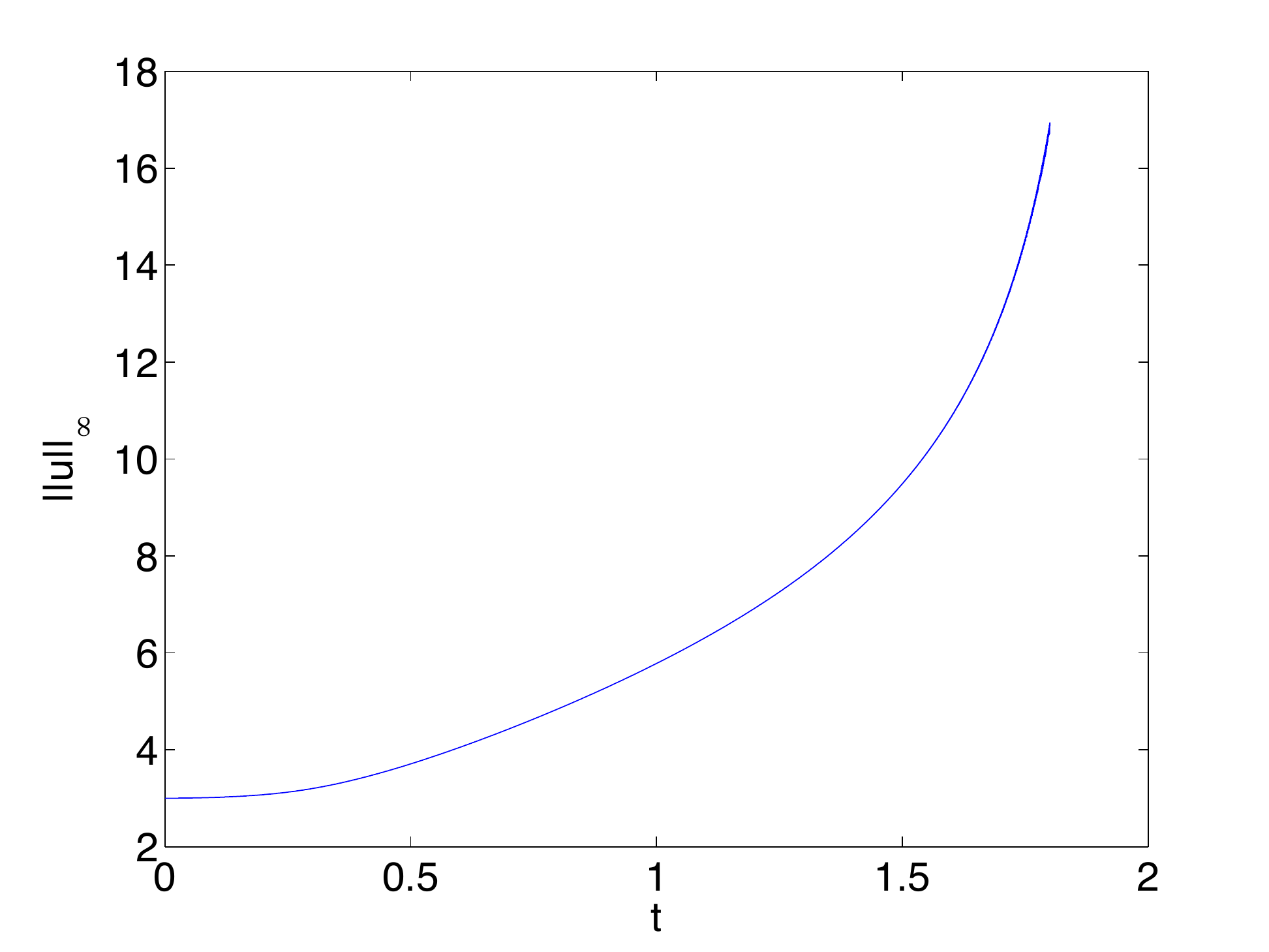}
  \includegraphics[width=0.49\textwidth]{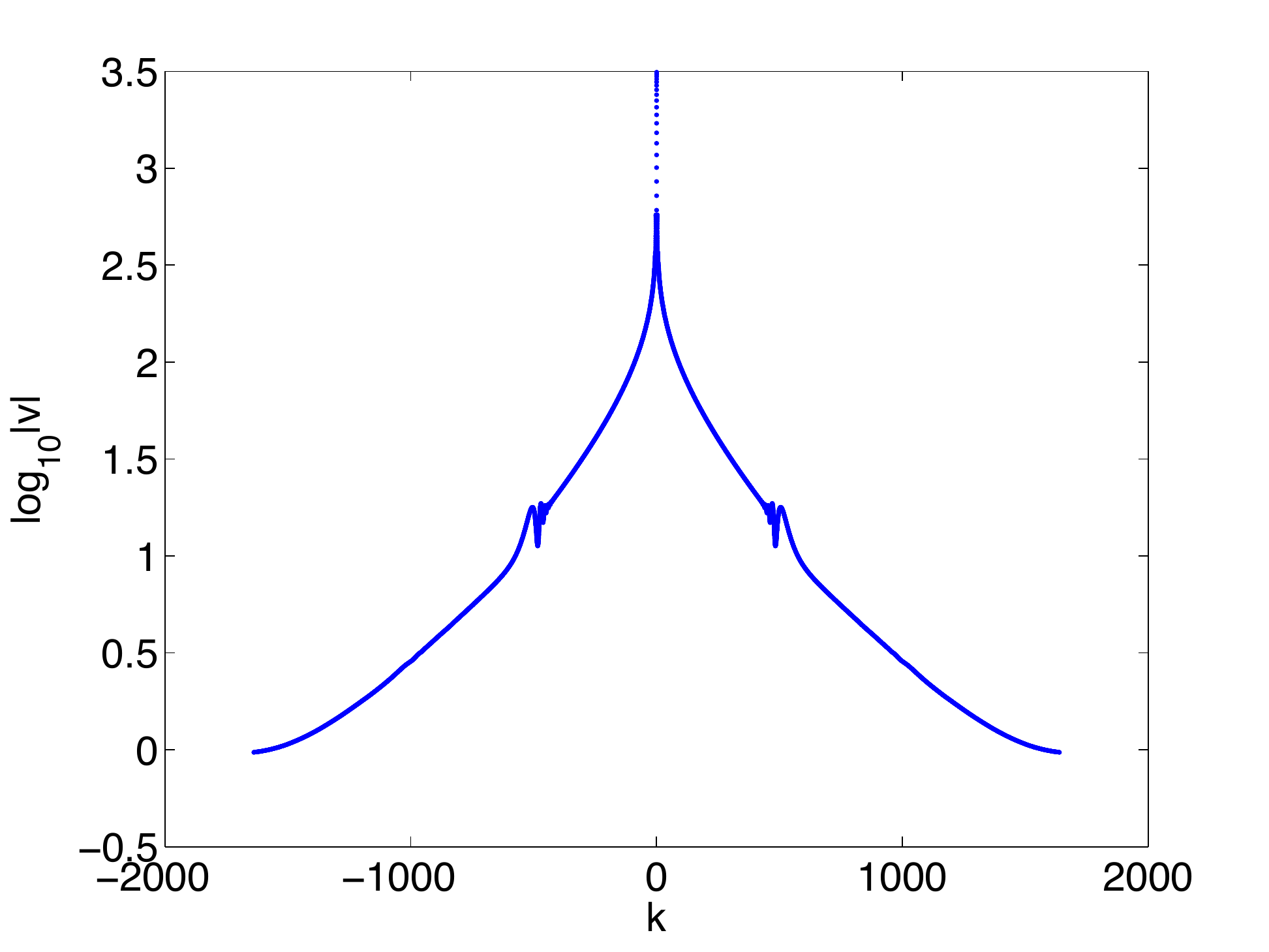}
 \caption{$L_{\infty}$ norm of the solution  
to the fKdV equation (\ref{Cauchy}) for 
 $\alpha=0.45$ and the initial data $u_{0}=3\mbox{sech}^{2}x$
 in dependence of time on the left, and the modulus of the Fourier 
 coefficients of the solution for $t=1.86$ on the right.}
 \label{gBO3sechalpha045fourier}
\end{figure}

The code is stopped 
at $t=1.86$ since the conservation of the numerically computed energy 
drops below $10^{-3}$. As can be seen in 
Fig.~\ref{gBO3sechalpha045fourier}, this is due to a lack of 
resolution in Fourier space. The $L_{\infty}$ norm of the solution in 
the same figure also indicates a blow-up.

\begin{figure}[htb!]
   \includegraphics[width=0.49\textwidth]{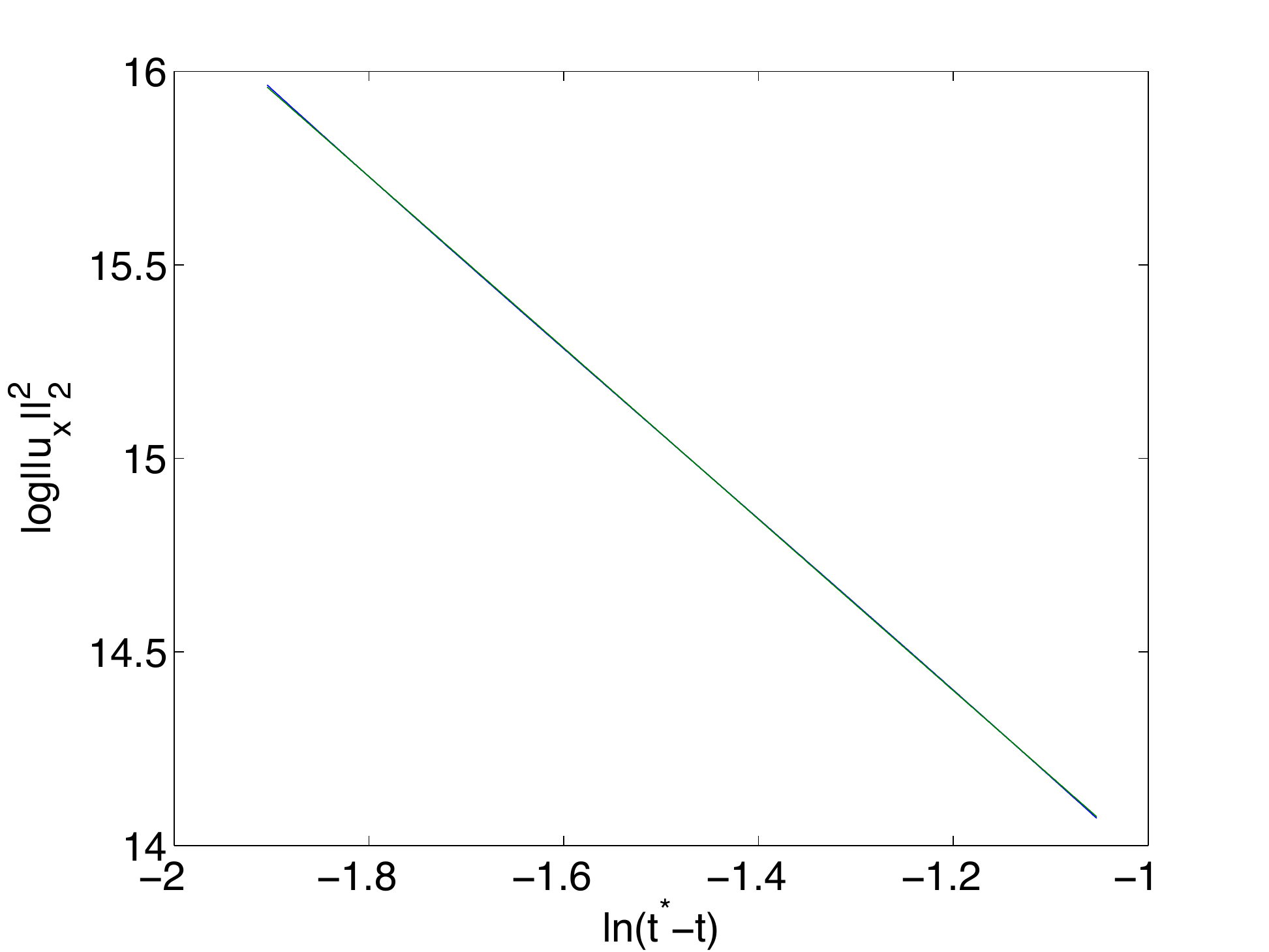}
  \includegraphics[width=0.49\textwidth]{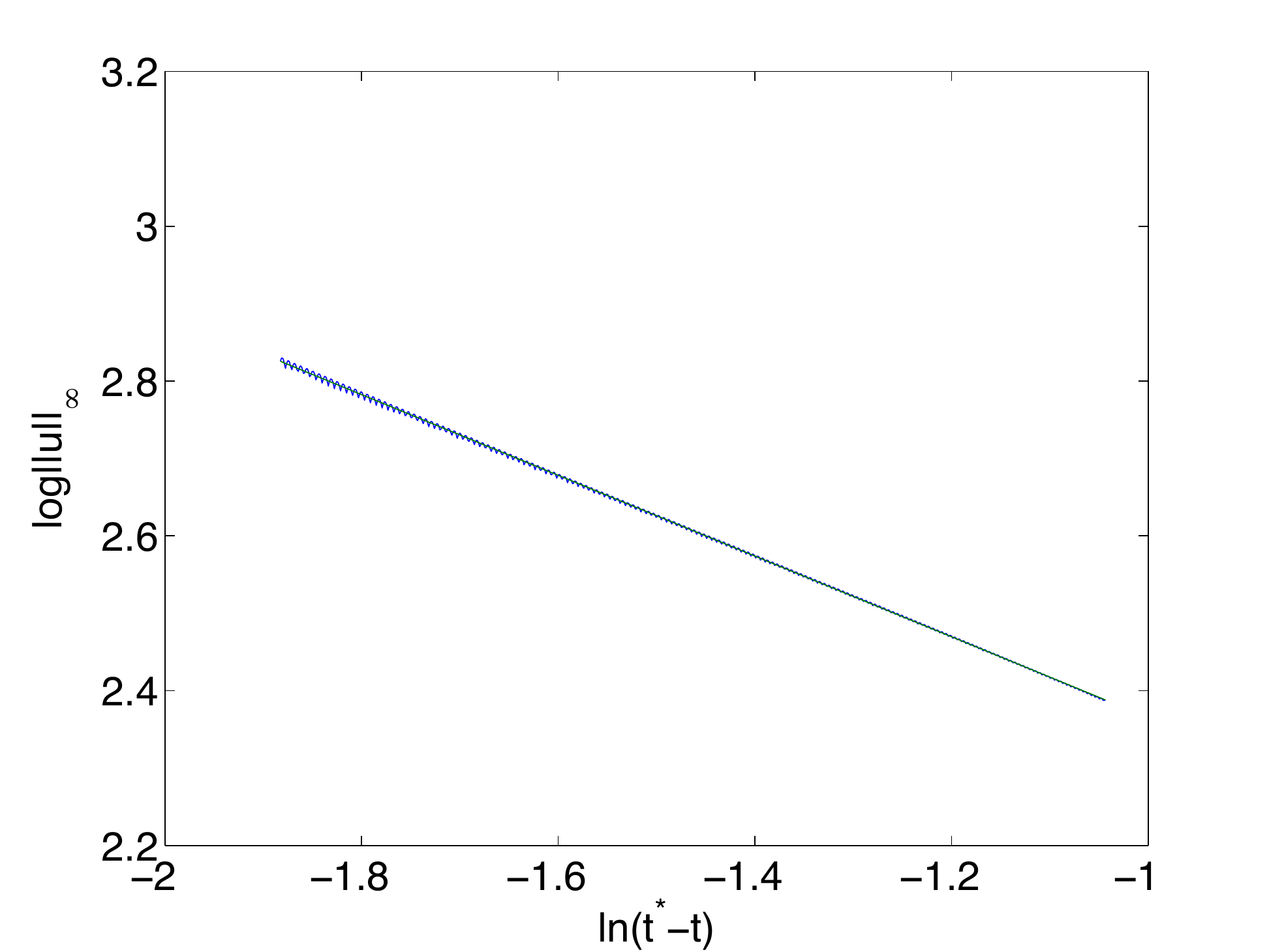}
 \caption{$L_{2}$ norm of the gradient (left) and 
 $L_{\infty}$ norm of the solution  
to the fKdV equation (\ref{Cauchy}) for 
 $\alpha=0.45$ and the initial data $u_{0}=3\mbox{sech}^{2}x$ (right) 
 for the last 1000 time steps in blue and the fitted lines $\kappa_{1}\ln 
 (t^{*}-t)+\kappa_{2}$ in green.}
 \label{gBO3sechalpha045fit}
\end{figure}

As for the $L_{2}$ critical case $\alpha=0.5$ in 
Fig.~\ref{gBO3sechalpha05fit}, we can fit various norms of the 
solution close to the blow-up to the formulae (\ref{L2scal}) and 
(\ref{genscal}) which can be seen in 
Fig.~\ref{gBO3sechalpha045fit}.  Fitting $||u_{x}||^{2}_{2}$ to  $\kappa_{1}\ln 
(t^{*}-t)+\kappa_{2}$, we find $t^{*}= 1.9489$, $\kappa_{1}=-2.2137$ 
and    $\kappa_{2}=11.7436$, and similarly we get  for $||u||_{\infty}$   
$t^{*}= 1.9522$, $\kappa_{1}=-0.5231$ and    $\kappa_{2}=1.8441$. 
Thus the fitted values for $t^{*}$ are consistent and in accordance 
with the computation. The agreement with expectation is better for 
the $L_{\infty}$ norm for 
$\kappa_{1}$ which should be according to (\ref{genscal}) $1.31$ and 
$0.31$ respectively.  

\begin{figure}[htb!]
   \includegraphics[width=0.49\textwidth]{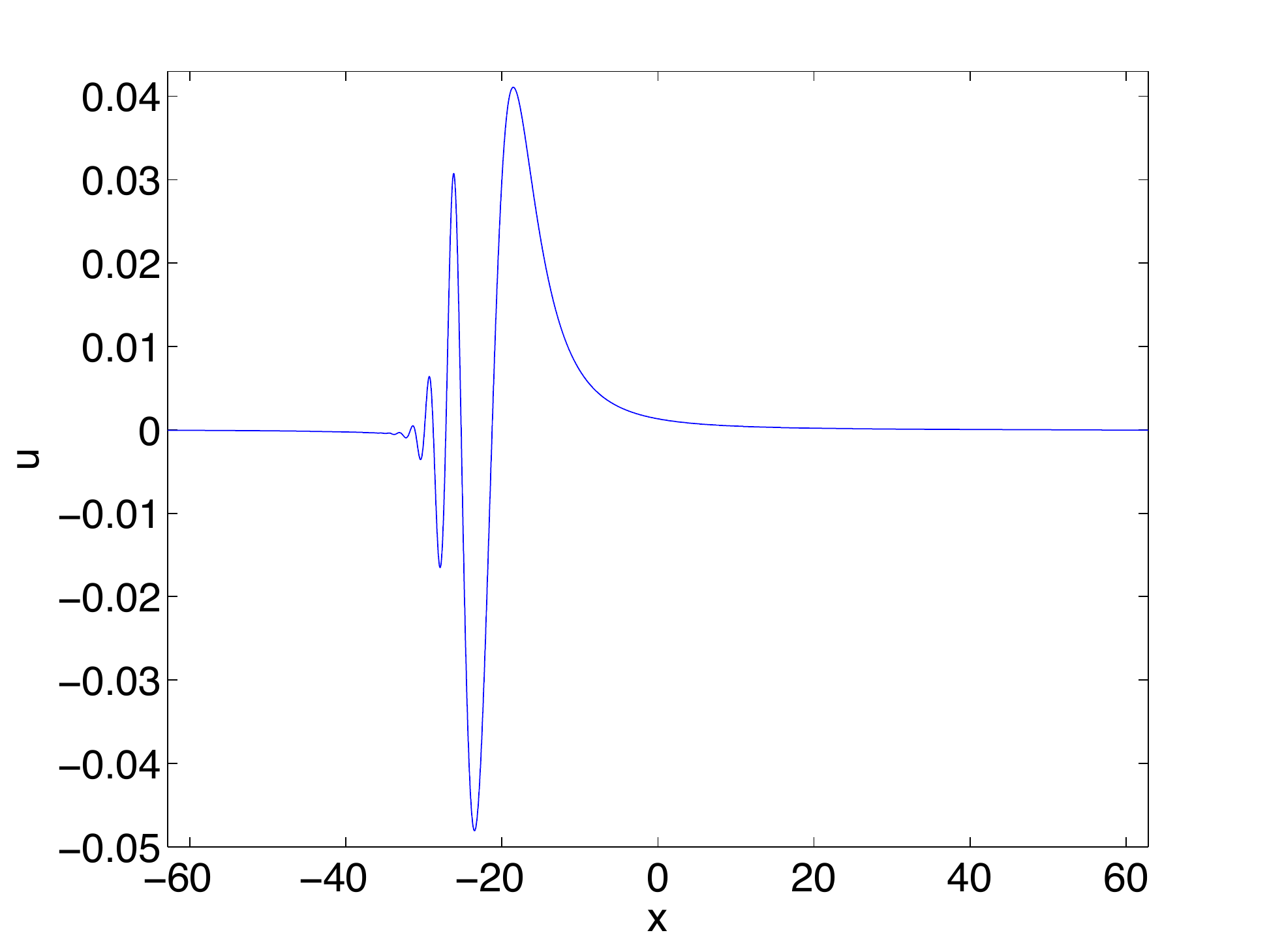}
  \includegraphics[width=0.49\textwidth]{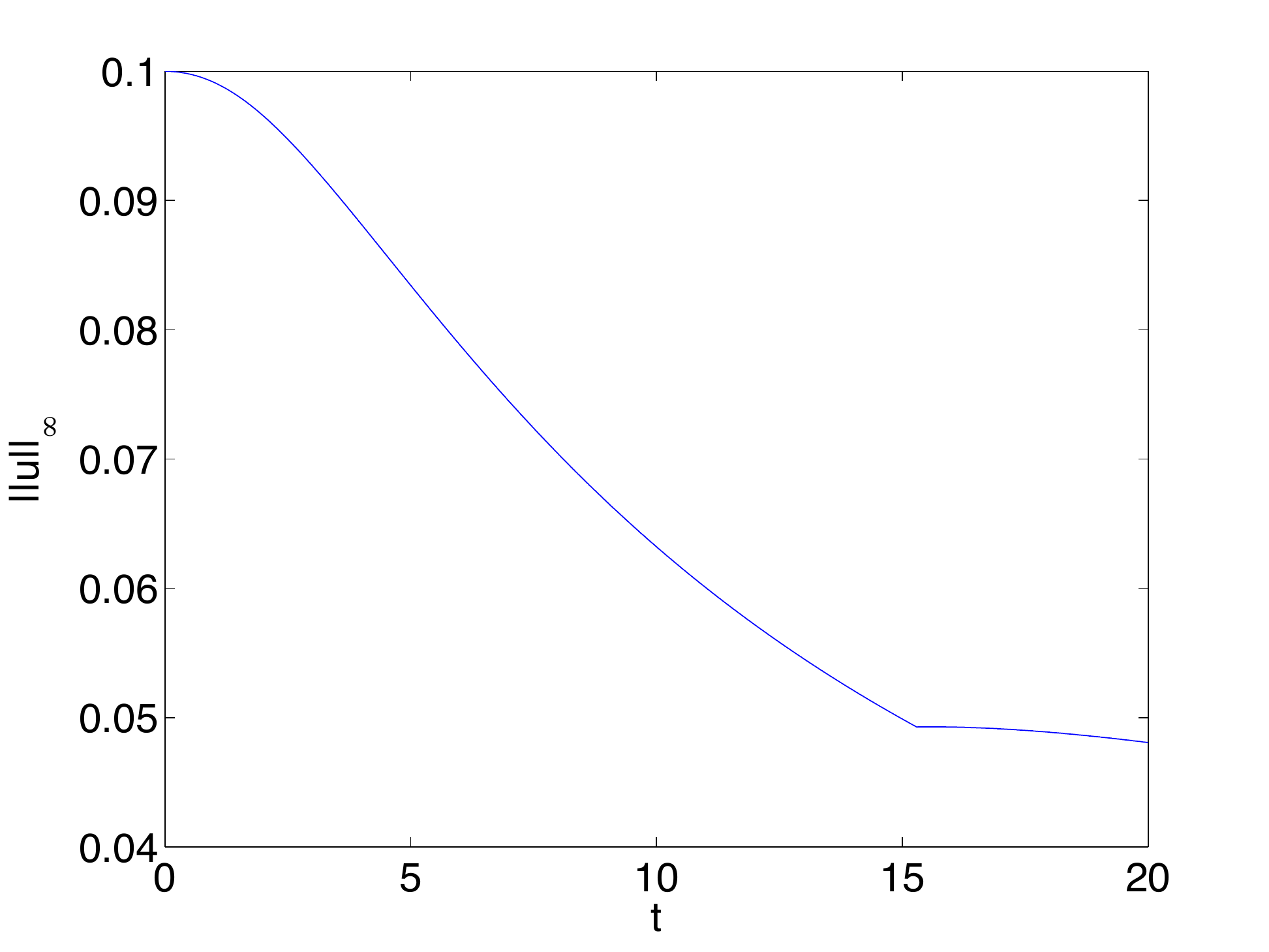}
 \caption{Solution  
to the fKdV equation (\ref{Cauchy}) for 
 $\alpha=0.2$ and the initial data $u_{0}=0.1\mbox{sech}^{2}x$ (right) 
 for $t=20$ on the left, and the $L_{\infty}$ norm of the solution in 
 dependence of $t$ on the right.}
 \label{gBO01sechalpha02}
\end{figure}

The exact criterion for the initial data to lead to a blow-up in 
finite time is not clear. For $1/3<\alpha<1/2$ it could be related to 
the mass or energy of the soliton. But since there are no solitons 
for $\alpha\leq 1/3$, the soliton for the given value of $\alpha$ 
cannot provide a criterion unless there is always blow-up.
However, this is not the case for smaller initial data  as can be seen for 
instance in Fig.~\ref{gBO01sechalpha02} where we show the solution 
to the fKdV equation for the initial data 
$u_{0}=0.1\mbox{sech}^{2}x$. The solution will be just radiated away to 
infinity for large times even if $\alpha<1/3$.  The computation is carried out with 
$N=2^{14}$ Fourier modes for $x\in 20[-\pi,\pi]$ with $N_{t}=10^{4}$ 
time steps for $t\leq 20$. It can be seen in the same figure that the 
$L_{\infty}$ norm of the solution is monotonically decreasing. Thus 
there is no indication of blow-up for initial data with small enough 
mass. Obviously there can be 
globally regular solution for all $t$ in this case, which is in strong contrast with 
the Burgers equation for which any localized initial data leads to a shock formation.

\begin{figure}[htb!]
  \includegraphics[width=\textwidth]{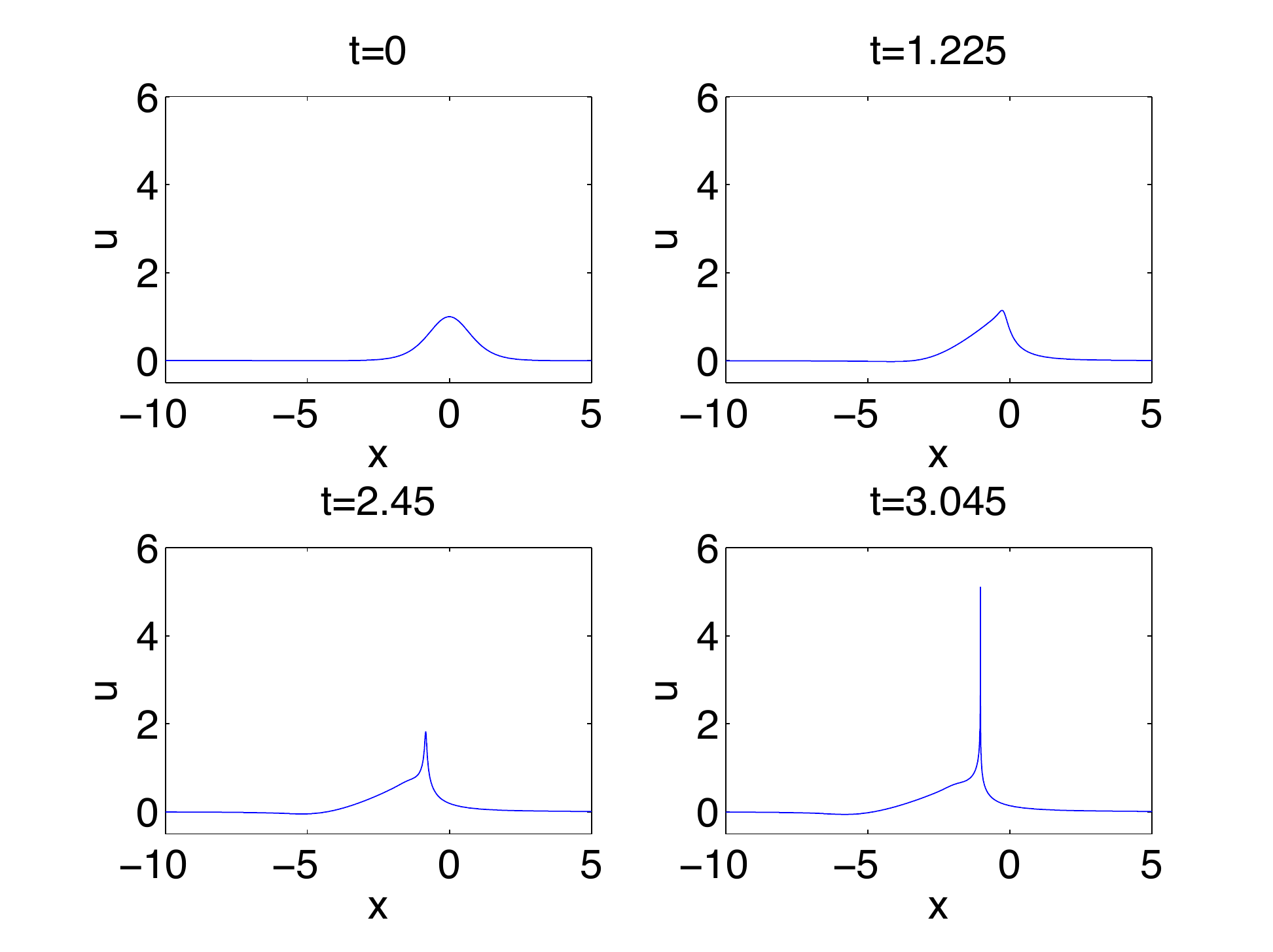}
 \caption{Solution to the fKdV equation (\ref{Cauchy}) for 
 $\alpha=0.2$ and the initial data $u_{0}=\mbox{sech}^{2}x$ 
 for several values of $t$.}
 \label{gBOsechalpha024t}
\end{figure}

The situation is different for larger data. For $\alpha=0.2$, the energy of the initial data 
$u_{0}=\mbox{sech}^{2}x$ is positive, but there seems to be a much 
more pronounced blow-up than in the blow-up cases above
as can be seen in Fig.~\ref{gBOsechalpha024t}. In this 
case there are almost no dispersive oscillations visible since the 
dispersion is much weaker, and the whole initial hump seems to be 
increasingly peaked. Obviously there can be no soliton forming which 
then blows up. But there is a rapid increase in the 
$L_{\infty}$ norm at $t=3.045$ where the code ceases to converge. 
Thus for $\alpha<1/3$, blow-up is observed also for initial data with 
positive energy. Note the difference in the type of 
blow-up in the energy supercritical case here and for $\alpha=0.45$ 
in Fig.~\ref{gBO3sechalpha045}. There a soliton forms which is 
unstable against blow-up. Here the solution follows initially the 
pattern of a solution to the Burgers' equation, a steepening almost 
leading to a shock. But instead of reaching a point of gradient 
catastrophe, the maximum of the solution appears to turn into an 
$L_{\infty}$ blow-up.

\begin{figure}[htb!]
   \includegraphics[width=0.49\textwidth]{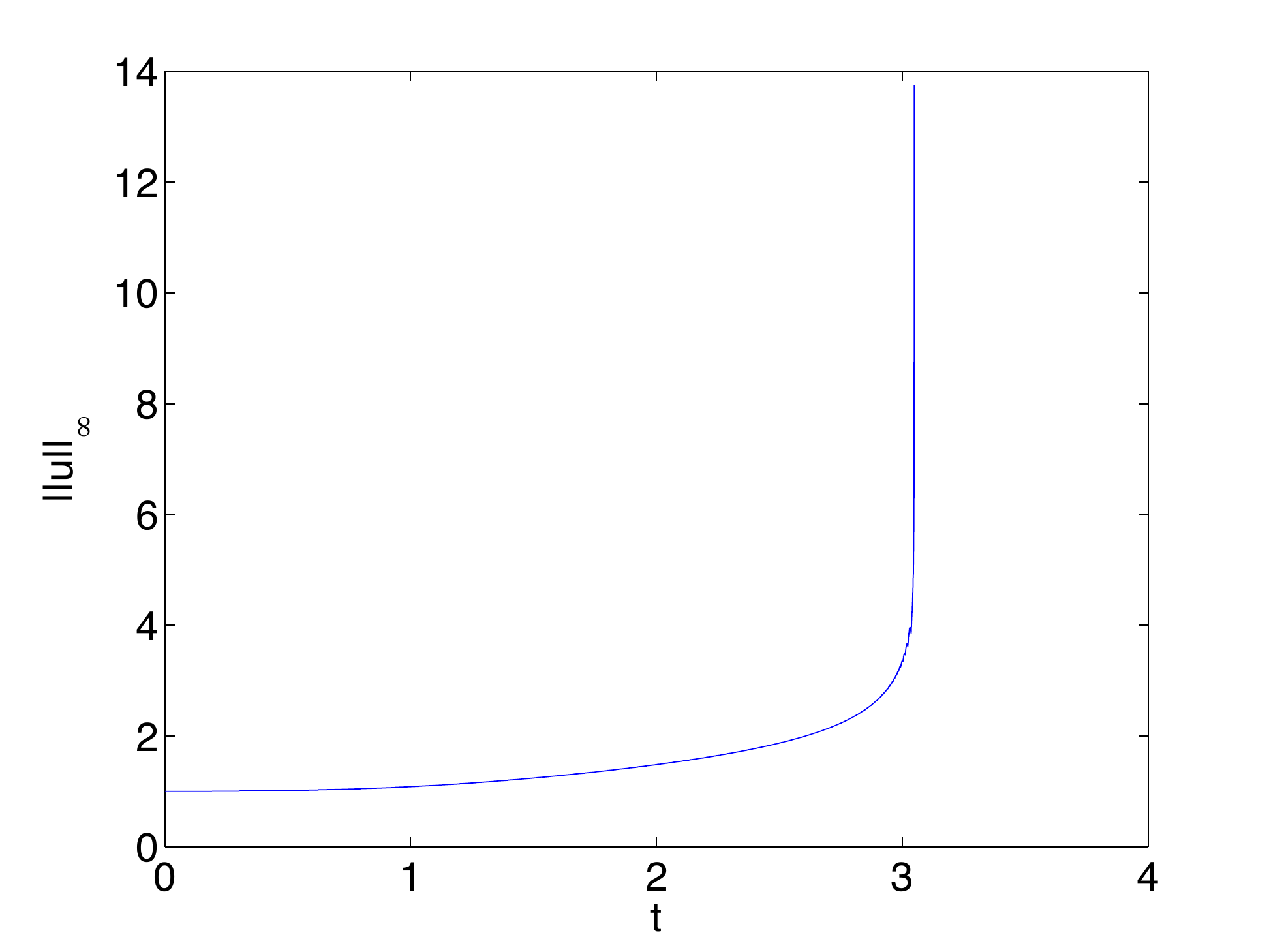}
  \includegraphics[width=0.49\textwidth]{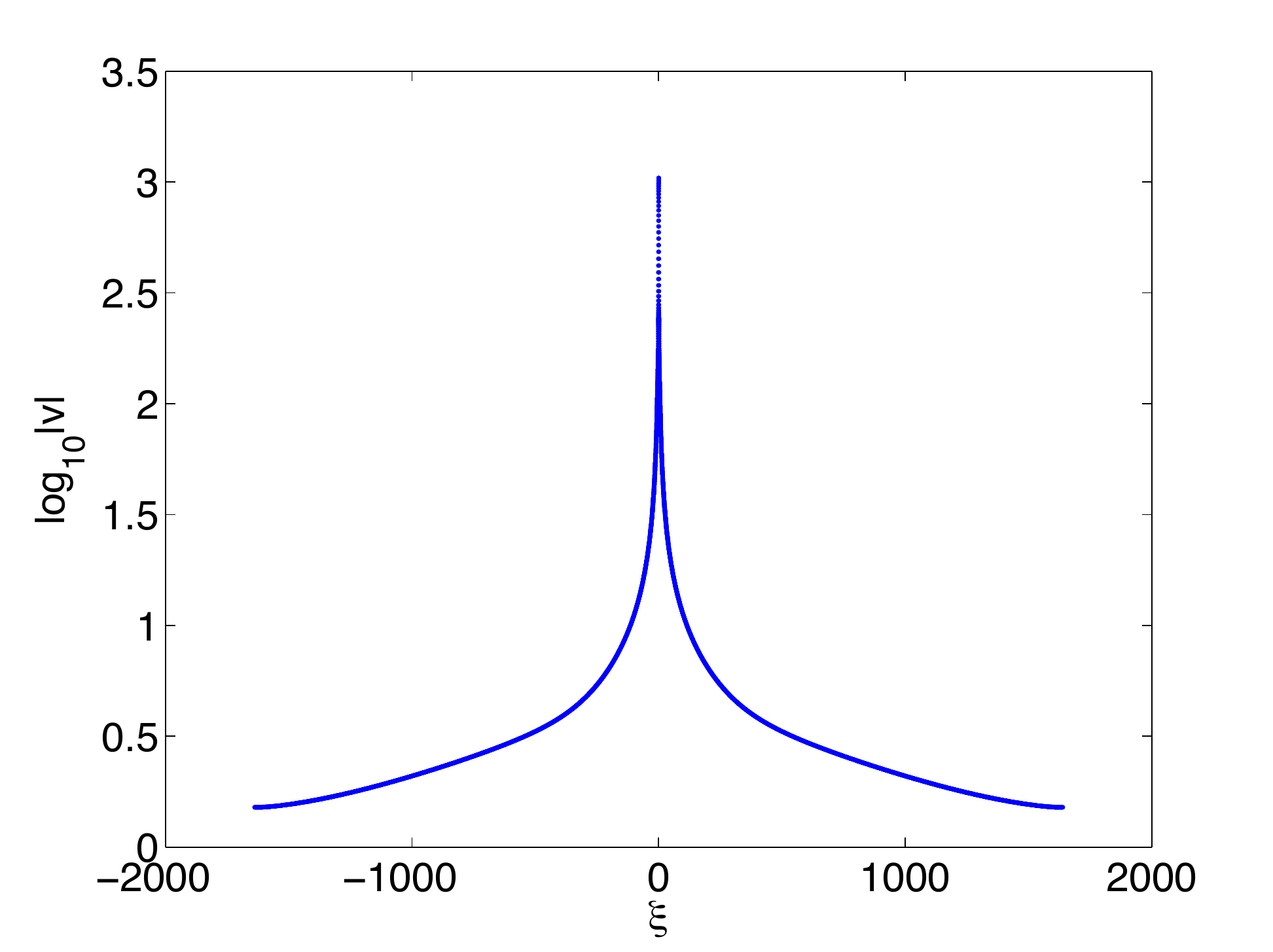}
 \caption{$L_{\infty}$ norm of the solution  
to the fKdV equation (\ref{Cauchy}) for 
 $\alpha=0.2$ and the initial data $u_{0}=\mbox{sech}^{2}x$
 in dependence of time on the left, and the modulus of the Fourier 
 coefficients of the solution for $t=3.045$ on the right.}
 \label{gBOsechalpha02}
\end{figure}

This rapid increase in the $L_{\infty}$ norm is even more visible in 
Fig.~\ref{gBOsechalpha02}. It can be also seen in this figure that 
the failure of the code to converge is due to a lack of resolution in 
Fourier space (the computed relative energy is still conserved to the 
order of $10^{-10}$). 

\begin{figure}[htb!]
   \includegraphics[width=0.49\textwidth]{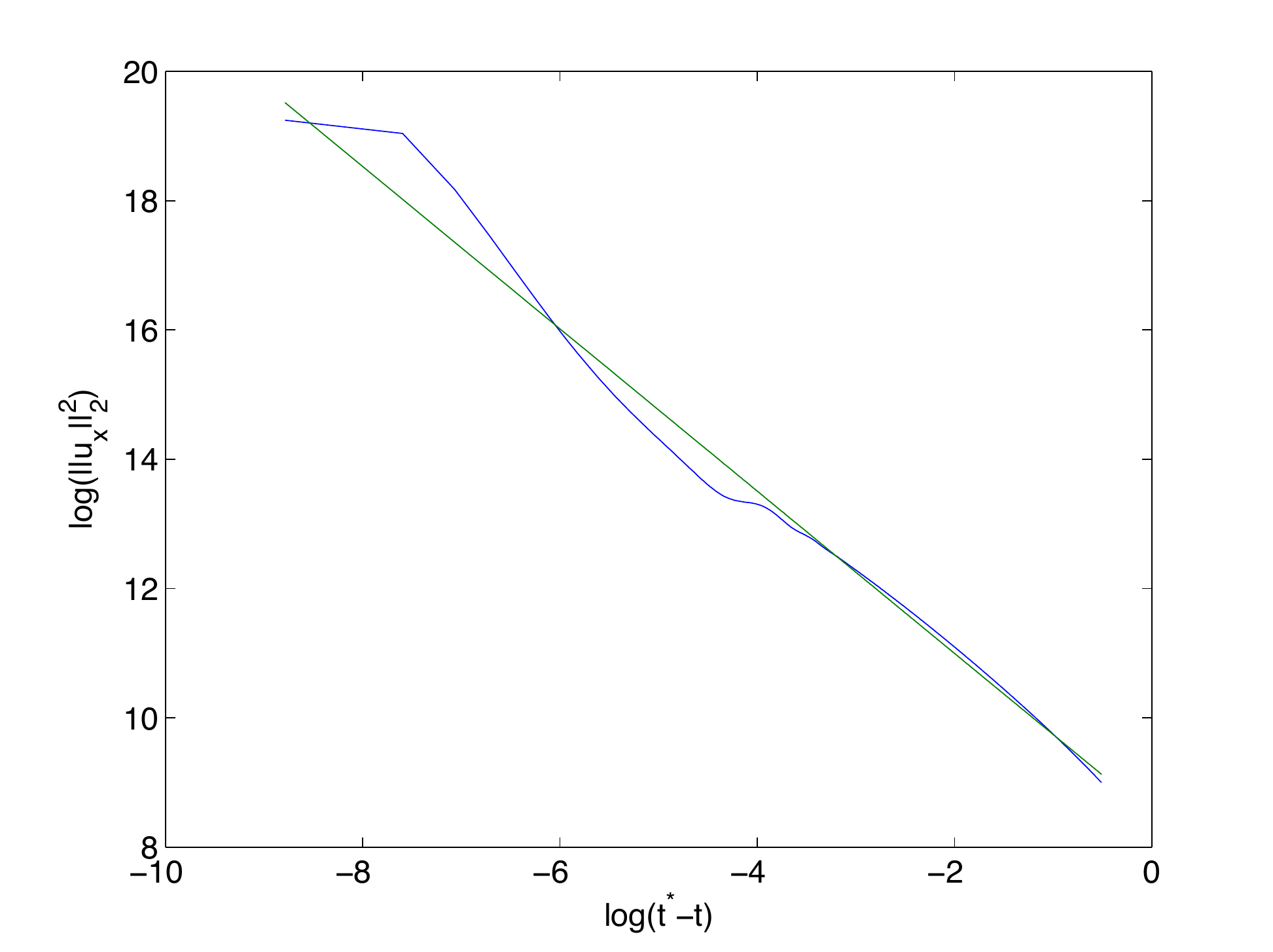}
  \includegraphics[width=0.49\textwidth]{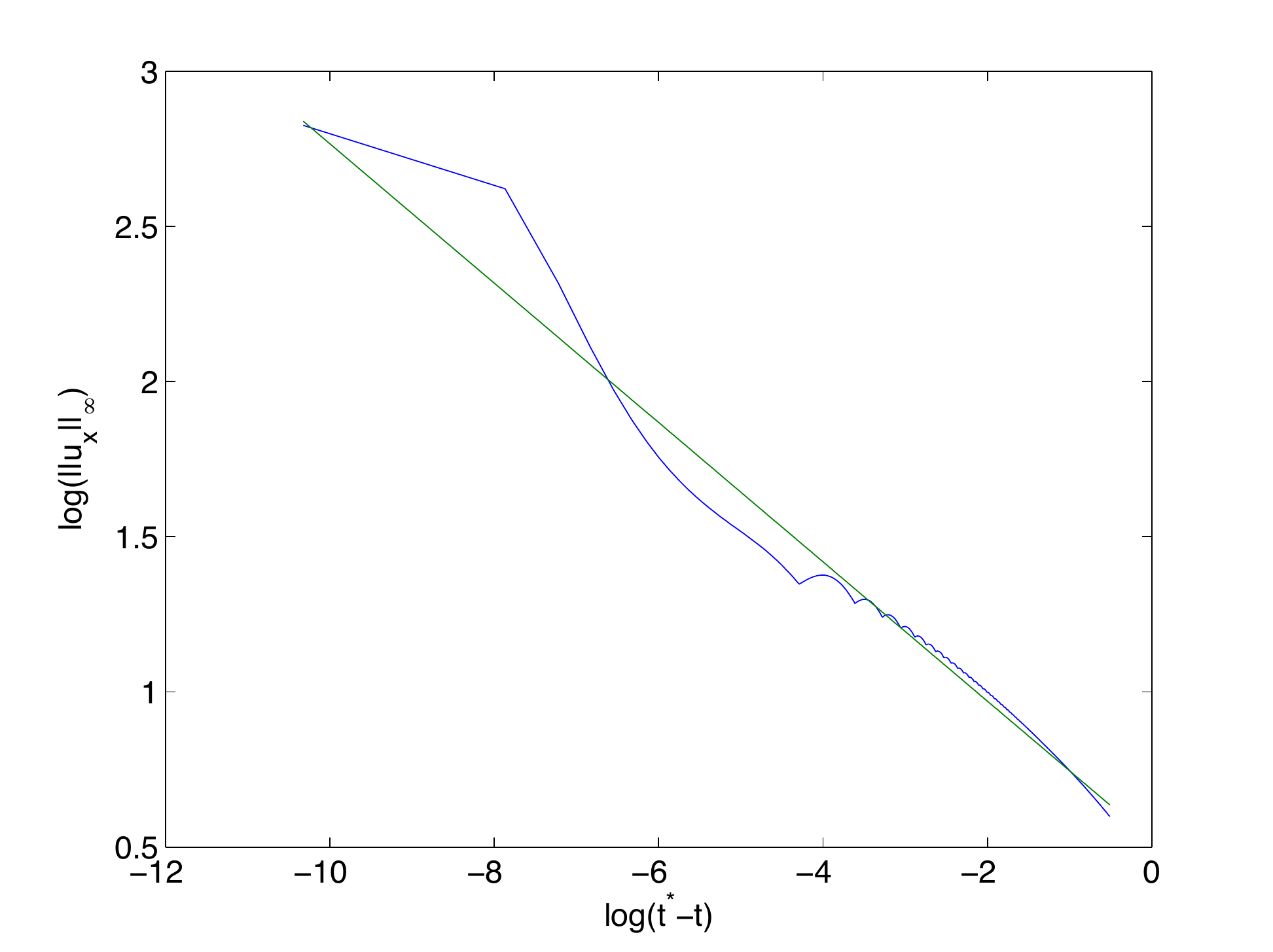}
 \caption{$L_{2}$ norm of the gradient (left) and 
 $L_{\infty}$ norm of the solution  
to the fKdV equation (\ref{Cauchy}) for 
 $\alpha=0.2$ and the initial data $u_{0}=\mbox{sech}^{2}x$ (right) 
 for $t>2.4497$ in blue and the fitted lines $\kappa_{1}\ln 
 (t^{*}-t)+\kappa_{2}$ in green.}
 \label{gBOsechalpha02fit}
\end{figure}

As for the above blow-up cases, we can fit various norms of the 
solution close to the blow-up to the formulae (\ref{L2scal}) and 
(\ref{genscal}). The rapid divergence of the norms indicates already 
an exponential (in the rescaled time $\tau$) blow-up, and this is 
confirmed by the fitting procedure. However, such an exponential 
behavior implies a less good fitting since it will be dominant only 
very close to the blow-up where the code is no longer fully reliable. 
Thus we fit on a larger interval, $2.4497<t<3.045$ and note that the 
fitting parameters get closer to the theoretically expected ones if 
the lower bound is increased. It can be seen in 
Fig.~\ref{gBOsechalpha02fit} that the fitting is very sensitive to 
the fitted time $t^{*}$, and that it is problematic for $t\sim 
t^{*}$.  Fitting $||u_{x}||^{2}_{2}$ to  $\kappa_{1}\ln 
(t^{*}-t)+\kappa_{2}$, we find $t^{*}= 3.0490$, $\kappa_{1}=-1.2552$ 
and    $\kappa_{2}=8.4845$, and similarly we get  for $||u||_{\infty}$   
$t^{*}= 3.0489$, $\kappa_{1}=-0.2245$ and    $\kappa_{2}=0.5210$. 
Thus the fitted values for $t^{*}$ are consistent and in accordance 
with the computation. This also applies within numerical precision to
the values of 
$\kappa_{1}$ which should be according to (\ref{genscal}) $7/6$ and 
$1/6$ respectively.  

\begin{remark}
    In the cases with blow-up, we can test whether the Fourier 
    coefficients can be fitted to the asymptotic formula 
    (\ref{fourasymp}) which should indicate the formation of a 
    singularity on the real axis via the vanishing of the parameter 
    $\delta$. For the example of an $L_{2}$ critical blow-up in 
    Fig.~\ref{gBO3sechalpha05} we get $\delta=0.007$ and 
    $\mu+1=0.36$. This indicates in accordance with the fitting of 
    the norms in Fig.~\ref{gBO3sechalpha05solfit} that the blow-up is 
    not fully reached. The exact value of $\mu$ is probably not 
    reliable (it is less accurate since the fitting to an exponential 
    is less prone to errors than to a linear dependence), 
    but its sign indicates that the singularity is a pole, 
    in accordance with what was found for the $L_{\infty}$ norm of 
    the solution. For the $L_{2}$ supercritical blow-up in 
    Fig.~\ref{gBO3sechalpha045}, we find $\mu+1=0.6535$ and 
    $\delta=0.0023$, where $\mu$ is once more less reliable. 
    This indicates that we are shortly before the 
    blow-up, and that this is again an $L_{\infty}$ blow-up. For the 
    energy supercritical blow-up in 
    Fig.~\ref{gBOsechalpha02} we find $\delta=-0.003$ and 
    $\mu+1=0.8077$. This indicates in accordance with the fitting of 
    the norms in Fig.~\ref{gBOsechalpha02fit} that blow-up is indeed 
    reached. The divergence of the solution appears to be 
    proportional to $|x-x^{*}|^{-\alpha}$.
\end{remark}

We summarize the numerical findings in this subsection in the following 
\begin{conj}
    Consider  smooth initial data $u_{0}\in L_{2}(\mathbb{R})$ 
with a single hump. Then for
    \begin{itemize}
\item $\alpha>0.5$: solutions to the fKdV equations with the initial 
data $u_{0}$ stay smooth for 
all $t$. For large $t$ they decompose asymptotically into solitons and 
radiation.

\item $0<\alpha\leq 0.5$: solutions to the fKdV equations with initial 
data $u_{0}$ of sufficiently small, but non-zero mass stay smooth for 
all $t$.

\item $\alpha=0.5$: solutions to the fKdV equations  with the initial 
data $u_{0}$ with negative energy and mass larger than the soliton 
mass
blow up at finite time $t^{*}$ and infinite $x^{*}$. The type 
of the blow-up for $t\nearrow t^*$ is characterized by 
\begin{equation}
    u(x,t)\sim \frac{1}{\sqrt{L(t)}}Q_{1}\left(\frac{x-x_{m}}{L(t)}\right), \quad
    L = c_{0}(t^{*}-t)
    \label{uL2},
\end{equation}
where $c_{0}$ is a constant, and where $Q_{1}$ is the solitary wave 
solution (\ref{solitarywave}) for $c=1$. In addition one has
\begin{equation}
    ||u_{x}||_{2}\sim \frac{1}{L^{2}(t)}
    \label{UxL2}.
\end{equation}
\item $1/3<\alpha<0.5$:  solutions to the fKdV equations  with the initial 
data $u_{0}$ and sufficiently large $L_{2}$ norm
blow up at finite time $t^{*}$ and finite $x=x^{*}$. A soliton-type 
hump separates from the initial hump and eventually blows up. The type 
of the blow-up for $t\nearrow t^*$ is characterized by 
\begin{equation}
    u(x,t)\sim \frac{1}{L^{\alpha}(t)}U\left(\frac{x-x_{m}}{L(t)}\right), \quad
    L = c_{1}(t^{*}-t)^{\frac{1}{1+\alpha}}
    \label{ualpha},
\end{equation}
where $c_{1}$ is a constant, and where $U$ is a
solution of equation (\ref{sautalpharinfty}) vanishing for $|y|\to\infty$ (if 
such a solution exists). In addition one has
\begin{equation}
    ||u_{x}||_{2}\sim \frac{1}{L^{2\alpha+1}(t)}
    \label{Uxalpha}.
\end{equation}

\item $0<\alpha<1/3$: solutions to the fKdV equations  with the initial 
data $u_{0}$ and sufficiently large $L_{2}$ norm
blow up at finite time $t^{*}$ and finite $x=x^{*}$. 
The nature of blow-up
is different from the previous one since  no solitary waves exist in 
this case, the maximum of the initial hump evolves directly into a 
blow-up.
Thus the blow-up seems to be different from that occurring 
in the supercritical gKdV equation
\eqref{gKdv} when $p>4.$ But the blow-up profile appears to be still 
given by (\ref{ualpha}).

 \end{itemize}
\end{conj}

\subsection{Numerical study of fractionary BBM equations}
In this subsection we study the behavior of solutions to the fractionary 
BBM (fBBM) equation (\ref{fracBBM}). We consider as before initial data of 
the form $u_{0}=\beta \mbox{sech}^{2}x$ with $\beta$ a real constant. 
We find that for $\alpha>1/3$, initial data of sufficient mass asymptotically 
decompose into solitons. For $\alpha\leq1/3$, such data lead to a 
blow-up in the form of a cusp at which the $L_{\infty}$ norm of the 
solution stays finite. 

\begin{figure}[htb!]
  \includegraphics[width=\textwidth]{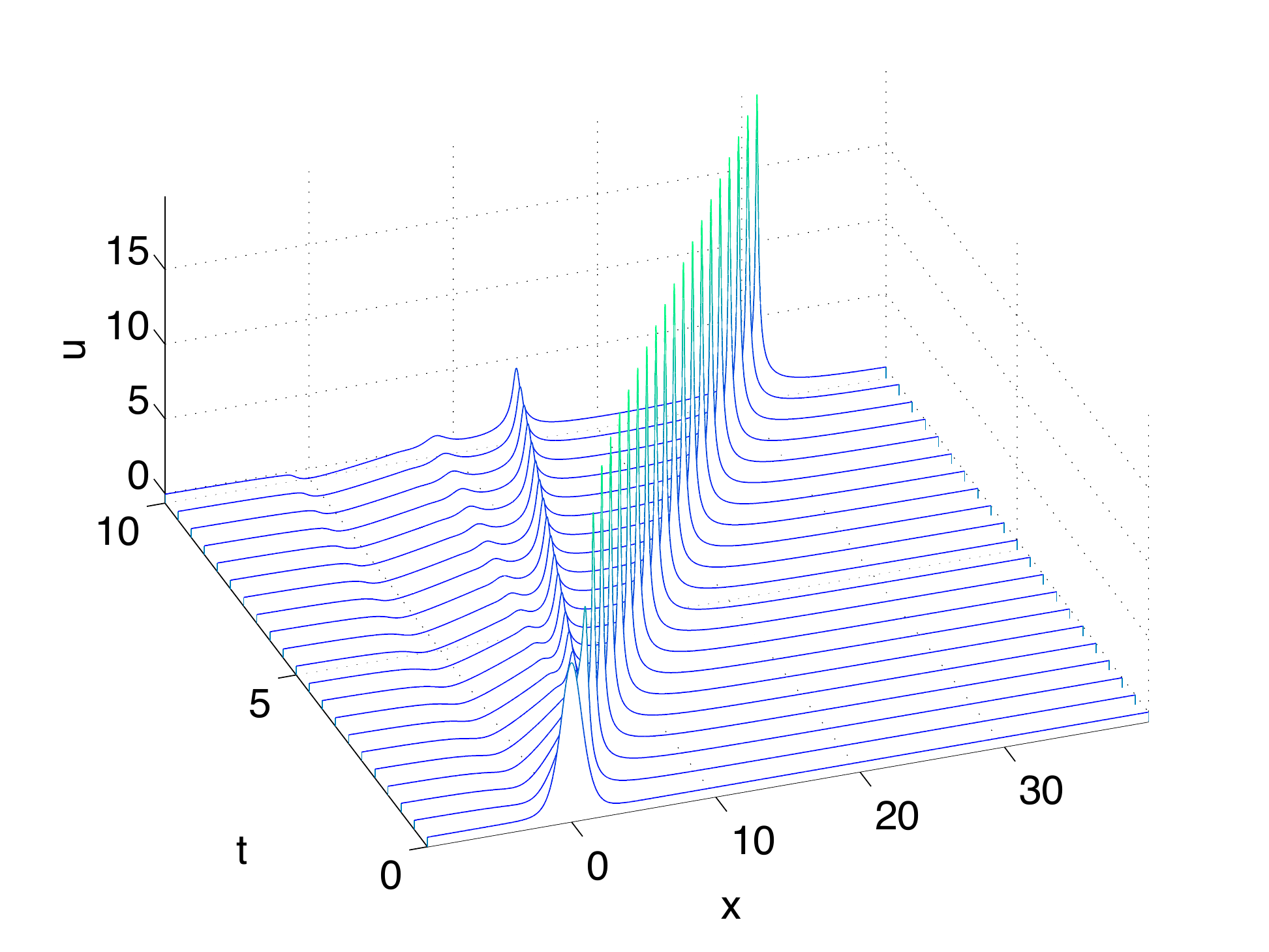}
 \caption{Solution to the fBBM equation (\ref{fracBBM}) for 
 $\alpha=0.5$ and the initial data $u_{0}=10\mbox{sech}^{2}x$ 
 in dependence  of $x$ and $t$.}
 \label{gBBM10sechalpha05water}
\end{figure}

We find a similar behavior of the solutions for $1/3<\alpha\leq1$. 
For $\alpha=0.5$ and $\beta=10$ we use $N=2^{14}$ Fourier modes for 
$x\in20[-\pi,\pi]$ and  $N_{t}=2*10^{4}$ time steps for $t\leq 10$ to obtain the 
solution shown in Fig.~\ref{gBBM10sechalpha05water}. It can be seen that 
the initial hump splits into several more strongly peaked humps which 
move like solitons, i.e., without changing their shape with 
essentially constant speed. This is similar to what was observed for 
fKdV solutions for $\alpha>1/2$. For large $t$ the solutions appear 
to be given by solitons and radiation.

\begin{figure}[htb!]
   \includegraphics[width=0.49\textwidth]{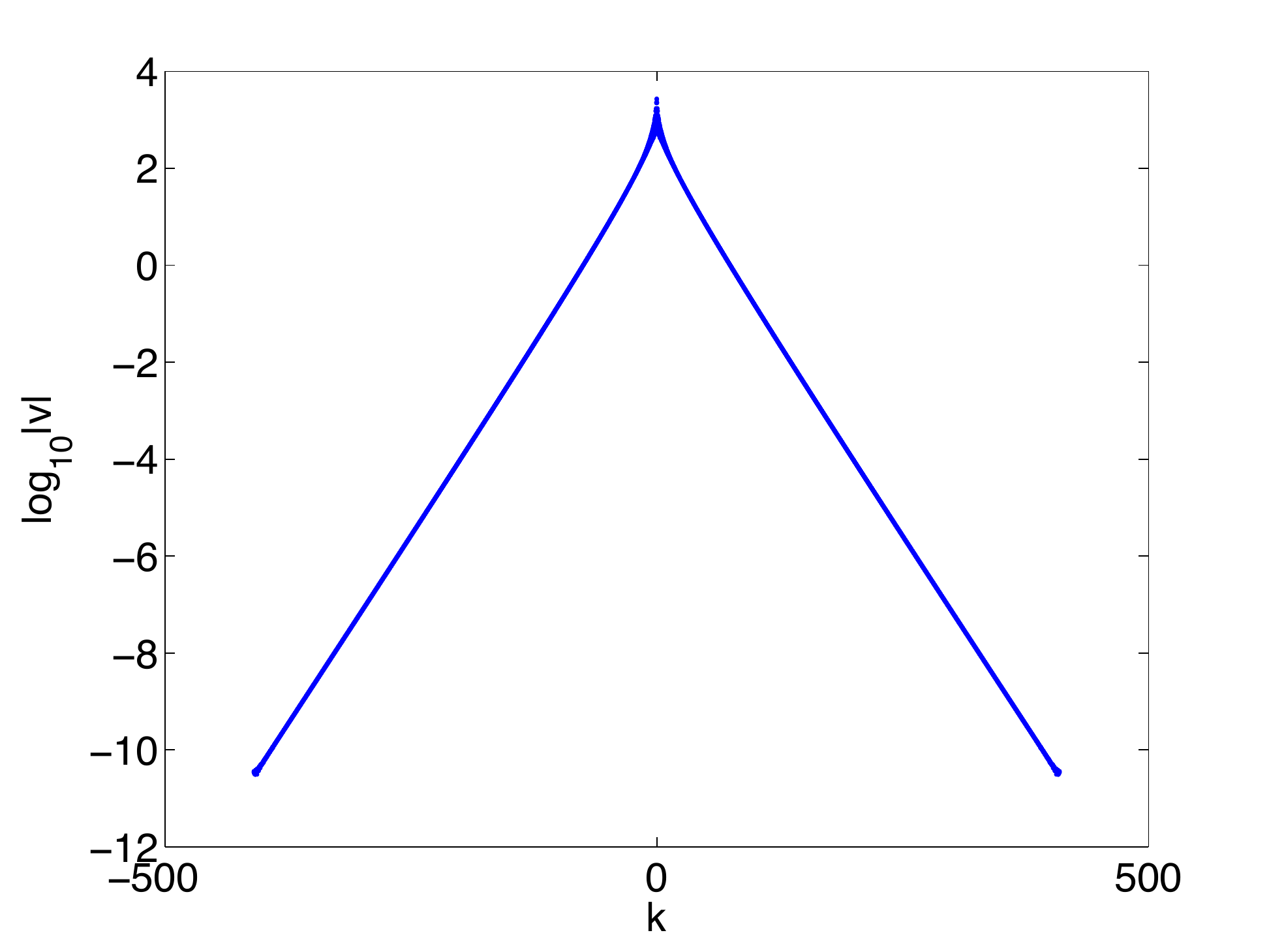}
  \includegraphics[width=0.49\textwidth]{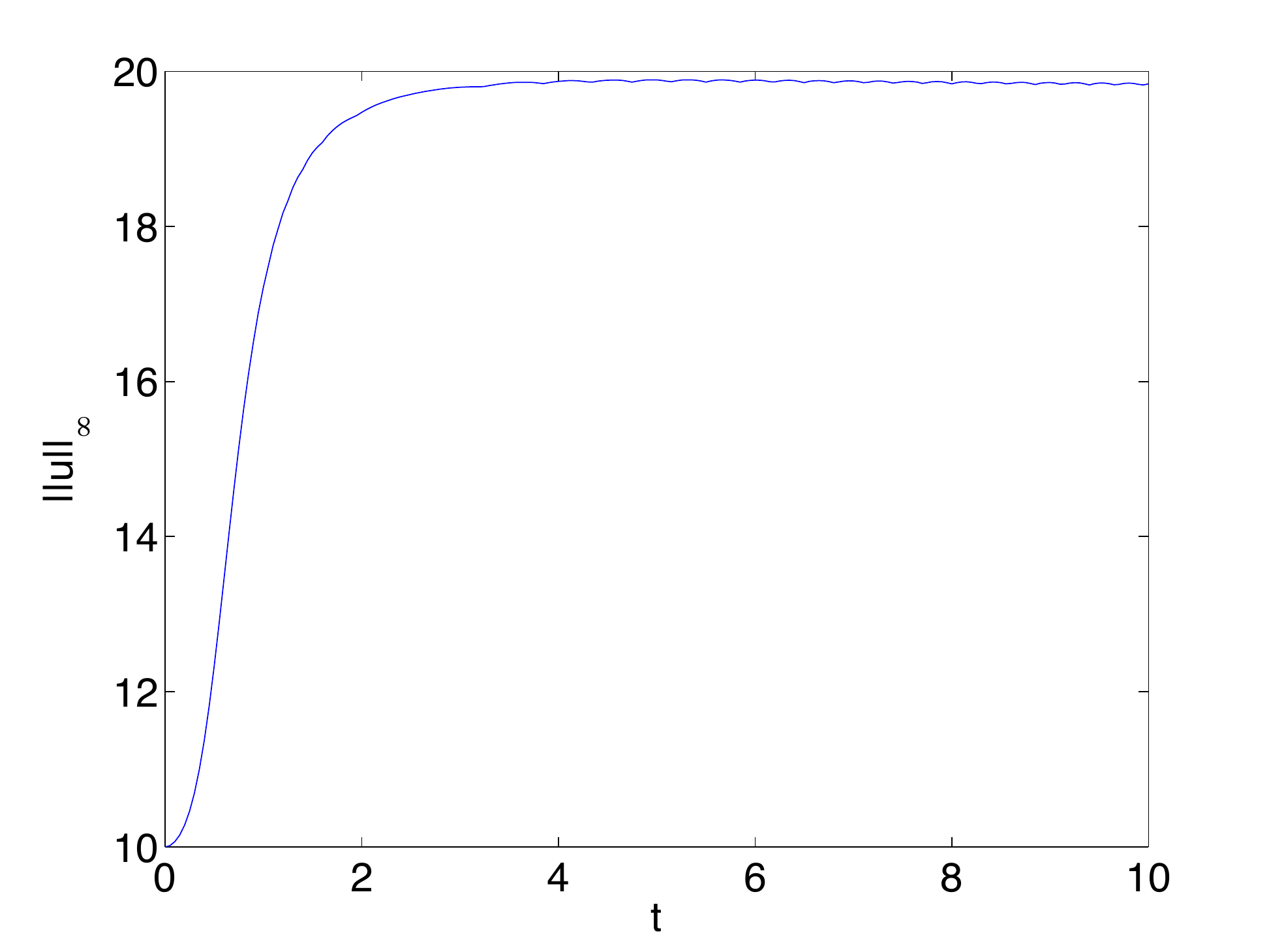}
 \caption{Solution to the fBBM equation (\ref{fracBBM}) for 
 $\alpha=0.5$ for the initial data $u_{0}=10\mbox{sech}^{2}x$  
 at $t=10$; on the left the Fourier coefficients, on the right the 
 $L_{\infty}$ norm in dependence of time.}
 \label{gBBM10sechalpha05fourier}
\end{figure}

The computation is carried out with a relative conservation of the computed 
energy of the order of $10^{-12}$. In 
Fig.~\ref{gBBM10sechalpha05fourier} it can be seen that the solution 
is well resolved in Fourier space. The $L_{\infty}$ norm of the 
solution confirms the `solitonic' appearance. After an initial 
increase of the norm, it stays essentially constant. The small 
oscillations are due to radiation propagating to the left which 
reenters the computational domain on the right because of the 
periodic boundary conditions.

\begin{figure}[htb!]
  \includegraphics[width=0.7\textwidth]{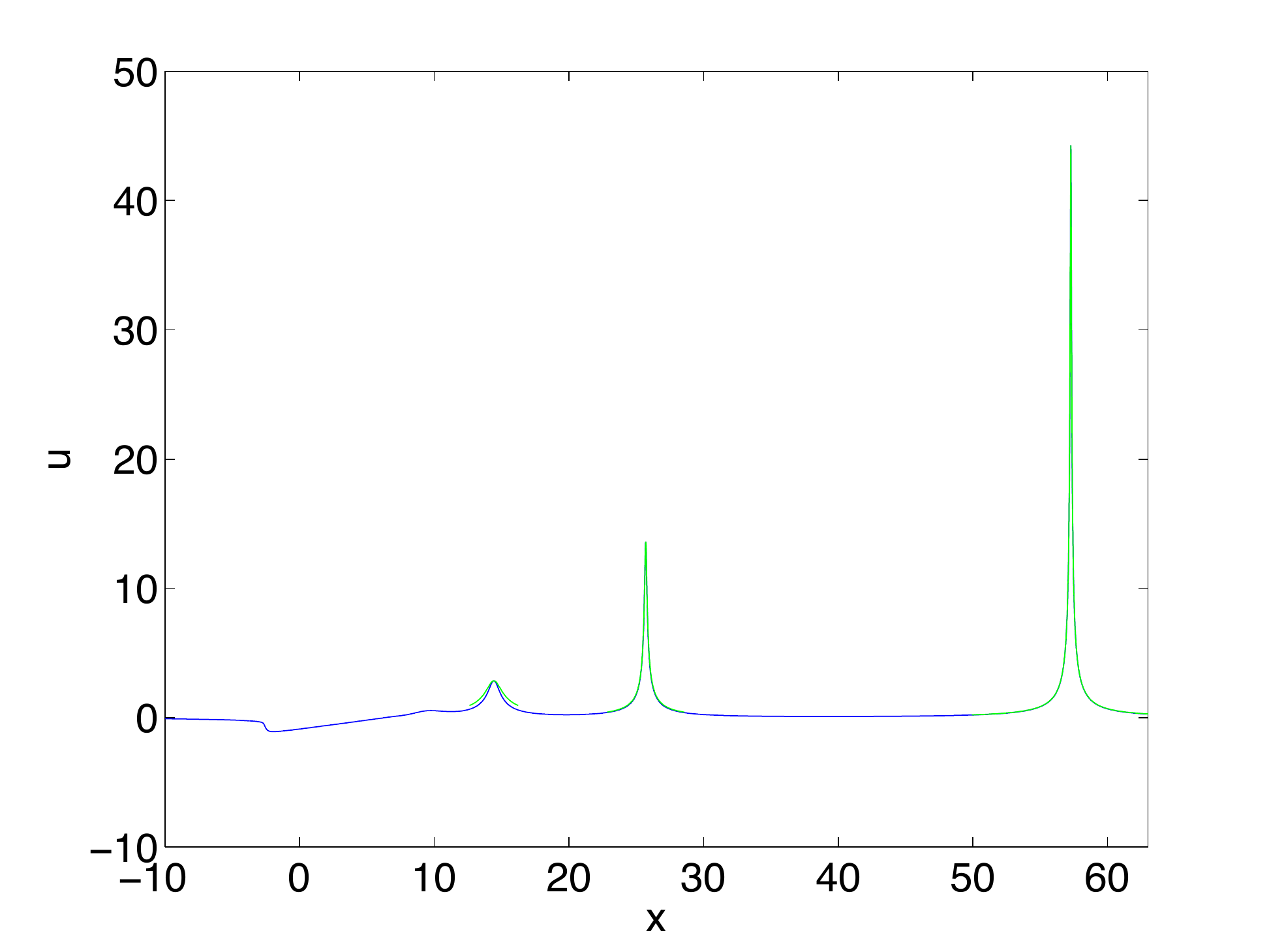}
 \caption{Solution to the fBBM equation (\ref{fracBBM}) for 
 $\alpha=0.5$ and the initial data $u_{0}=20\mbox{sech}^{2}x$  for 
 $t=10$ in blue, and fitted solitons (\ref{solitarywave}) rescaled 
 according to (\ref{solresc}) and (\ref{QfBBM}) in green.}
 \label{gBBM20sechalpha05t10}
\end{figure}

There is no indication of blow-up in this case, the solutions appears 
to be globally smooth in time. This does not change for initial data 
with larger $L_{2}$ norm and energy as can be seen in 
Fig.~\ref{gBBM20sechalpha05t10}, where the solution for 
$u_{0}=20\mbox{sech}^{2}x$ is shown at $t=10$. In this case the humps into 
which the initial data decompose have a larger maximum and propagate 
faster, but there is again no indication of blow-up. We show in  
Fig.~\ref{gBBM20sechalpha05t10} the fitting of the humps to rescaled 
(according to (\ref{solresc}) and (\ref{QfBBM}))
solitons. It can be seen that the largest and left most hump is 
almost indistinguishable from the soliton, whereas the smallest 
fitted hump is not yet in the asymptotic regime. Nonetheless the good 
agreement of the fitting to solitons indicates that the solution 
asymptotically for large $t$ decomposes into solitons and radiation. 
Note that this implies that the solitons are stable also for 
$\alpha\leq1/2$ for fBBM in contrast to fKdV. 

\begin{figure}[htb!]
  \includegraphics[width=\textwidth]{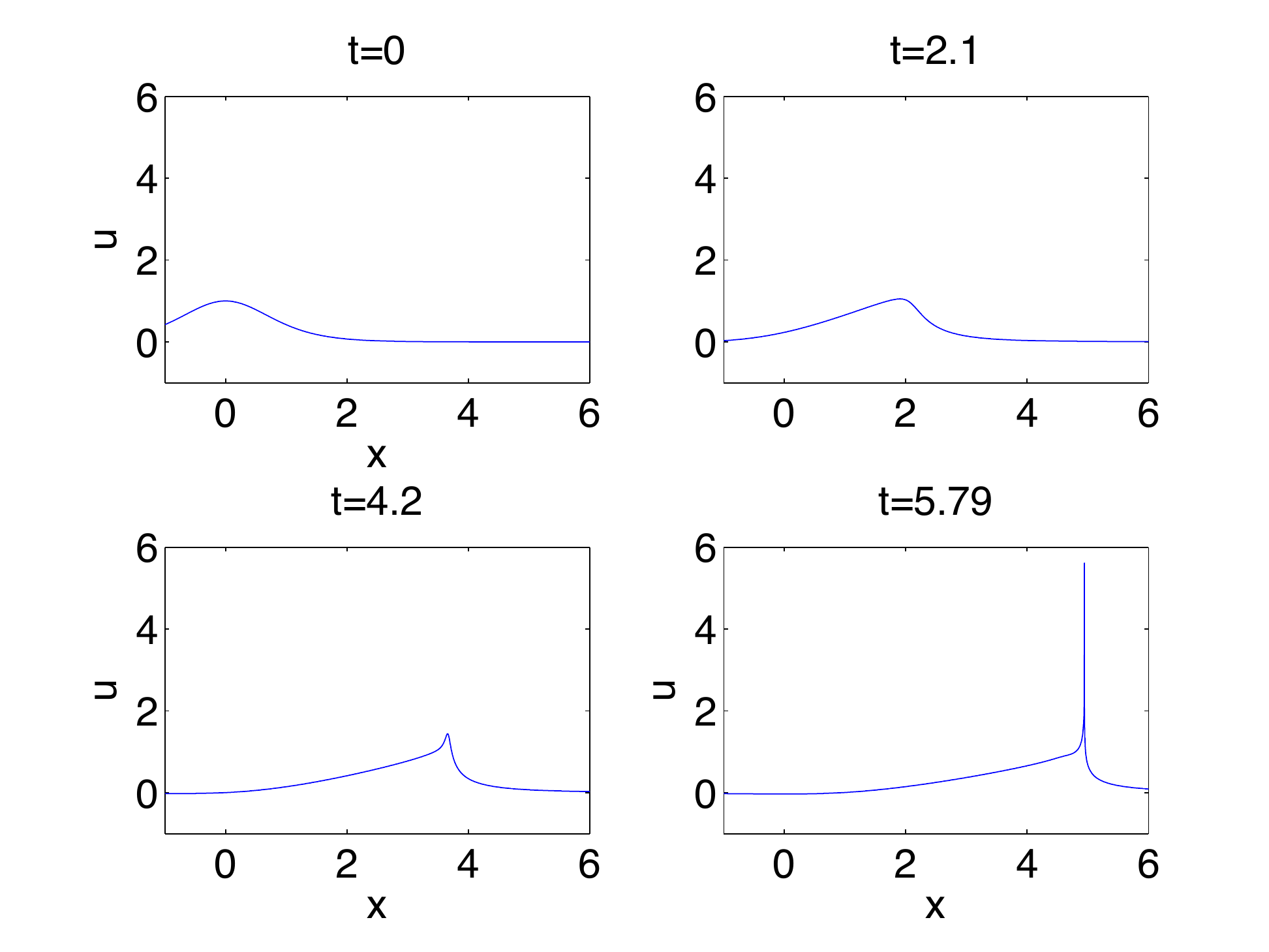}
 \caption{Solution to the fBBM equation (\ref{fracBBM}) for 
 $\alpha=0.2$ and the initial data $u_{0}=\mbox{sech}^{2}x$ 
 for several values of $t$.}
 \label{gBBMsechalpha024t}
\end{figure}

The situation is less clear for $\alpha<1/3$. We study the solution 
for $\alpha=0.2$ and the initial data $u_{0}=\mbox{sech}^{2}x$ with 
$N=2^{16}$ modes for $x\in3[-\pi,\pi]$ with $N_{t}=2*10^{4}$ time 
steps for $t\leq6$. At $t=5.85$ the relative computed energy drops 
below $10^{-3}$ at which point the code is stopped since the solution 
is no longer reliable.  In Fig.~\ref{gBBMsechalpha024t} we show the 
solution for several values of $t$. 

\begin{figure}[htb!]
   \includegraphics[width=0.49\textwidth]{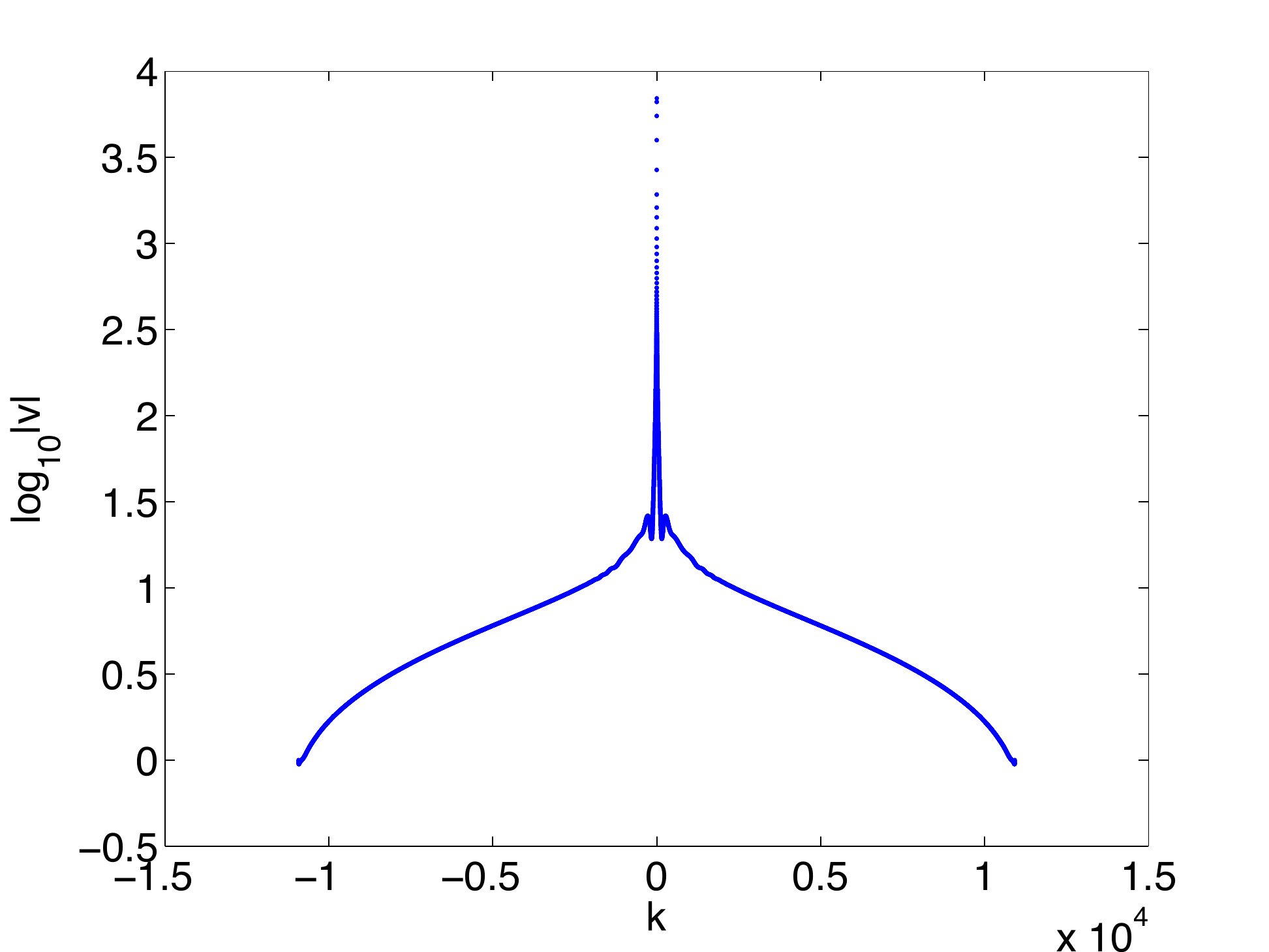}
  \includegraphics[width=0.49\textwidth]{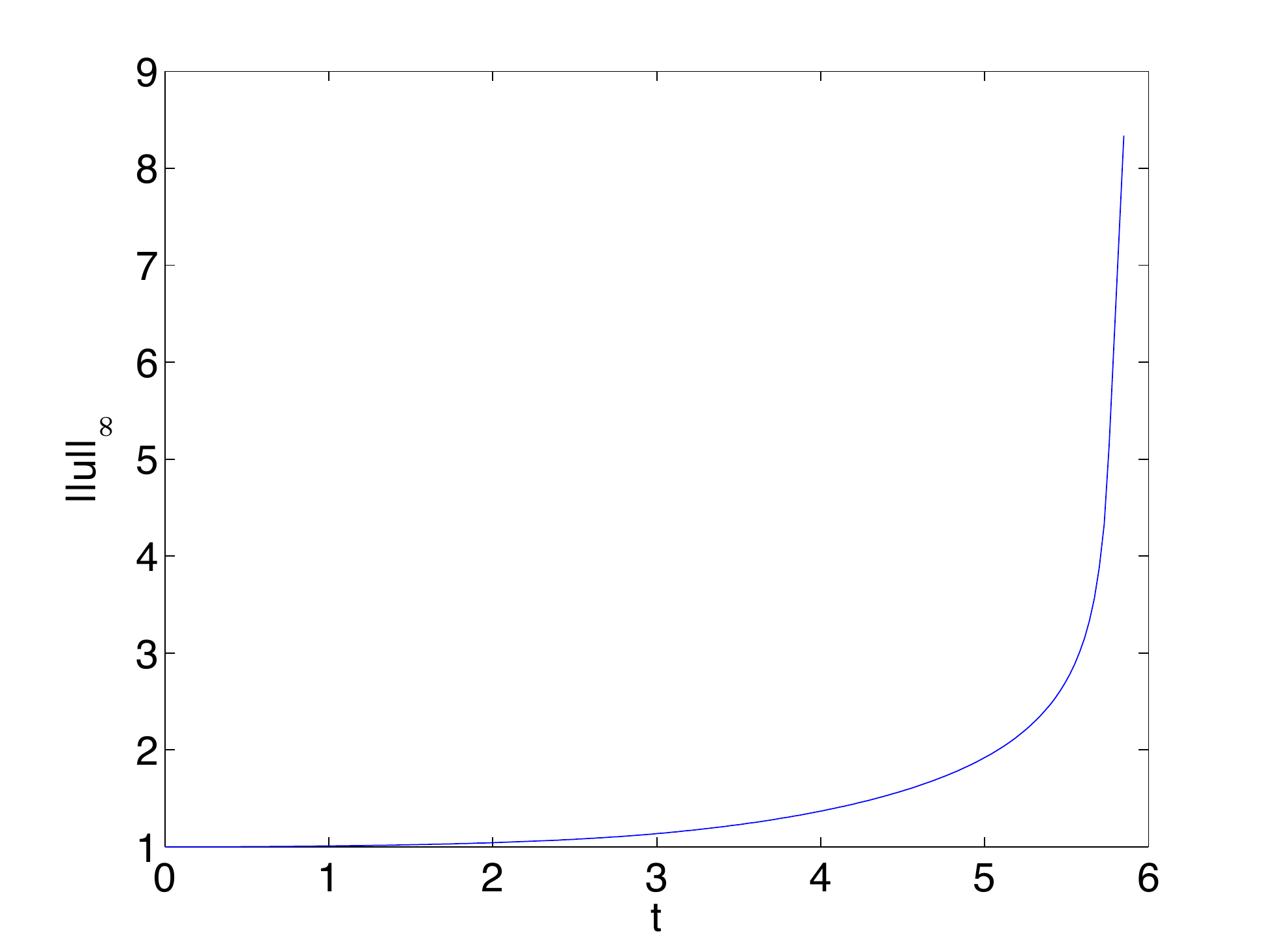}
 \caption{Solution to the fBBM equation (\ref{fracBBM}) for 
 $\alpha=0.2$ for the initial data $u_{0}=\mbox{sech}^{2}x$  
 at $t=5.85$; on the left the Fourier coefficients, on the right the 
 $L_{\infty}$ norm in dependence of time.}
 \label{gBBMsechalpha02fourier}
\end{figure}

The solution appears to have a blow-up in this case, 
but the type of the formed singularity is not obvious. The Fourier 
coefficients in Fig.~\ref{gBBMsechalpha02fourier} indicate a lack of 
resolution at $t=5.85$. The $L_{\infty}$ norm in the same figure also 
appears to diverge. If we rerun the code in this case with higher 
resolution in space and time, we get within numerical precision the 
same behavior which indicates that this is indeed a blow-up. 

A fitting of the Fourier coefficients at the last computed time 
according to the asymptotic formula (\ref{fourasymp}) yields 
$\delta=-0.003$ and $\mu+1=1.21$. The value of $\delta$ indicates 
that a singularity has formed. But there seems to be no 
$L_{\infty}$ blow-up in this case, but the formation of a cusp of the 
form $|x-x^{*}|^{\alpha}$. This is in clear contrast to corresponding 
fKdV examples studied in the previous subsection where the blow-up 
was proportional to $|x-x^{*}|^{-\alpha}$, i.e., an $L_{\infty}$ blow-up 
in contrast to the gradient catastrophe here.  Recall that the 
break-up in the Burgers solution is proportional to $(x-x_{c})^{1/3}$. 

This result is confirmed qualitatively by a study of the $L_{\infty}$ 
norm of the solution and the $L_{2}$ norm of the gradient as for fKdV 
in the previous subsection. We take the time where the singularity 
appears on the real axis as the blow-up time and study the dependence 
of the logarithms of the norms on $\ln(t^{*}-t)$. In both cases there 
is no linear dependence as in the fKdV case. Instead the $L_{\infty}$ 
norm of $u$ seems to approach a constant value in accordance with the 
expected cusp from the asymptotic behavior of the Fourier 
coefficients. 
\begin{figure}[htb!]
   \includegraphics[width=0.49\textwidth]{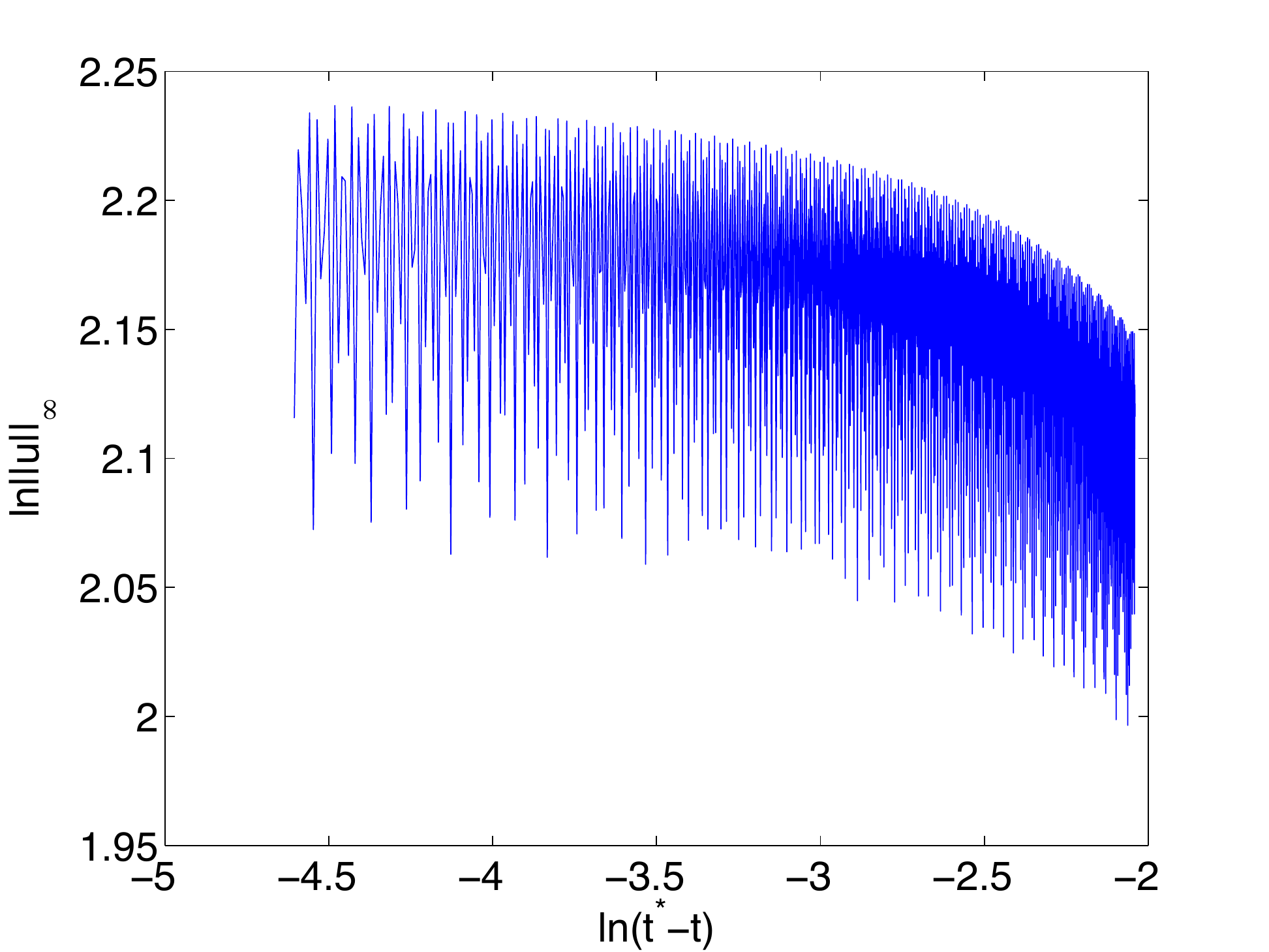}
  \includegraphics[width=0.49\textwidth]{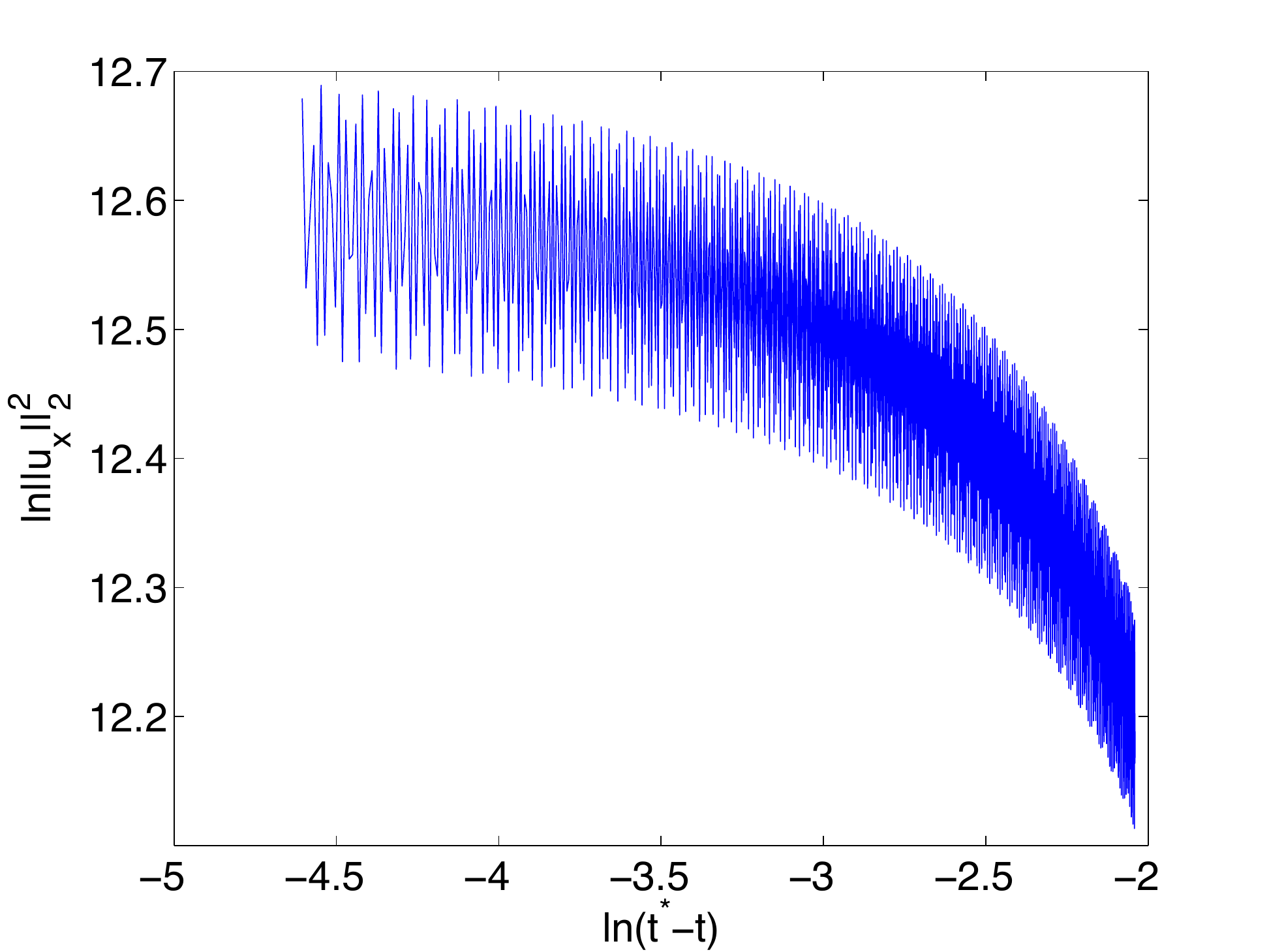}
 \caption{Norms of the solution to the fBBM equation (\ref{fracBBM}) for 
 $\alpha=0.2$ for the initial data $u_{0}=\mbox{sech}^{2}x$  
 at $t=5.85$; on the left the $L_{\infty}$ norm of $u$, on the right 
 the $L_{2}$ norm of $u_{x}$.}
 \label{gBBMsechalpha02fit}
\end{figure}

Note that the blow-up in fBBM appears for much later times than the 
blow-up for the same initial data and the same $\alpha$ for fKdV in 
Fig.~\ref{gBOsechalpha024t}. In both cases the time is larger than 
the breakup time of the solution for the Burgers equation for these 
initial data, see (\ref{tch}), which is in this case $t_{c}\sim 
1.299$.

\begin{figure}[htb!]
   \includegraphics[width=0.49\textwidth]{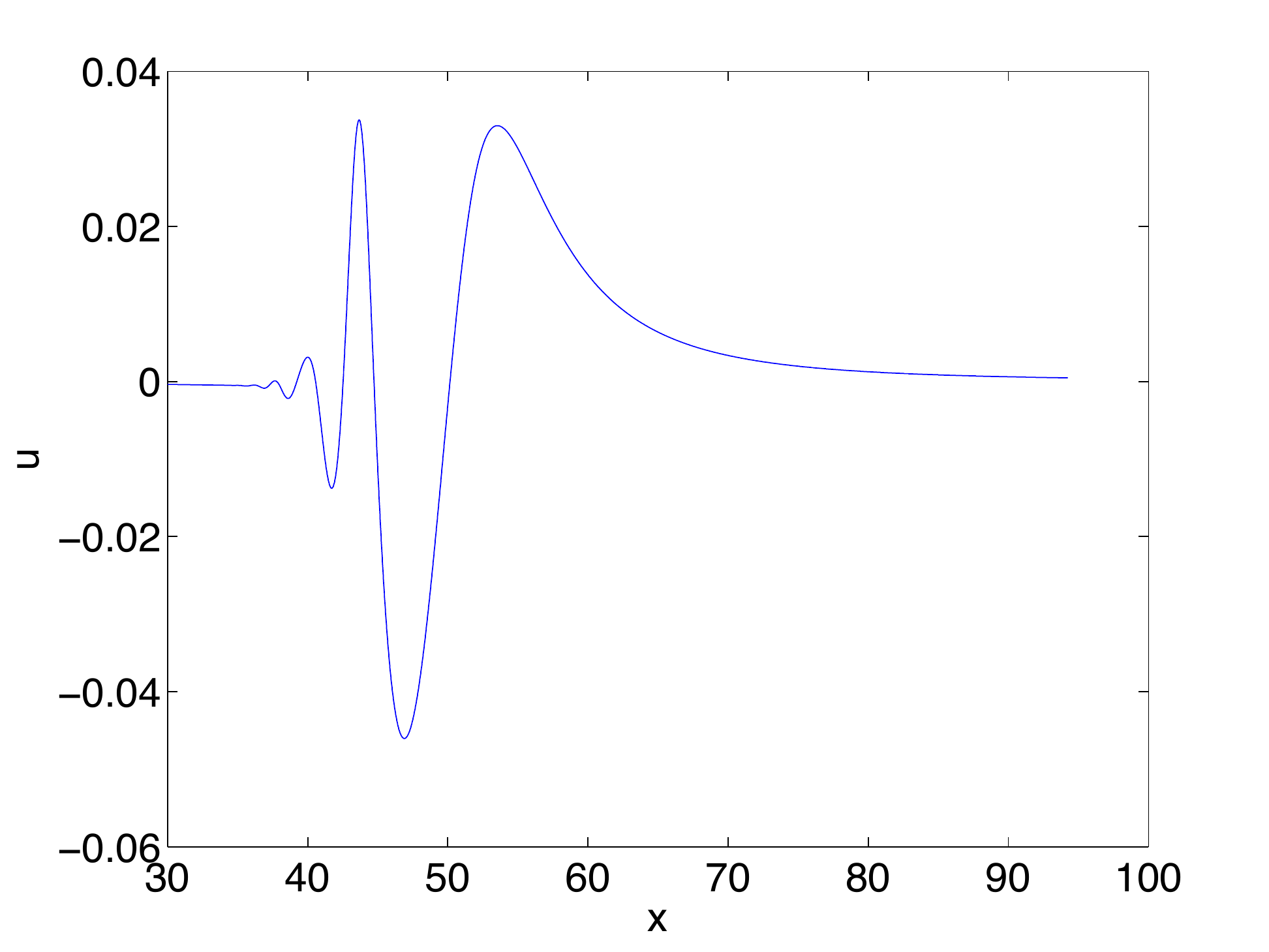}
  \includegraphics[width=0.49\textwidth]{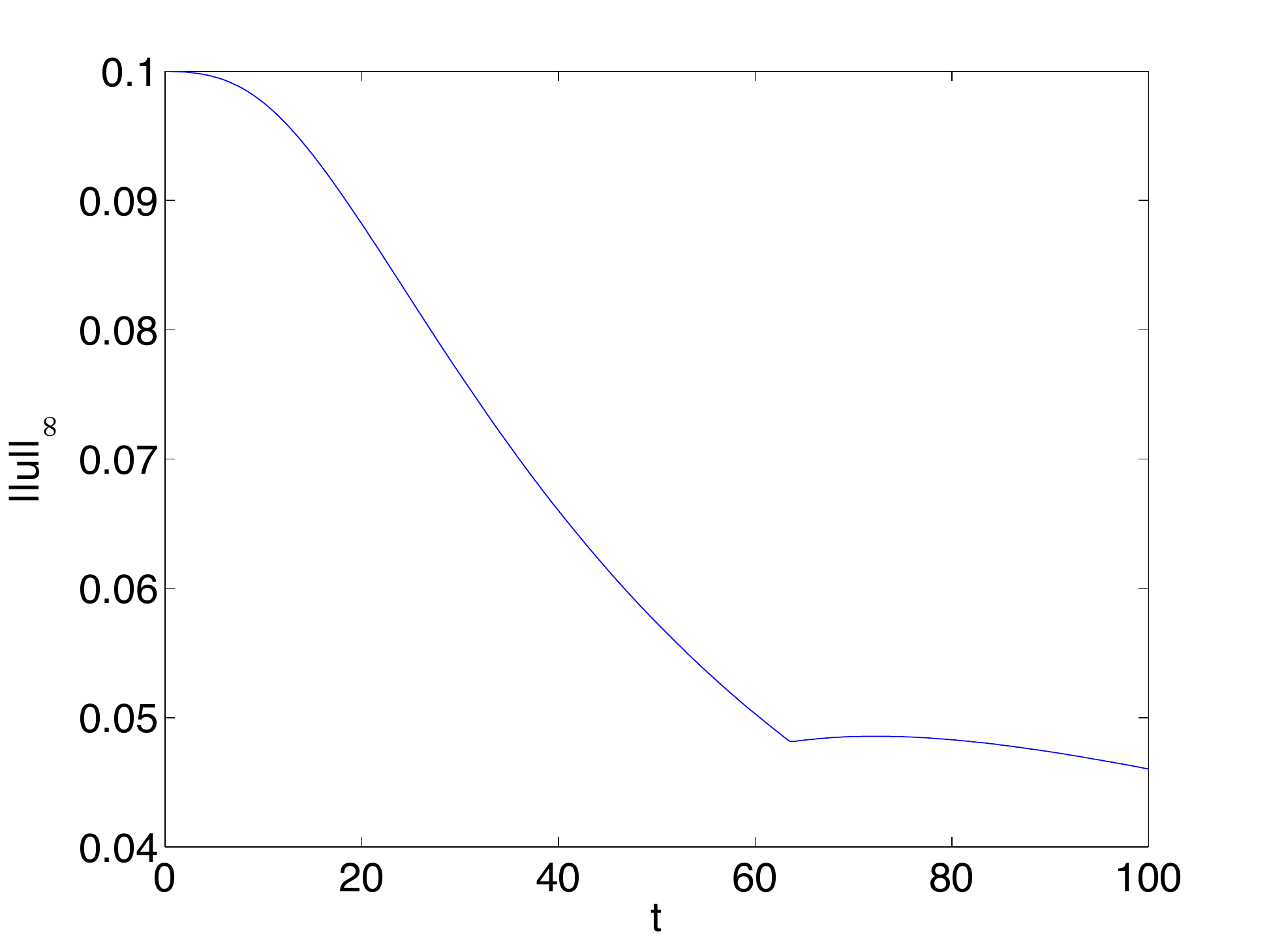}
 \caption{Solution to the fBBM equation (\ref{fracBBM}) for 
 $\alpha=0.2$ for the initial data $u_{0}=0.1\mbox{sech}^{2}x$  
 at $t=10$; on left the solution for $t=100$, on the right the 
 $L_{\infty}$ norm in dependence of time.}
 \label{gBBM01sechalpha02}
\end{figure}

For the same $\alpha$ and smaller $\beta$, the initial data appear to 
be just radiated away as can be seen in Fig.~\ref{gBBM01sechalpha02} 
for $t=100\gg t_{c}$, where $t_{c}\sim12.99$ is the critical time of Burgers
solution for the same initial data. 
The $L_{\infty}$ norm in the same figure decreases monotonically 
until the oscillations in the solution set in. 

We summarize the numerical findings in this subsection in the following 
\begin{conj}
    Consider  smooth initial data $u_{0}\in L_{2}(\mathbb{R})$ 
with a single hump. Then for
    \begin{itemize}
\item $\alpha>1/3$: solutions to the fBBM equations with the initial 
data $u_{0}$ stay smooth for 
all $t$. For large $t$ they decompose asymptotically into solitons and 
radiation.
\item $0<\alpha\leq 1/3$:  solutions to the fKdV equations  with the initial 
data $u_{0}$ and sufficiently large $L_{2}$ norm form a cusp of the form 
$|x-x^{*}|^{\alpha}$ at finite time $t^{*}$ and finite $x=x^{*}$. 
Solutions with sufficiently small initial data are global.

\item The fBBM solitons (\ref{solitarywave}) are stable for 
$\alpha>1/3$.
 \end{itemize}
\end{conj}

\begin{remark}
We note here a strong contrast  
between the gKdV equation \eqref{gKdv} and the 
generalized BBM equation
\begin{equation}\label{gbbm}
u_t+u_x+u^pu_x-u_{xxt}=0.
\end{equation}
For both \eqref{gKdv} and \eqref{gbbm}, the critical exponent for the stability of solitary waves is $p=4,$
though the explanation for instability when $p\geq 4$ is different since no blow-up occurs for \eqref{gbbm},
whatever $p.$

For the fKdV and fBBM equations the critical exponents seem to be respectively $\alpha =1/2$
and $\alpha =1/3.$
\end{remark}

\subsection{Numerical study of the Whitham equation}
In this subsection, we will study numerically solutions to the 
Whitham equation (\ref{Whit}) and for fKdV with negative $\alpha$, in particular 
$\alpha=-1/2$ which should show the same dispersion as the Whitham 
equation for the high wavenumbers. We consider again initial data of 
the form $u_{0}=\beta\mbox{sech}^{2}x$. However this time we will 
also consider negative $\beta$ since these are the initial 
data where the Whitham equation could develop solitons (see 
\cite{EGW}). For negative initial data, the 
solution will propagate to the left, and a gradient 
catastrophe of the Burgers solution is also expected to the left of the 
hump for such data. We find that the dispersion is 
strong enough for initial data of small norm to radiate the initial 
hump away to infinity. However, solutions for larger values of 
$|\beta|$ show a 
hyperbolic blow-up for $t^{*}>t_{c}$, where $t_{c}$ is the time of 
shock formation of the Burgers solution for the same initial data. 
For positive initial data, a cusp formation appears possible even 
for $t<t_{c}$. 

\begin{figure}[htb!]
   \includegraphics[width=0.49\textwidth]{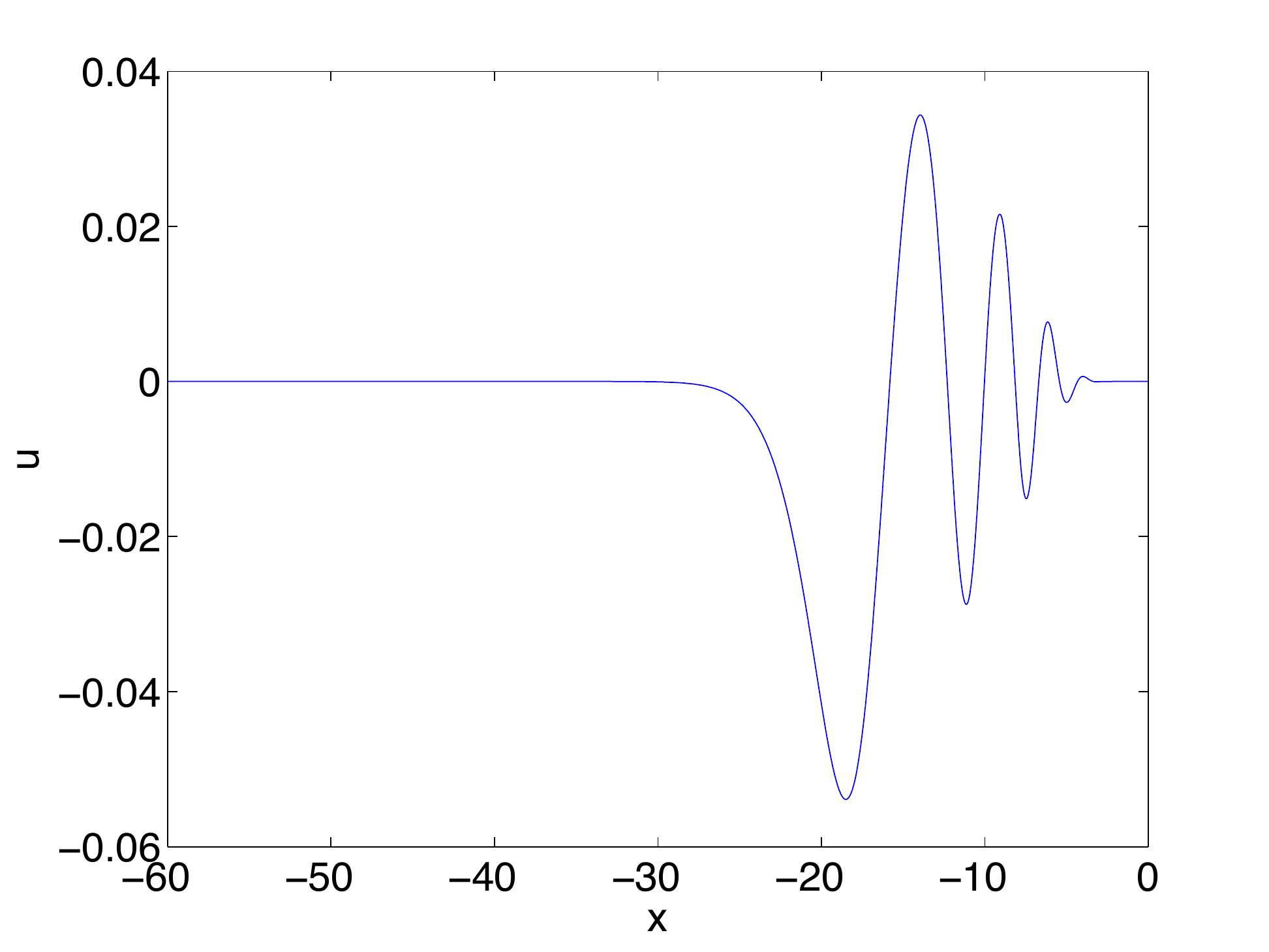}
  \includegraphics[width=0.49\textwidth]{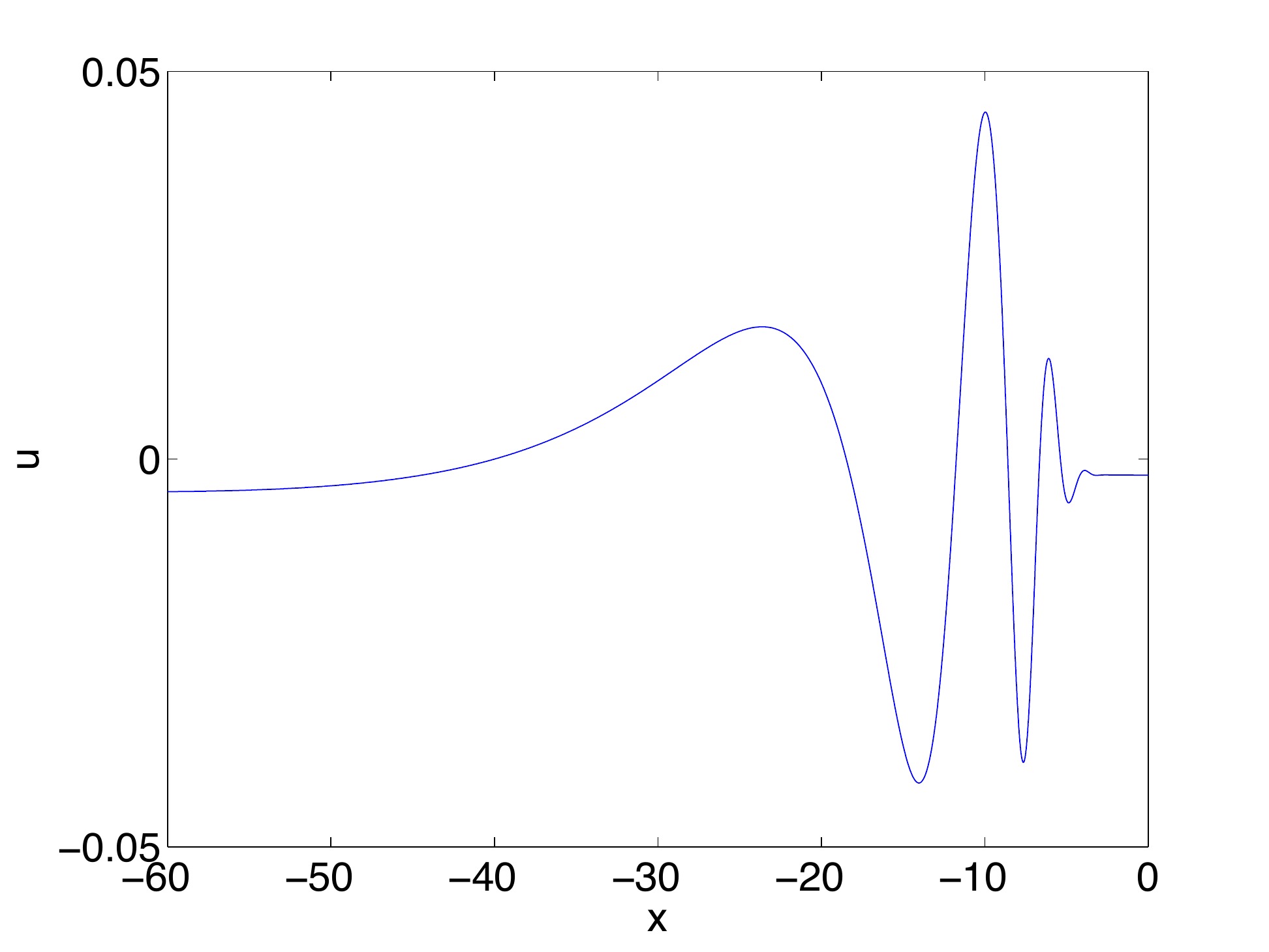}
 \caption{Solution to the Whitham equation (\ref{Whit})  for the 
 initial data $u_{0}=-0.1\mbox{sech}^{2}x$  
 at $t=20$ on left, and the solution to the fKdV equation for 
 $\alpha=-1/2$ for the same initial data at the same time on the right.}
 \label{sautwhitham01}
\end{figure}

We first study initial data $u_{0}=-0.1\mbox{sech}^{2}x$. The break-up 
time of the Burgers solution for these data is $t_{c}\sim12.99$. We use 
$N=2^{14}$ Fourier modes for $x\in20[-\pi,\pi]$ with $N_{t}=10^{4}$ 
time steps for $t<20$, i.e., larger than $t_{c}$. It can be 
seen in Fig.~\ref{sautwhitham01} that both the solution to the 
Whitham  and the fKdV equation for $\alpha=-1/2$ appear to be 
simply radiated away to infinity. The more sophisticated dispersion 
leads to more oscillations for the solution to the Whitham equation.

\begin{figure}[htb!]
   \includegraphics[width=0.49\textwidth]{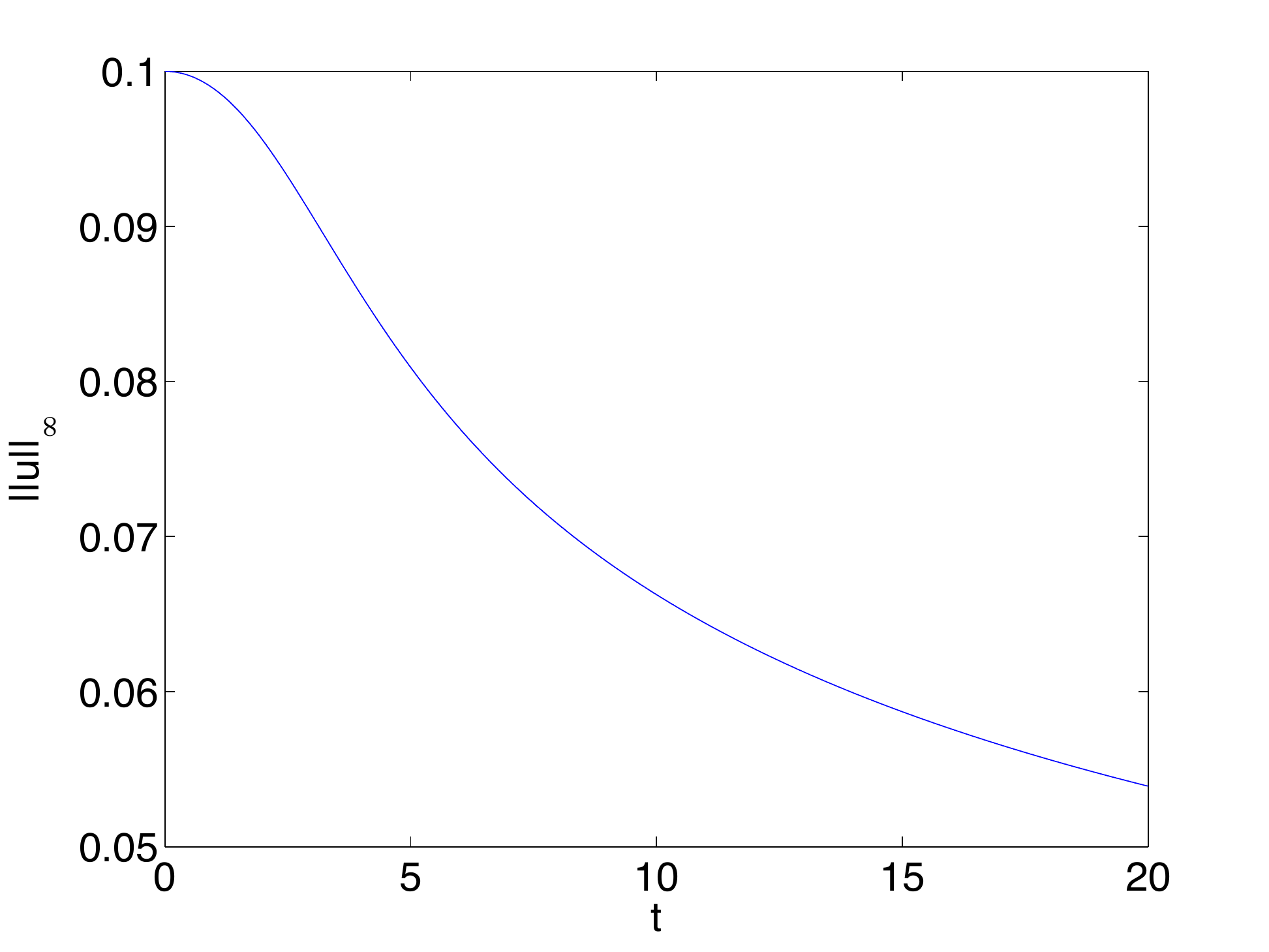}
  \includegraphics[width=0.49\textwidth]{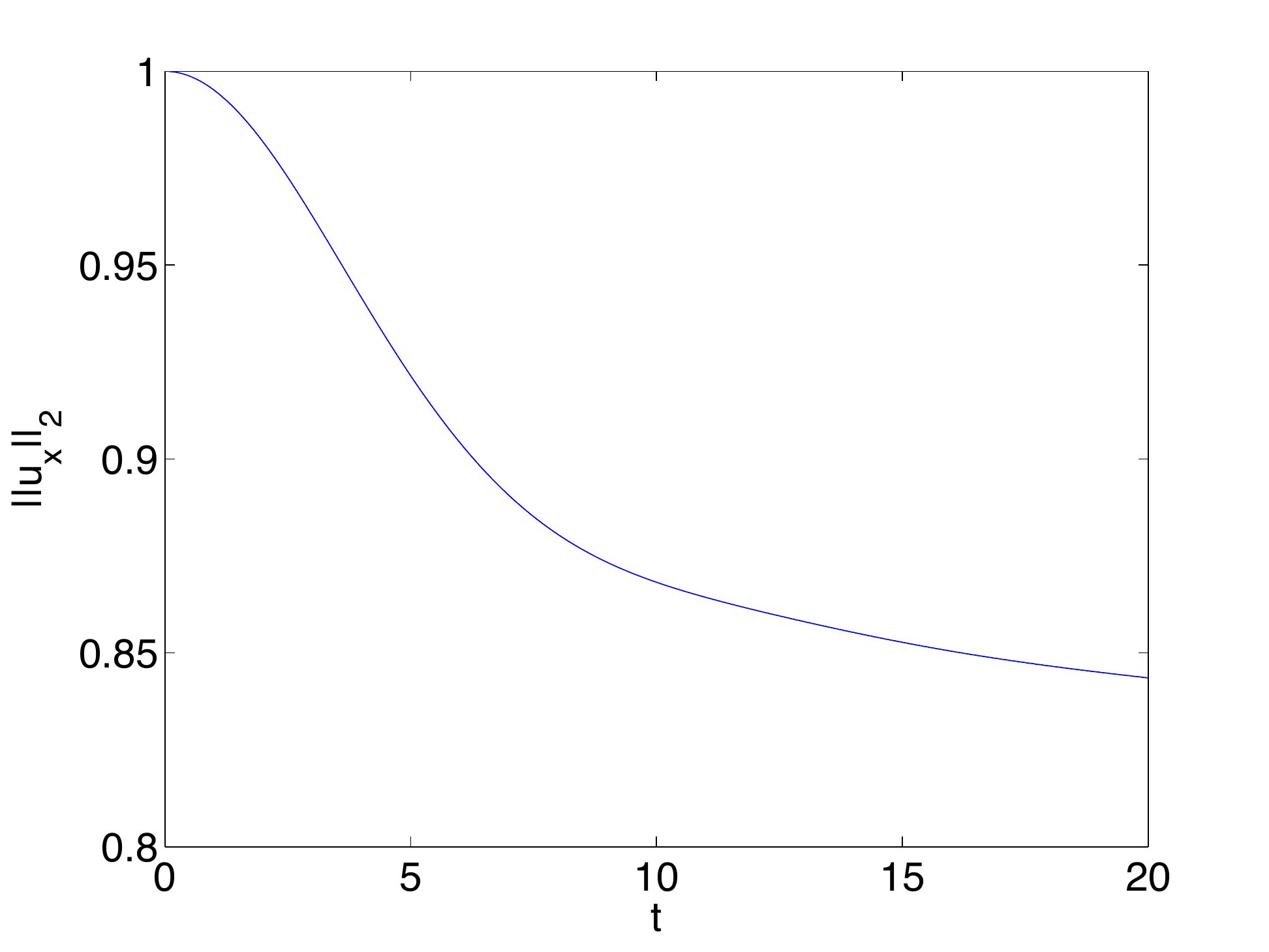}
 \caption{$L\infty$ norm of the Whitham solution of 
 Fig.~\ref{sautwhitham01} in dependence of time on the left, and the 
 corresponding
 $L_{2}$ norm of $u_{x}$ (normalized to 1 for $t=0$) on the right.}
 \label{sautwhitham01norm}
\end{figure}

There is no indication for a blow-up in this case as is even more 
obvious from various norms of the solution. In 
Fig.~\ref{sautwhitham01norm} we show the $L_{\infty}$ norm of the 
Whitham solution and the $L_{2}$ 
norm of the gradient of the solution, which are both monotonically decreasing. The corresponding norms of the 
fKdV solution in Fig.~\ref{sautwhitham01} are  not presented here, 
but are also decreasing.

\begin{figure}[htb!]
   \includegraphics[width=\textwidth]{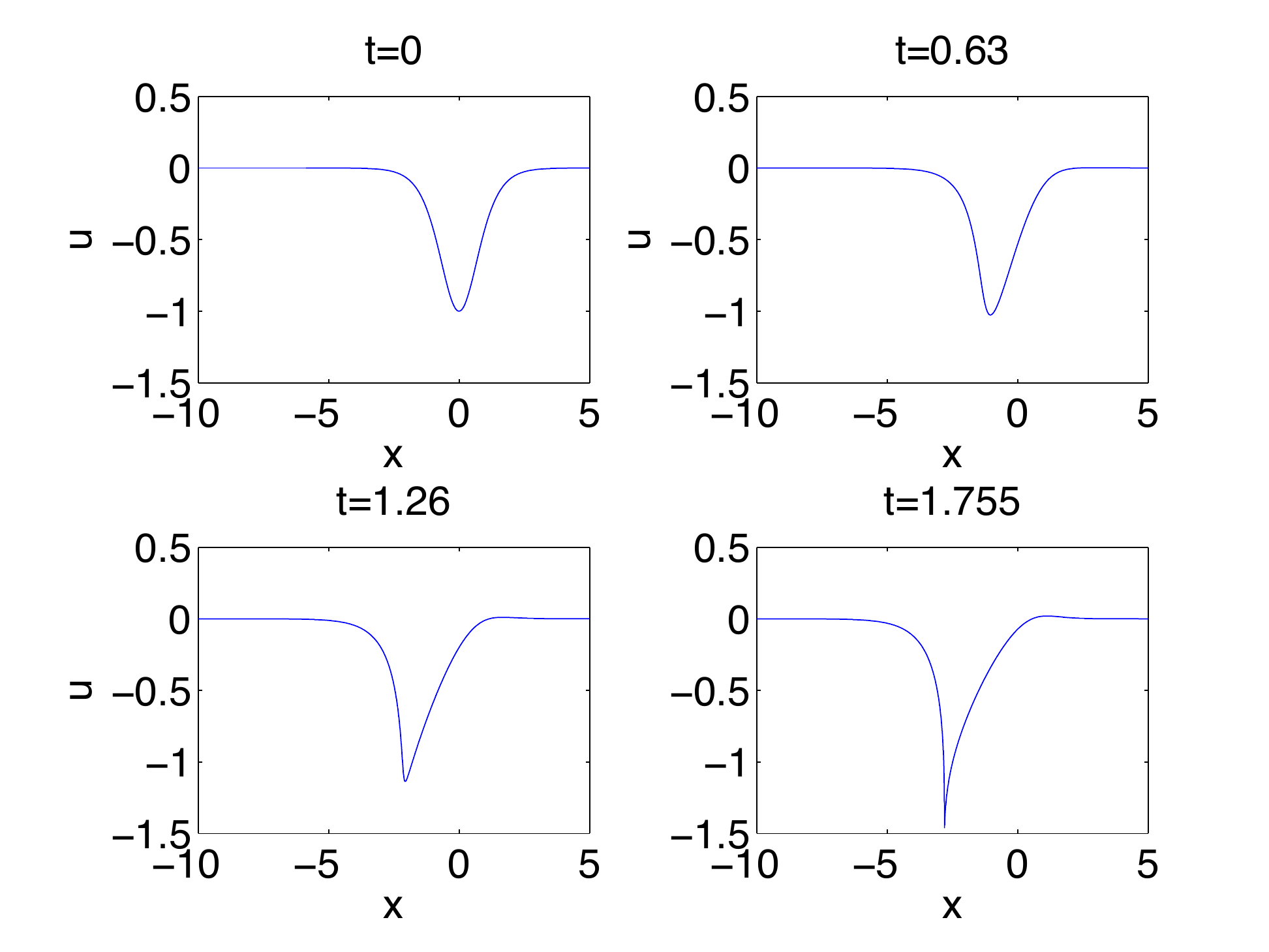}
 \caption{Solution to the Whitham equation (\ref{Whit})  for the 
 initial data $u_{0}=-\mbox{sech}^{2}x$  
 for several values of $t$.}
 \label{sautwhitham}
\end{figure}

The picture changes if we consider the same situation for the 
initial data $u_{0}=\mbox{sech}^{2}x$. We use $N=2^{16}$ Fourier 
modes for $x\in5[-\pi,\pi]$ and $N_{t}=20000$ time steps for $t<2.1$ 
(the break-up time of the Burgers solution in this case is $t_{c}\sim 
1.299$). 
Note that for the high wavenumbers, both the Whitham and the fKdV 
equation with $\alpha=-1/2$ show a weaker dispersion ($\propto 
|k|^{1/2}$) than a first order derivative which could be eliminated 
via a Galilean transformation. 
As can be seen in 
Fig.~\ref{sautwhitham}, the solution shows for early times the same 
behavior as the corresponding Burgers solution, a steepening of one 
front of the solution with an increasing gradient. At a given point 
the maximum of the solution turns into a  cusp forming for $t>t_{c}$. 
This is a different type of singularity as observed in the Burgers shock 
where the inflection point becomes singular.
The forming of the cusp leads to an increase of the modulus of the Fourier coefficients 
for the high wavenumbers.    We fit the 
Fourier coefficients to the formula (\ref{fourasymp}). The code is 
stopped when $\delta\sim 10^{-6}$, where the 
numerically computed energy at the last shown time in 
Fig.~\ref{sautwhitham} is still of the order of $10^{-10}$. We find 
$\mu+1\sim1.36$ which indicates a cusp of the form $u\sim 
|x-x^{*}|^{1/3}$, i.e., the same one would find for the Burgers 
shock. 
Note that no solitary wave seems to be 
forming in the solution, but the existence proof in
\cite{EGW} does not provide an estimate for the velocity of the solitary wave which is obtained as a Lagrange multiplier.

\begin{figure}[htb!]
   \includegraphics[width=0.49\textwidth]{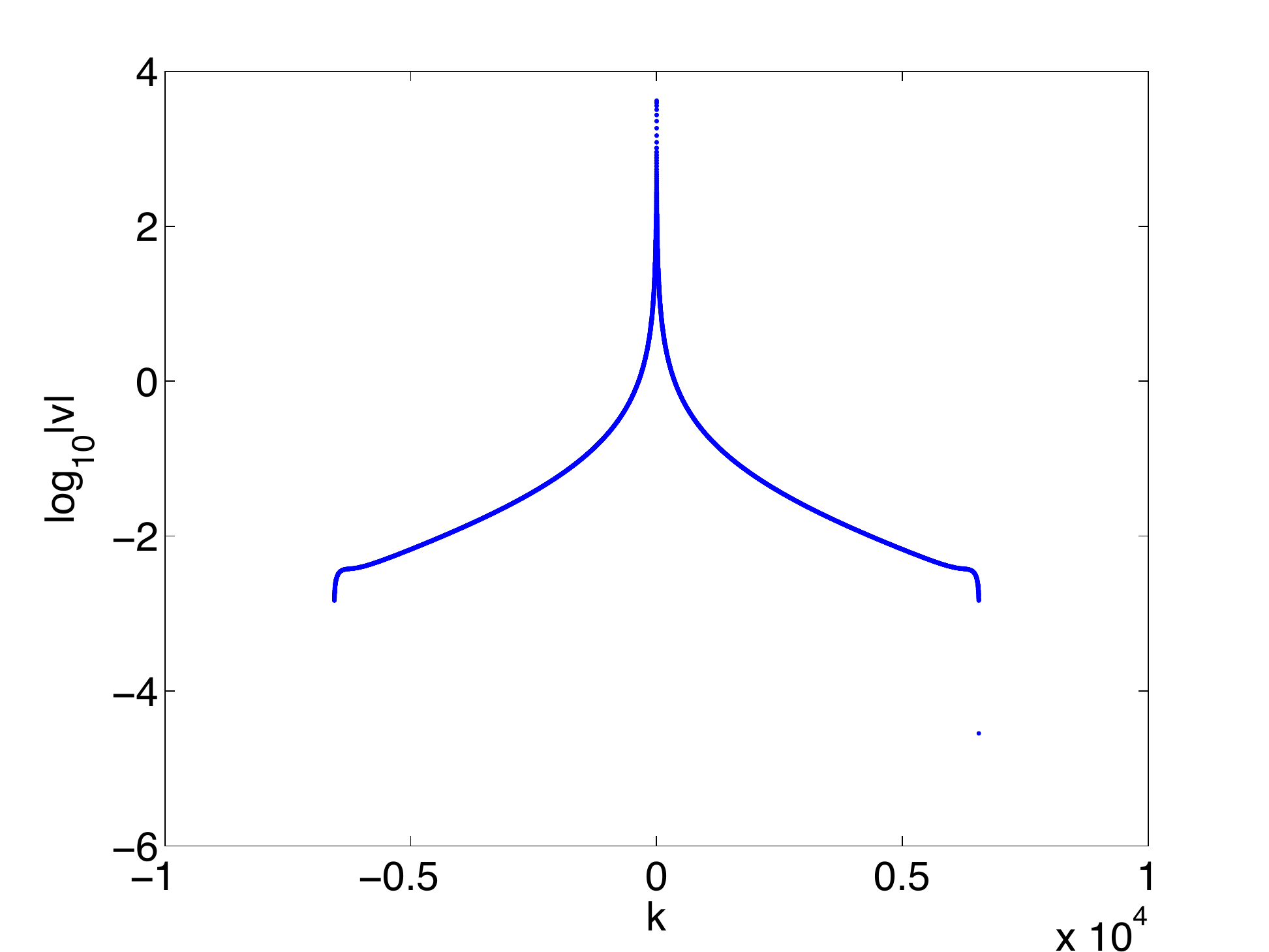}
  \includegraphics[width=0.49\textwidth]{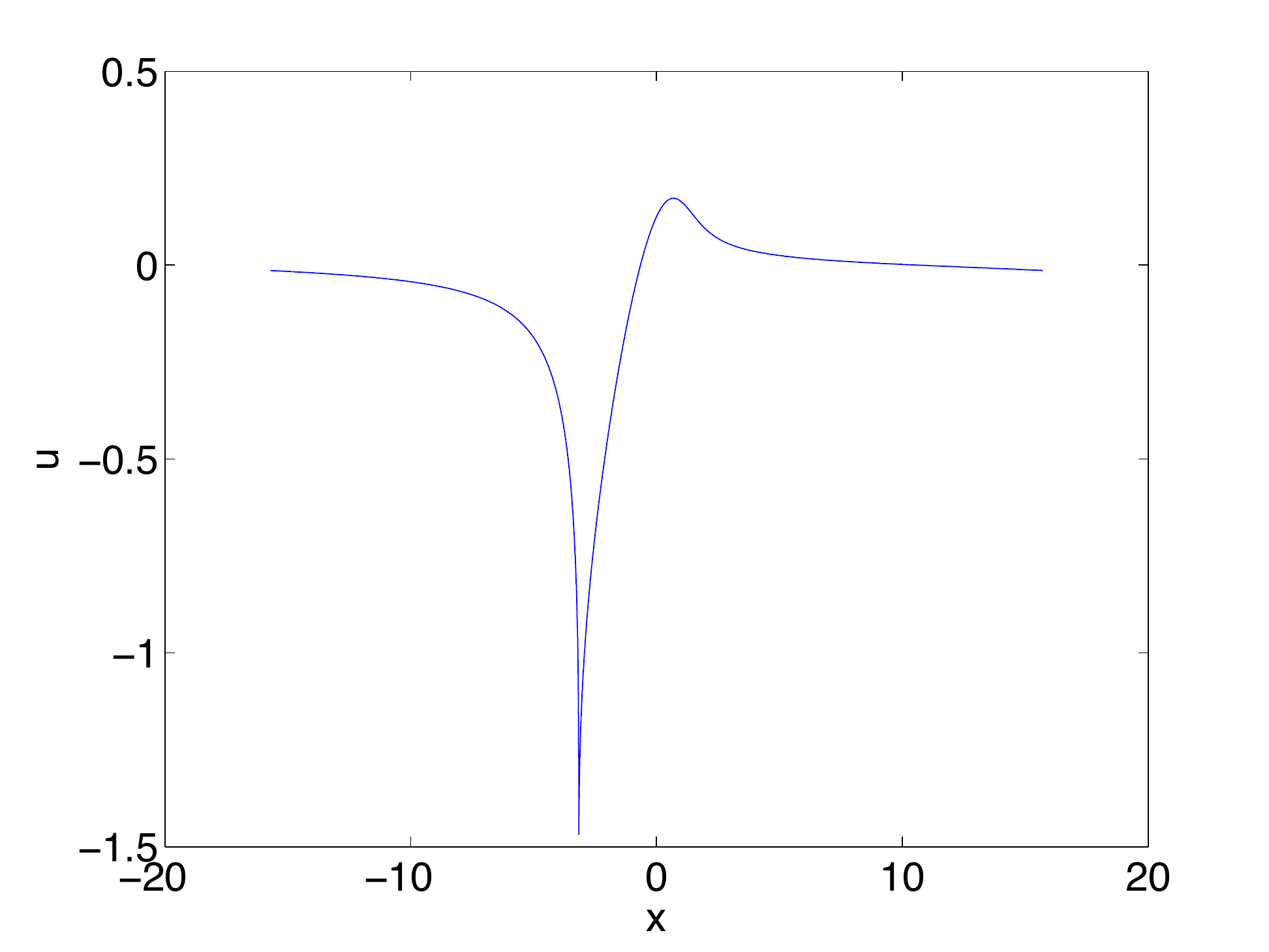}
 \caption{Modulus of the Fourier coefficients of the solution to the 
 Whitham equation for the initial data $u_{0}=-\mbox{sech}^{2}x$ for 
 $t=1.764$ on the left, and the solution for the fKdV equation for 
 the same initial data for $t=2.001$ on the right.}
 \label{sautwhithamfourier}
\end{figure}

The Fourier coefficients at the last shown time in 
Fig.~\ref{sautwhitham} can be seen in Fig.~\ref{sautwhithamfourier}. 
The solution to the fKdV equation 
for $\alpha=-1/2$ for the same initial data can be seen in the same 
figure. Here the solution becomes singular (as indicated by a 
vanishing for the parameter $\delta$ in (\ref{fourasymp})) at a 
later time, and we find $\mu+1\sim1.515$.  Thus there seem to be the same reasons for the 
singularity formation  in solutions to the Whitham and the fKdV 
equation for $\alpha=-1/2$.

\begin{figure}[htb!]
   \includegraphics[width=0.49\textwidth]{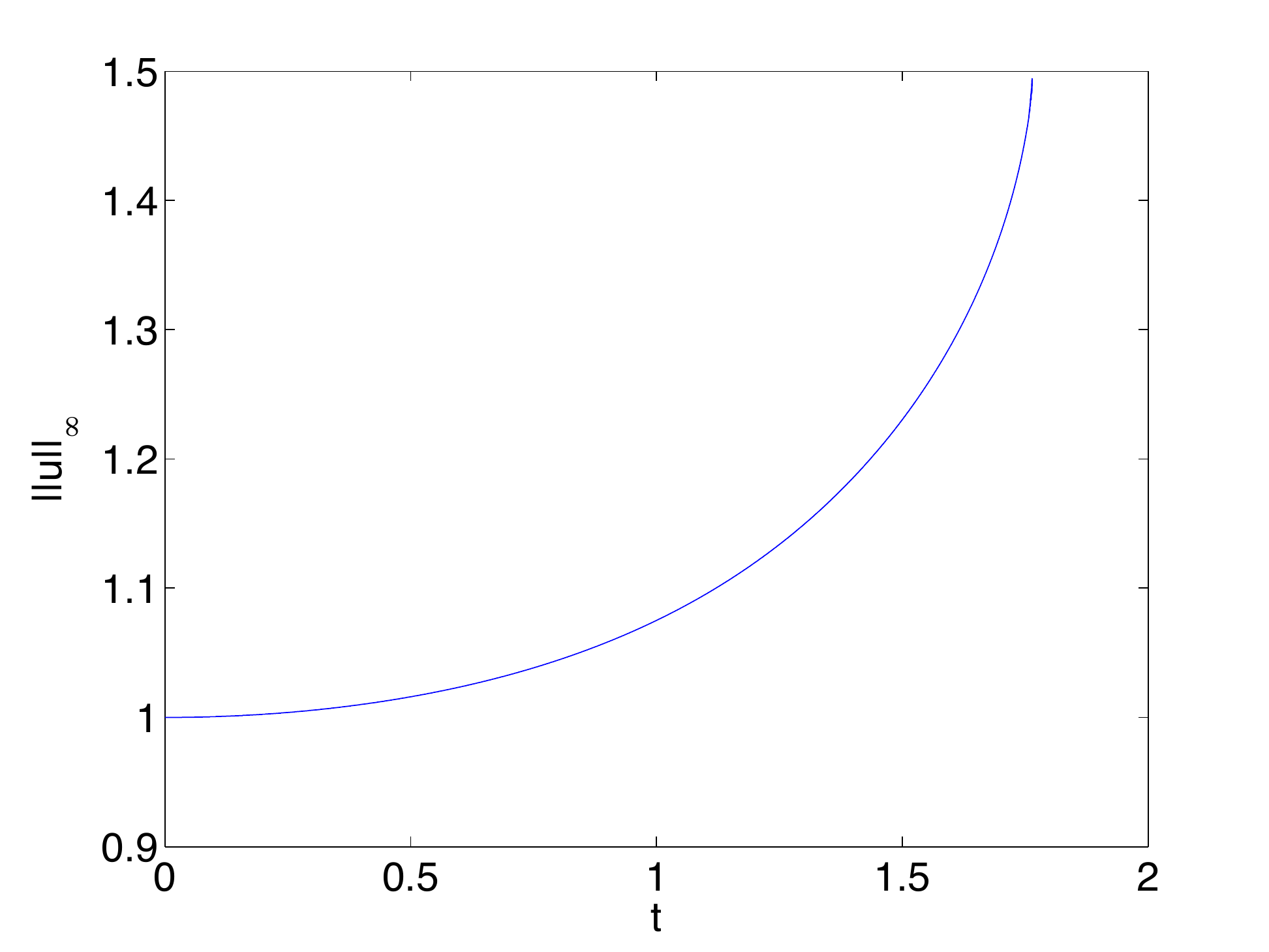}
  \includegraphics[width=0.49\textwidth]{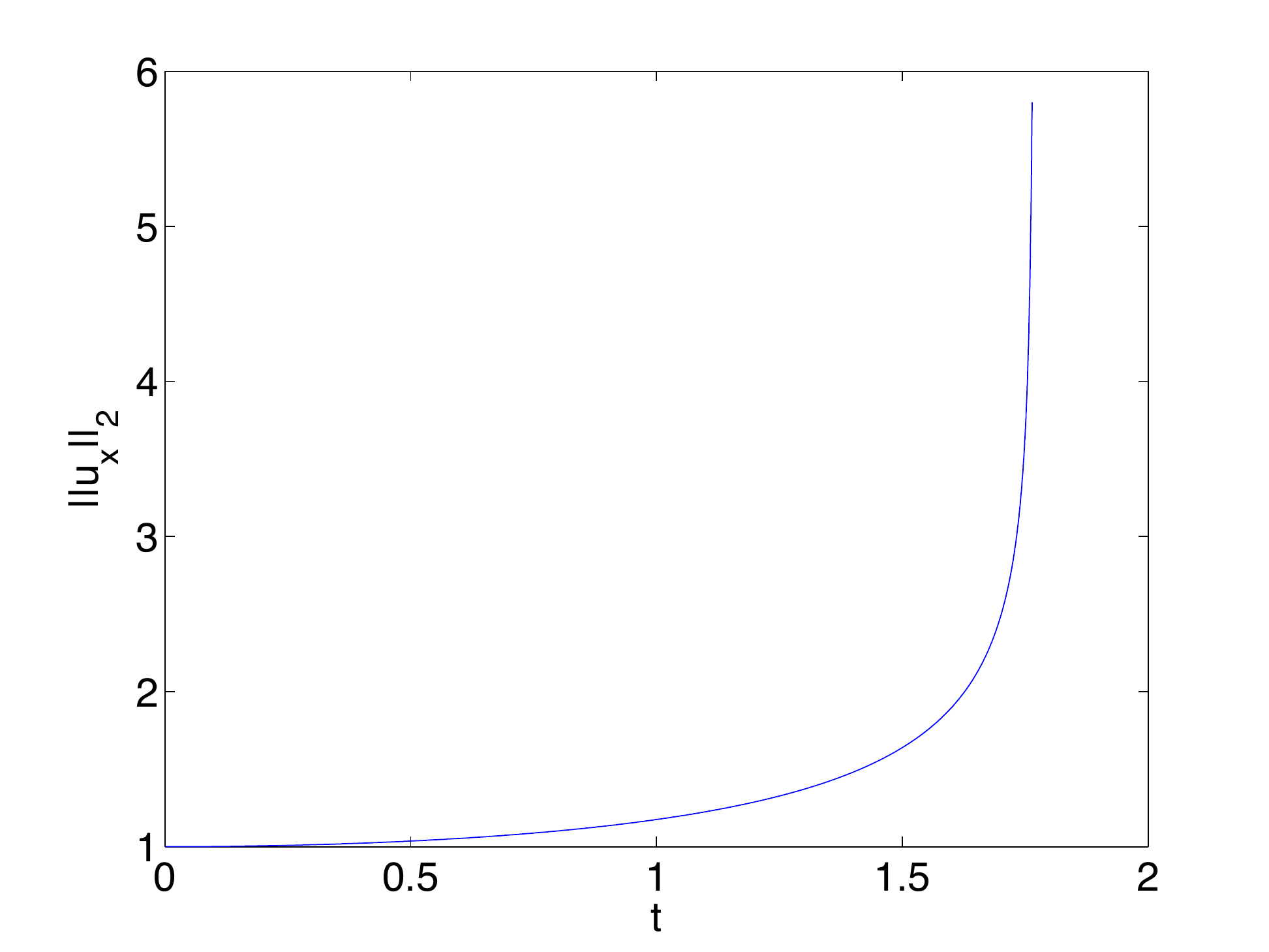}
 \caption{$L\infty$ norm of the Whitham solution of 
 Fig.~\ref{sautwhitham} in dependence of time on the left, and the 
 corresponding
 $L_{2}$ norm of $u_{x}$ on the right.}
 \label{sautwhithamnorm}
\end{figure}

The type of the singularity is also confirmed by the norms of the 
solution to the Whitham equation (the corresponding norms for the 
fKdV solution are very similar and thus not shown) in 
Fig.~\ref{sautwhithamnorm}. Both norms are only moderately increasing 
at the time $t^{*}$. 

\begin{figure}[htb!]
   \includegraphics[width=0.49\textwidth]{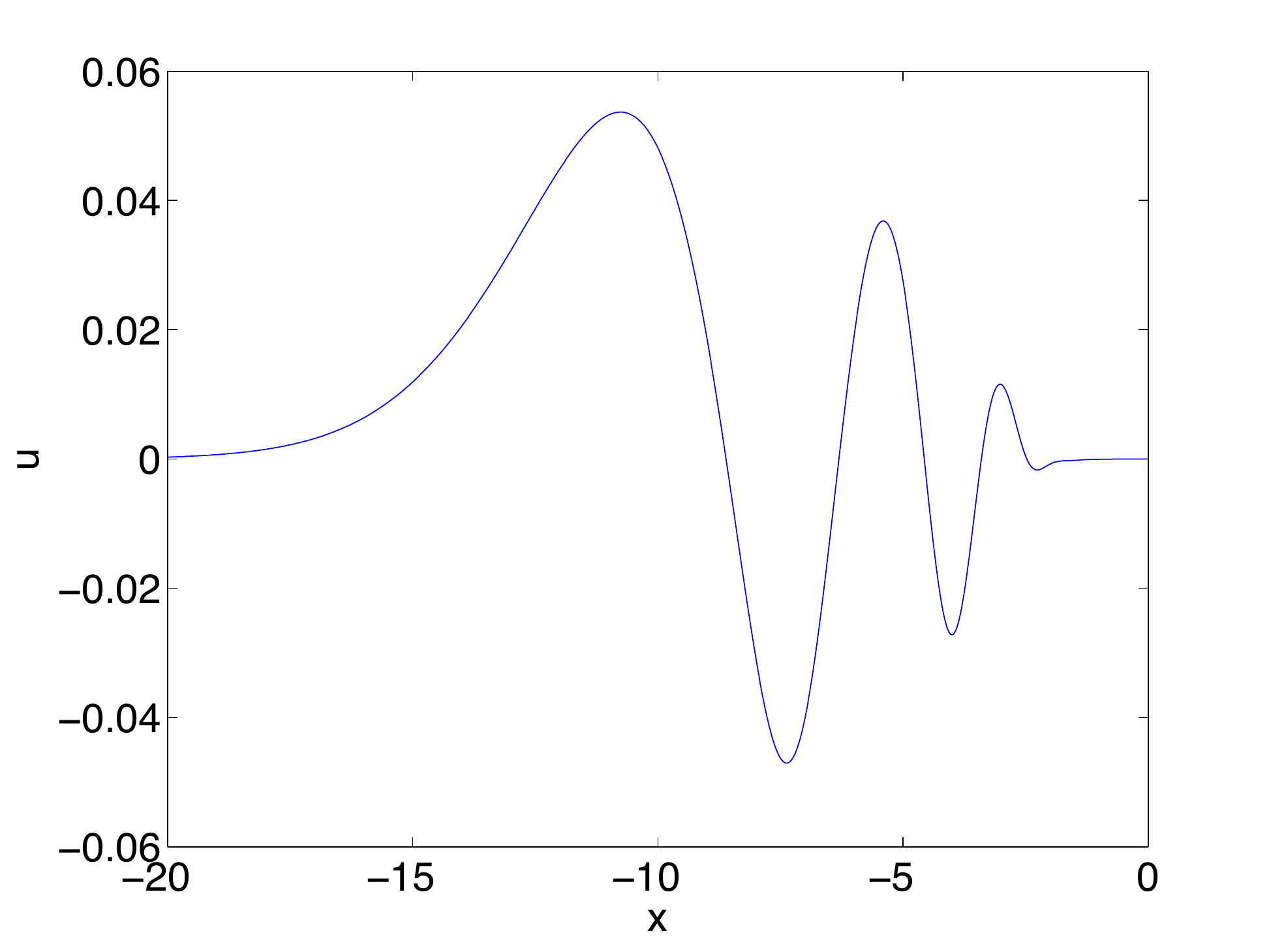}
  \includegraphics[width=0.49\textwidth]{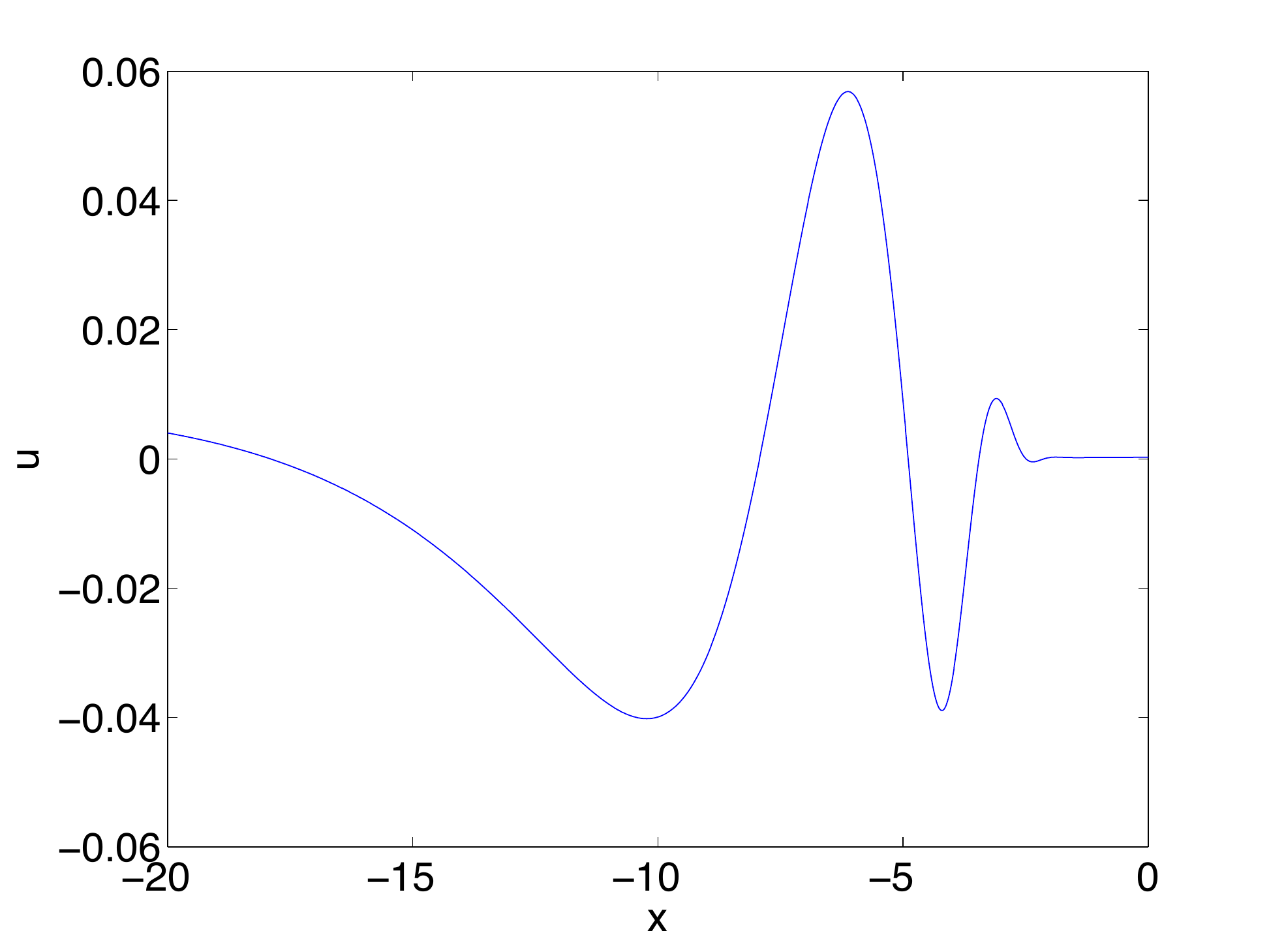}
 \caption{Solution to the Whitham equation (\ref{Whit})  for the initial data $u_{0}=0.1\mbox{sech}^{2}x$  
 at $t=13$ on left, and the solution to the fKdV equation for 
 $\alpha=-1/2$ for the same initial data at the same time on the right.}
 \label{sautwhitham01p}
\end{figure}

Although there are no solitons in this case, 
we also study positive initial data to allow for a comparison with 
the results for fKdV in the previous subsections. Initial data as 
$u_{0}=0.1\mbox{sech}^{2}x$ of small mass again appear to be just 
radiated away for both the Whitham equation and fKdV with $\alpha=-1/2$
as can be seen in Fig.~\ref{sautwhitham01p}. There is also no 
indication for a blow-up in this case from the $L_{\infty}$ norm of $u$ 
or the $L_{2}$ norm of $u_{x}$. 

\begin{figure}[htb!]
   \includegraphics[width=\textwidth]{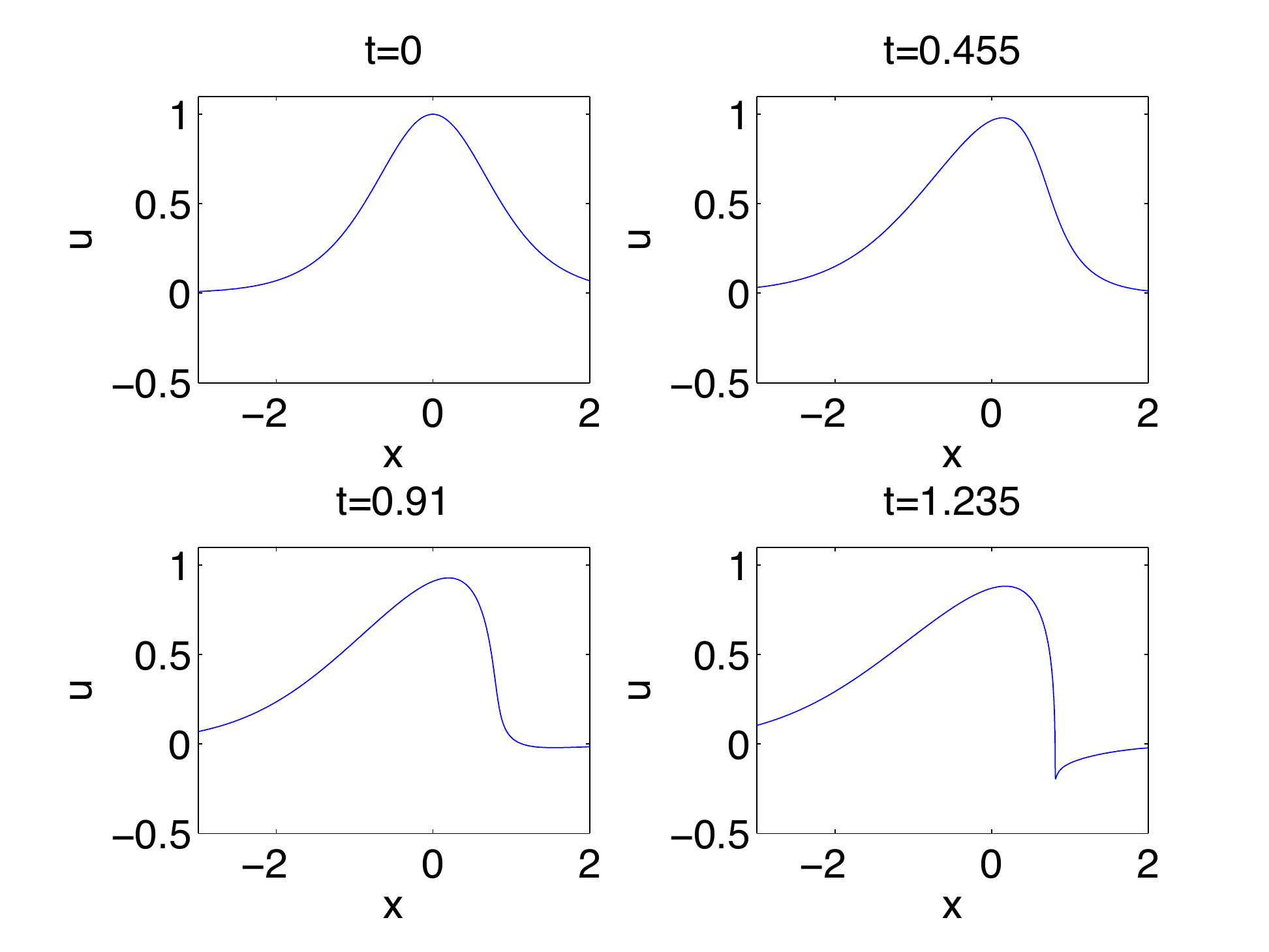}
 \caption{Solution to the Whitham equation (\ref{Whit})  for the initial data $u_{0}=\mbox{sech}^{2}x$  
 for several values of $t$.}
 \label{sautwhithamp}
\end{figure}

Initial data with more mass as 
$u_{0}=\mbox{sech}^{2}x$ lead to a behavior different from both 
the small mass and the blow-up case for negative initial data. We use $N=2^{16}$ Fourier 
modes for $x\in5[-\pi,\pi]$ and $N_{t}=20000$ time steps for $t<1.3$ 
(the break-up time of the Burgers solution in this case is $t_{c}\sim 
1.299$). This is the only time in this paper that the code breaks due 
to \emph{aliasing errors}, i.e., due to a growing of the Fourier modes for 
the high wave numbers. This can be understood in a hand-waving manner 
as follows: for the high wavenumbers, both the Whitham and the fKdV 
equation with $\alpha=-1/2$ show a weaker dispersion ($\propto 
|k|^{1/2}$) than a first order derivativen. As can be seen in 
Fig.~\ref{sautwhitham}, the solution shows for early times the same 
behavior as the corresponding Burgers solution, a steepening of one 
front of the solution with an increasing gradient. At a given point 
there is a cusp forming for $t<t_{c}$ to the right of the point where 
the point of gradient catastrophe of the Burgers solution would appear. 
This leads to an increase of the modulus of the Fourier coefficients 
for the high wavenumbers. The weak dispersion together with 
unavoidable numerical errors leads to an 
amplification of this phenomenon which would eventually break the 
code. 
To address these aliasing problems we use dealiasing according to the 
$2/3$ rule, i.e., we put the Fourier coefficients corresponding to 
the $1/3$ highest wavenumbers equal to zero.   We fit the remaining 
Fourier coefficients to the formula (\ref{fourasymp}). The code is 
stopped when $\delta\sim 10^{-6}$, where the 
numerically computed energy at the last shown time in 
Fig.~\ref{sautwhithamp} is still of the order of $10^{-11}$. We find 
$\mu+1\sim1.516$ which indicates a cusp of the form $u\sim 
|x-x^{*}|^{1/2}$.

\begin{figure}[htb!]
   \includegraphics[width=0.49\textwidth]{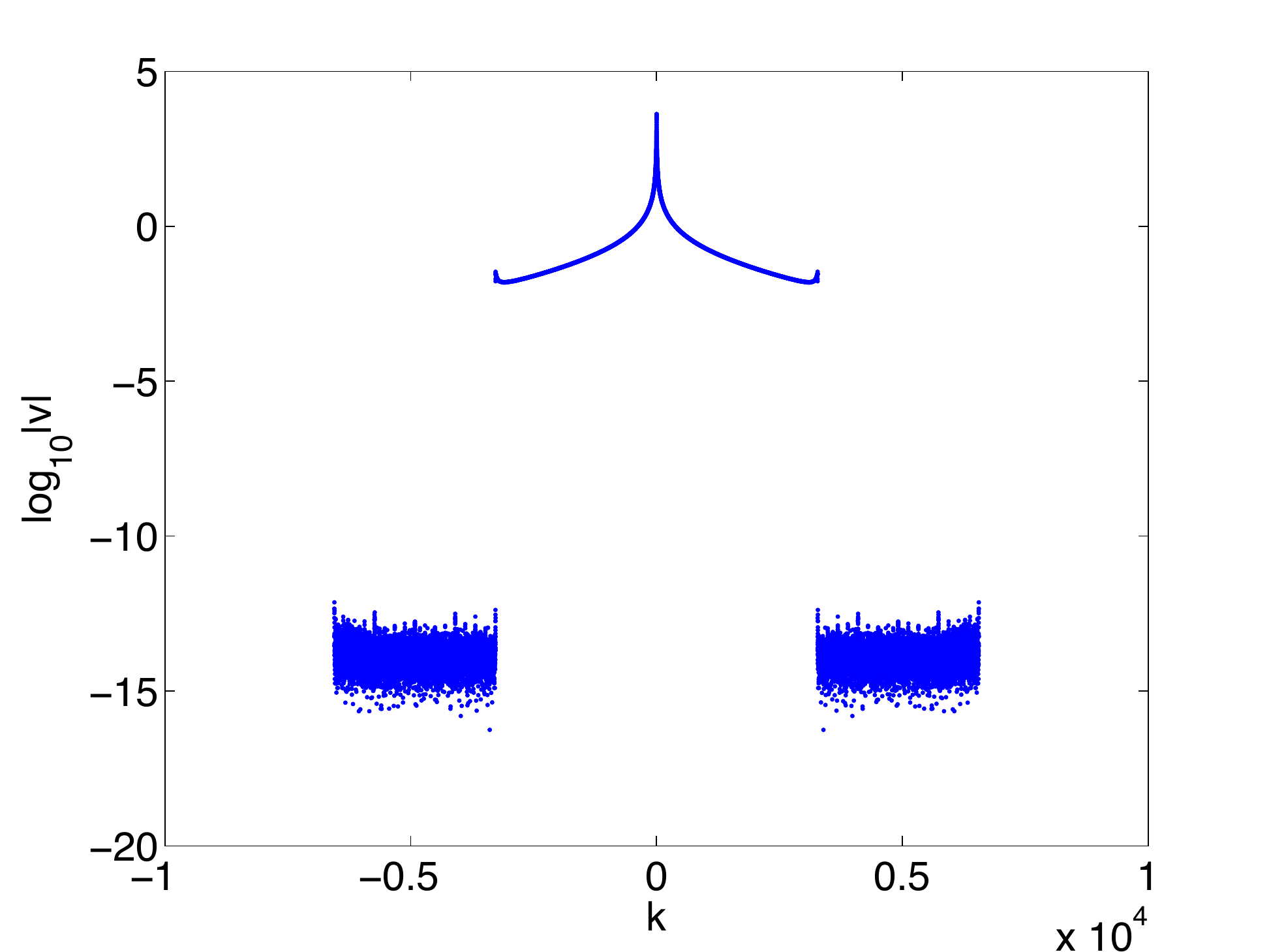}
  \includegraphics[width=0.49\textwidth]{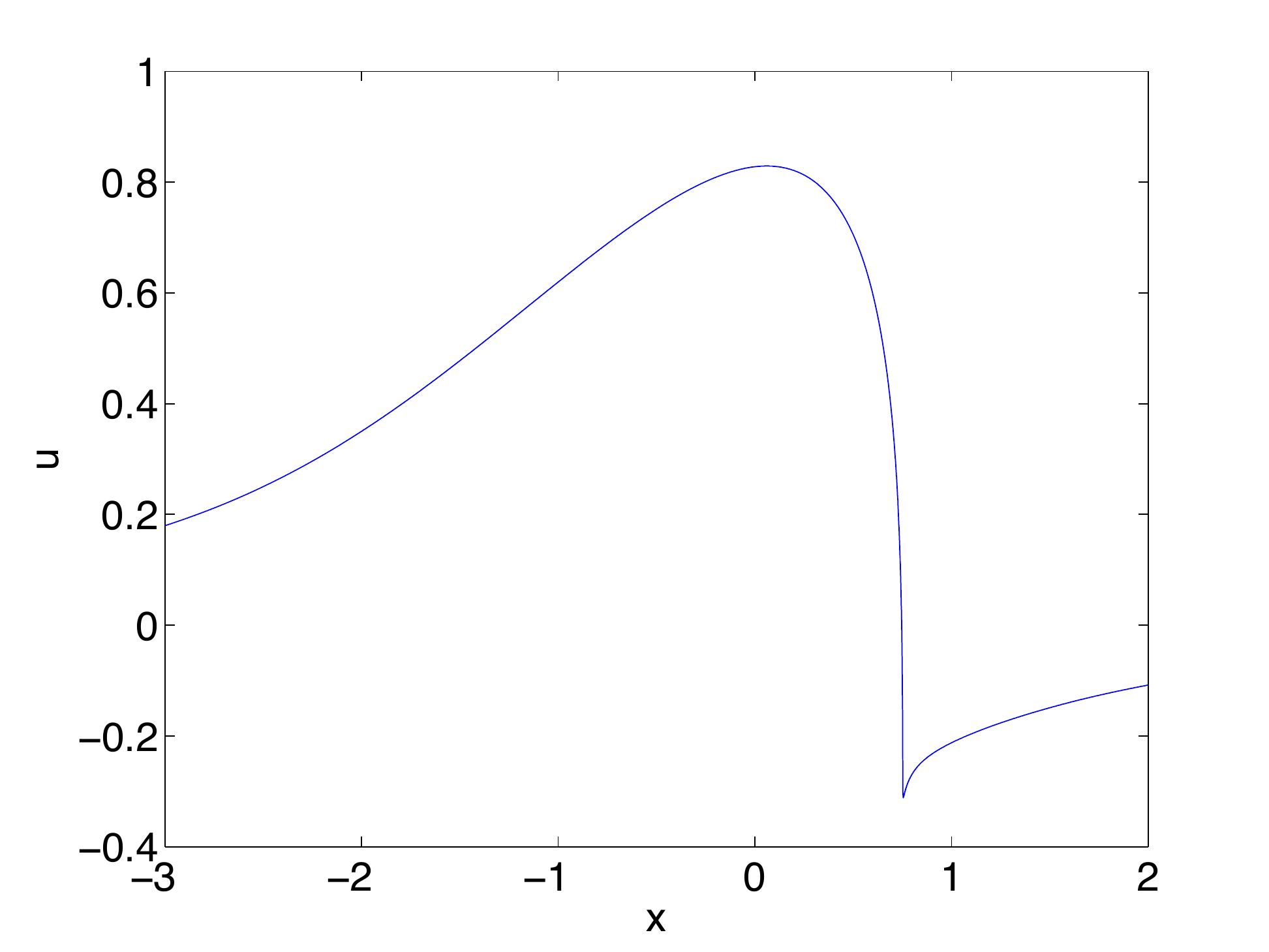}
 \caption{Modulus of the Fourier coefficients of the solution to the 
 Whitham equation for the initial data $u_{0}=\mbox{sech}^{2}x$ for 
 $t=1.235$ on the left, and the solution for the fKdV equation for 
 the same initial data for $t=1.2109$ on the right.}
 \label{sautwhithamfourierp}
\end{figure}

The Fourier coefficients at the last shown time in 
Fig.~\ref{sautwhitham} can be seen in Fig.~\ref{sautwhithamfourierp}. 
The dealiasing is clearly visible. The solution to the fKdV equation 
for $\alpha=-1/2$ for the same initial data can be seen in the same 
figure. Here the solution becomes singular (as indicated by a 
vanishing for the parameter $\delta$ in (\ref{fourasymp})) at an even
earlier time, and we find $\mu+1\sim1.515$.  Thus there seem to be the same reasons for the 
singularity formation  in solutions to the Whitham and the fKdV 
equation for $\alpha=-1/2$ also for positive initial data.

\begin{figure}[htb!]
   \includegraphics[width=0.49\textwidth]{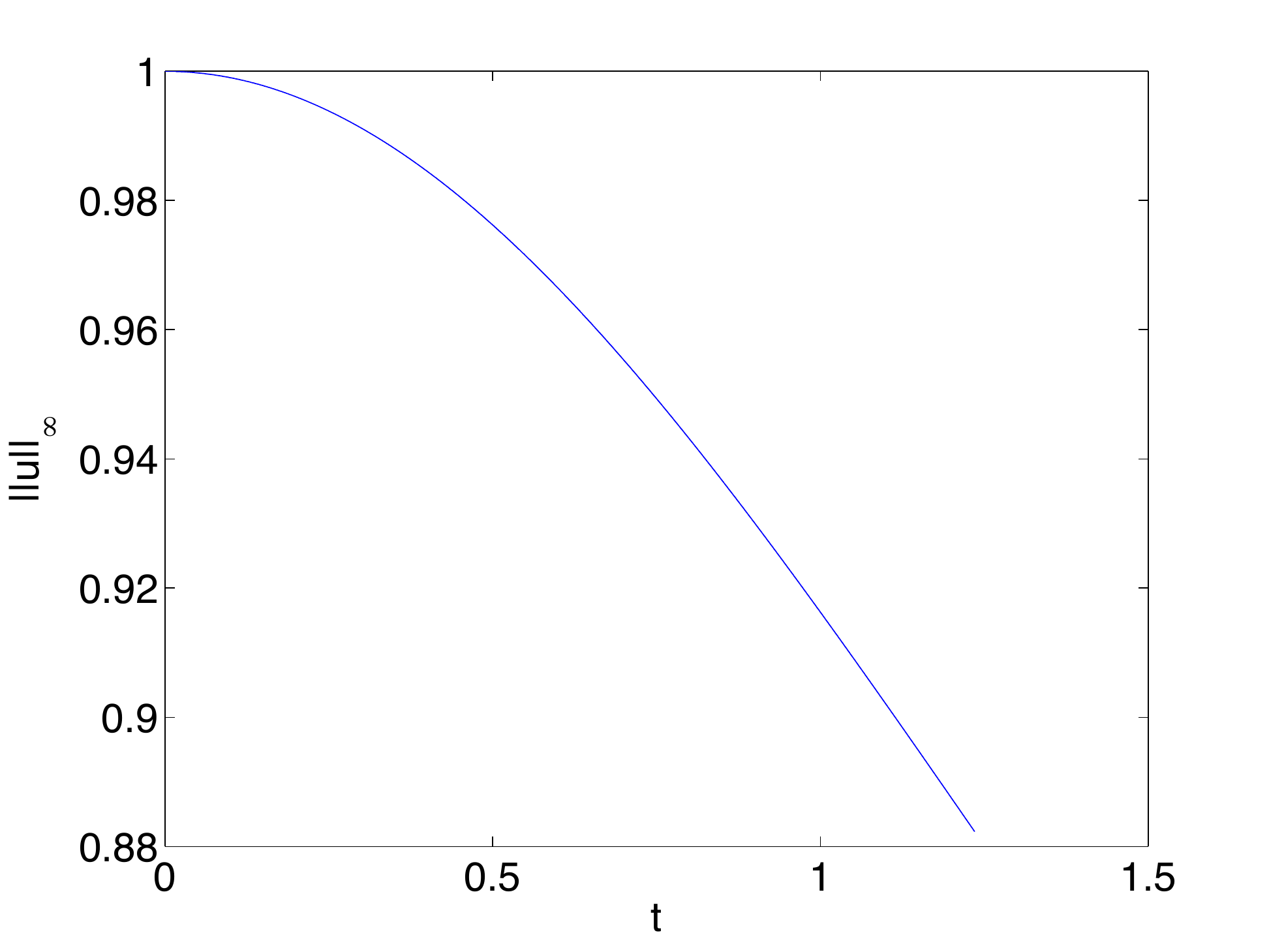}
  \includegraphics[width=0.49\textwidth]{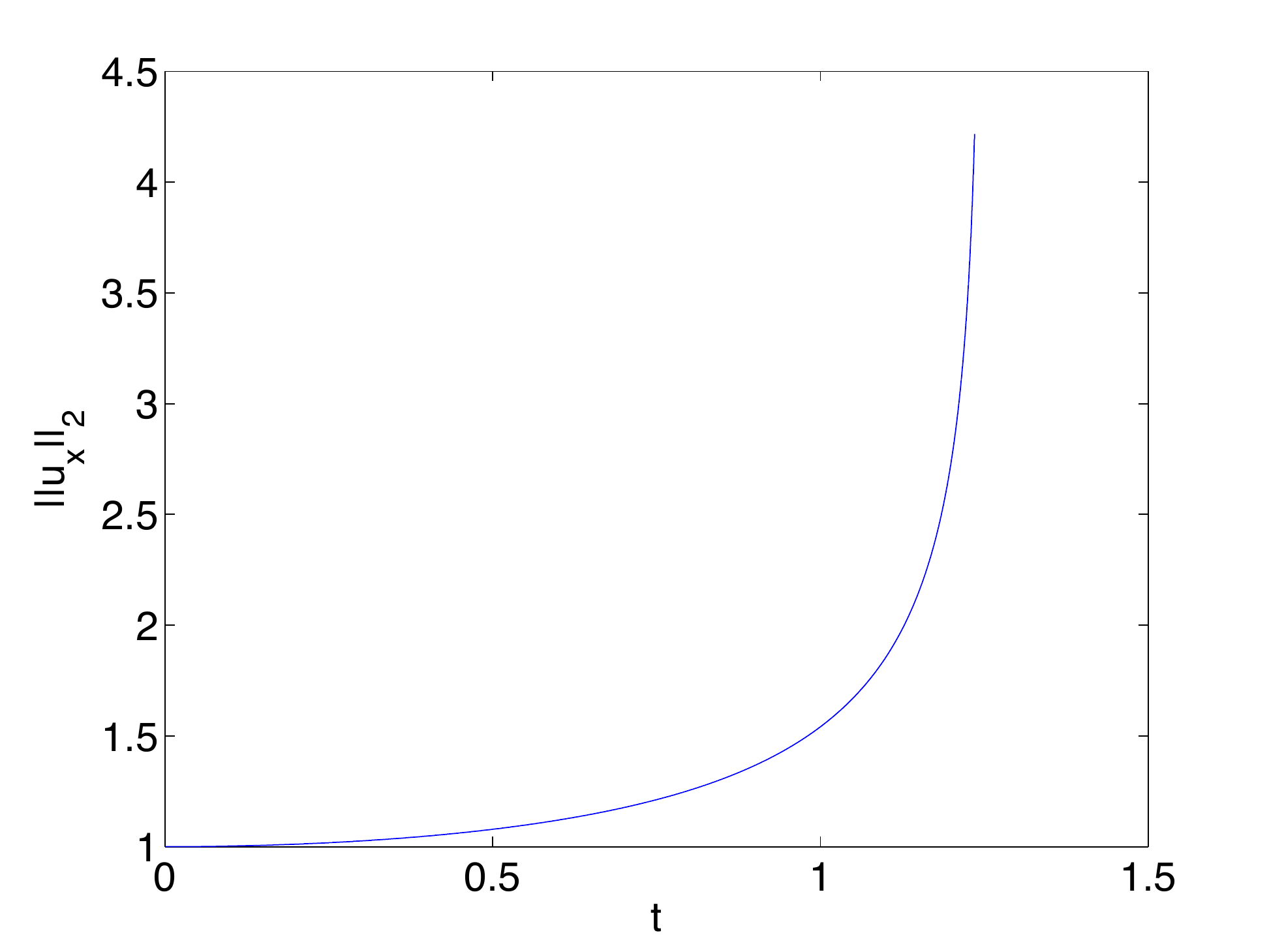}
 \caption{$L\infty$ norm of the Whitham solution of 
 Fig.~\ref{sautwhitham} in dependence of time on the left, and the 
 corresponding
 $L_{2}$ norm of $u_{x}$ on the right.}
 \label{sautwhithamnormp}
\end{figure}

The type of the singularity is once more confirmed by the norms of the 
solution to the Whitham equation (the corresponding norms for the 
fKdV solution are very similar and thus not shown) in 
Fig.~\ref{sautwhithamnormp}. The $L_{\infty}$ norm of the solution is 
monotonically increasing, the $L_{2}$ norm of the gradient has a 
blow-up. But as the Fourier coefficients indicate, this is not a 
cubic but a square root singularity in this case. 

The results of this subsection can be summarized in the following
\begin{conj}
    Consider  smooth initial data $u_{0}\in L_{2}(\mathbb{R})$ 
with a single negative hump. Then 
    \begin{itemize}
\item solutions to the Whitham equation (\ref{Whit}) and to fKdV 
equations with $-1<\alpha<0$ for initial 
data $u_{0}$ of sufficiently small mass stay smooth for 
all $t$ and will be radiated away.

\item solutions to the Whitham equation (\ref{Whit}) and to the fKdV 
equation with $\alpha=-1/2$ for negative initial 
data $u_{0}$ of sufficiently large mass will develop a cusp at 
$t^{*}>t_{c}$ 
of the form $|x-x^{*}|^{1/3}$. The sup norm of the solution remains bounded at the blow-up point.

\item solutions to the Whitham equation (\ref{Whit}) and to the fKdV 
equation with $\alpha=-1/2$ for positive initial 
data $u_{0}$ of sufficiently large norm mass will develop a cusp at $t^{*}<t_{c}$ 
of the form $|x-x^{*}|^{1/2}$.
 \end{itemize}
\end{conj}

\subsection{Large time behavior}
In this subsection we study the long time behavior of solutions to 
the fKdV (\ref{Cauchybis}) and fBBM equation (\ref{BBMbis}) in the 
case of nonlinearity and dispersion of order 
$\epsilon$ for the initial data $u_{0}=\beta\mbox{sech}^{2}x$. 

Introducing the slow time variable $\tau=t/\epsilon,$ the fKdV case reduces to the standard fKdV equation

\begin{equation}\label{sanseps}
u_\tau+uu_x-D^\alpha u_x=0,
\end{equation}

which  implies the dichotomy: either the solution is global or the blow-up time scales like the 
break-up time as $t\sim1/\epsilon$. The numerical study here will be therefore performed only as a test of the numerical approach.

Alternatively, introducing the function $\tilde{u}=\epsilon u$, we get the equations
\begin{equation}
    \tilde{u}_{t}+\tilde{u}\tilde{u}_{x}-\epsilon 
    D^{\alpha}\tilde{u}_{x}=0,\quad \tilde{u}_0=\epsilon u_0
    \label{fKdVt}
\end{equation}
and
\begin{equation}
    \tilde{u}_{t}+\tilde{u}_{x}+\tilde{u}\tilde{u}_{x}+
    \epsilon D^{\alpha}\tilde{u}_{t}=0, \quad \tilde{u}_0=\epsilon u_0.
    \label{fBBMt}
\end{equation}

Both (\ref{fKdVt}) and (\ref{fBBMt}) 
will be solved for the initial data $\tilde{u}_{0}=\epsilon \beta
\mbox{sech}^{2}x$ for several values of $\epsilon$. For $\epsilon=0$, 
both equations reduce to the Burgers equation. For 
the initial data $\tilde{u}_{0}$, the solutions to the Burgers equation 
have a point of gradient catastrophe, a hyperbolic blow-up, at the 
critical time 
\begin{equation}
    t_{c}=\frac{3^{3/2}}{4\beta\epsilon}
    \label{tch}.
\end{equation}
Obviously the critical time $t_{c}$ is of order $1/\epsilon$. We will 
now study for several $\epsilon$ for examples, for which we 
observed blow-up in the previous sections, how the  time 
$t^{*}$ scales with $\epsilon$. 

We first consider fKdV with $\alpha=0.2$ for the initial data 
$\epsilon \mbox{sech}^{2}x$ for the values 
$\epsilon=0.01,0.02,\ldots,0.1$. Since the blow-up is exponential, we 
just take the time for which the code stops converging as the blow-up 
time. This gives a very good approximation to $t^{*}$, and what is 
important since we are interested in the $\epsilon$ dependence, a 
consistent one for all $\epsilon$. The computation is carried out 
with $N=2^{16}$ Fourier modes for $x\in[-20\pi,20\pi]$ and 
$N_{t}=10^{4}$ time steps. Doing a linear regression analysis for 
$\log_{10}t^{*}\sim a\log_{10}\epsilon+b$, we find that 
the blow-up time $t^{*}$ is as the time $t_{c}$ (\ref{tch}) of order 
$1/\epsilon$. More precisely we get $a=-0.9998$, $b=0.4845$ with 
standard deviation $\sigma_{a}=2.1*10^{-4}$ and
correlation coefficient $r = 9*10^{-8}$. Thus the expected dependence 
$t^{*}\propto 1/\epsilon$ is numerically well observed. 

A corresponding analysis for fBBM is much more involved for two 
reasons: first the initial peak travels to the right before blowing 
up, the farther the smaller $\epsilon$, whereas blow-up for fKdV 
happens close to the initial maximum. And secondly the blow-up is not 
as pronounced in fBBM as in fKdV for $\alpha<0.5$. Thus it is less 
obvious to obtain an estimate for the blow-up time than for fKdV. The 
first problem is addressed as in (\ref{resc}) and (\ref{vx}) by using 
a frame commoving with the location of the maximum $x_{m}(t)$. Thus 
we solve
$$U_{t}-VU_{y}+(1+\epsilon D^{\alpha})^{-1}(U_{y}+UU_{y})=0,$$
where 
$$V = \left.\frac{(1+\epsilon 
D^{\alpha})^{-1}(U_{y}+UU_{y})}{U_{yy}}\right|_{y=0}.$$
To determine the blow-up time, we use the approach \cite{SSF,KR2013} 
based on a fitting of the asymptotic behavior of the Fourier 
coefficients.  We again use the 
initial data $U_{0}=\epsilon \mbox{sech}^{2}x$ for the values 
$\epsilon=0.01,0.02,\ldots,0.1$.  The computation is carried out 
with $N=2^{15}$ Fourier modes for $x\in[-10\pi,10\pi]$ and 
$N_{t}=2*10^{4}$ time steps. Doing a linear regression as for fKdV, we find that 
the blow-up time $t^{*}$ is again as the time $t_{c}$ of order 
$1/\epsilon$. More precisely we get $a=-0.9677$, $b=    0.5521$ with 
standard deviation $\sigma_{a}=0.005$ and
correlation coefficient $r = 0.9999$. Thus we find again that the 
blow-up time $t^{*}$ is proportional to $1/\epsilon$ as the break-up  
time $t_{c}$ for the corresponding Burgers solution. Note that we always 
have $t^{*}>t_{c}$.

The above results therefore illustrates the fact   that the life span of solutions 
of the Burgers equation 
\begin{equation}\label{Burgeps}
u_t+\epsilon uu_x=0
\end{equation}
is not enhanced by a weak dispersive 
perturbation having the same small coefficient $\epsilon$. On the other hand, as 
we have seen in the previous subsection the blow-up is of a different 
nature than the 
shock type one of the Burgers equation.   

When the small parameter $\epsilon$ affects only the quadratic term 
in (\ref{sanseps}), we have already noticed that 
the resulting equation reduces to the standard fKdV equation with initial data of order $\epsilon$ and our simulations suggest that
for any $0<\alpha<1$ the solution is global for $\epsilon$ small enough.


\section{Outlook}
The numerical results in this paper suggest that blow-up in 
solutions to fKdV 
equations when $1/3<\alpha\leq 1/2$ is similar to blow-up in solutions to gKdV equations. In 
the $L^{2}$ critical case ($\alpha =1/2$), the blow-up profile appears to be given by 
a rescaled soliton. In the $L^2$ supercritical case $1/3<\alpha<1/2$, the blow-up profile 
should be given by an asymptotically decreasing solution to the 
fractionary equation (\ref{sautalpharinfty}). To obtain further insight 
into this case, the asymptotics of solutions to this equation have to 
be worked out. It has to be shown whether there is a unique solution 
which is asymptotically decreasing. Then a numerical scheme has to 
be developed to compare the blow-up in fKdV solutions to a rescaled 
form of this particular solution to (\ref{sautalpharinfty}). This 
will be subject of further work. 

Though the blow-up in the energy supercritical case 
$0<\alpha<1/3$ appears to be given by the same asymptotic profile, it 
is of a different nature due to the absence of solitary waves.

Another important question is related to the long-time behavior of 
solutions to fBBM equations for small nonlinearity as studied in the 
previous subsection. It would be worthwhile to investigate whether or 
not the found results persist in a more general context of dispersive 
perturbations of quasilinear hyperbolic systems and of course for 
relevant water waves models such as the Boussinesq systems. In the 
latter case things might be more subtle since the presence of a ``BBM" 
term in the system weakens the nonlinearity in the corresponding 
equation. 

\begin{merci}
CK thanks for financial support by the ANR via the program 
ANR-09-BLAN-0117-01, JCS acknowledges partial support from the ANR project GEODISP.
\end{merci}

\bibliographystyle{amsplain}

\end{document}